\documentclass[12pt, reqno]{amsart}

\usepackage{
amsmath,
amssymb,
amsthm,
mathrsfs,
amsfonts,
ytableau,
enumitem,
comment,
upgreek,
soul,
yhmath,
mathtools,
shuffle,
kotex,
array,
caption,
tabularx
}
\usepackage{makecell}
\usepackage{multirow}

\usepackage{aliascnt}
\usepackage[linktocpage]{hyperref}

\hypersetup{
    colorlinks=true,
    linkcolor = blue,
    citecolor=magenta,
    urlcolor=cyan,
}

\usepackage[capitalise, nameinlink]{cleveref}

\crefname{equation}{}{}
\crefname{figure}{{\sc Figure}}{{\sc Figure}}
\crefname{subsection}{Section}{Sections}

\usepackage[euler]{textgreek}

\usepackage{tikz,tikz-cd}

\usepackage[colorinlistoftodos, textwidth = 2.3cm]{todonotes}

\ytableausetup{mathmode, boxsize=1.2em}

\theoremstyle{plain}

\newtheorem{theorem}{Theorem}[section]

\newaliascnt{proposition}{theorem}
\newtheorem{proposition}[proposition]{Proposition}
\aliascntresetthe{proposition}

\newaliascnt{lemma}{theorem}
\newtheorem{lemma}[lemma]{Lemma}
\aliascntresetthe{lemma}

\newaliascnt{corollary}{theorem}
\newtheorem{corollary}[corollary]{Corollary}
\aliascntresetthe{corollary}

\newaliascnt{conjecture}{theorem}

\aliascntresetthe{conjecture}

\newaliascnt{claim}{theorem}

\aliascntresetthe{claim}

\newaliascnt{question}{theorem}

\aliascntresetthe{question}

\newtheorem*{claim*}{Claim}

\theoremstyle{definition}

\newaliascnt{algorithm}{theorem}

\aliascntresetthe{algorithm}

\newaliascnt{example}{theorem}
\newtheorem{example}[example]{Example}
\aliascntresetthe{example}

\newaliascnt{definition}{theorem}
\newtheorem{definition}[definition]{Definition}
\aliascntresetthe{definition}

\newaliascnt{remark}{theorem}
\newtheorem{remark}[remark]{Remark}
\aliascntresetthe{remark}

\crefname{theorem}{Theorem}{Theorems}
\crefname{proposition}{Proposition}{Propositions}
\crefname{lemma}{Lemma}{Lemmas}
\crefname{corollary}{Corollary}{Corollaries}
\crefname{conjecture}{Conjecture}{Conjectures}
\crefname{claim}{Claim}{Claims}
\crefname{question}{Question}{Questions}
\crefname{algorithm}{Algorithm}{Algorithms}
\crefname{example}{Example}{Examples}
\crefname{definition}{Definition}{Definitions}
\crefname{remark}{Remark}{Remarks}

\numberwithin{equation}{section}
\numberwithin{figure}{section}
\numberwithin{table}{section}

\def\Z{\mathbb Z}
\def\C{\mathbb C}

\newcommand{\nc}{\newcommand}

\nc\doublequotient[3]{{_{\displaystyle #1}}\backslash{^{\displaystyle #2}}/{_{\displaystyle #3}}}
\nc\rightquotient[2]{{_{\displaystyle #1}}\backslash{^{\displaystyle #2}}}
\nc\leftquotient[2]{{^{\displaystyle #1}}/{_{\displaystyle #2}}}
\nc{\boxleq}{\leq_{\lambda/\mu}}
\nc{\boxl}{<_{\lambda/\mu}}
\nc{\restarrow}{\hspace{-0.5ex}\downarrow}

\nc{\rmev}{\mathrm{ev}}
\nc{\tilS}{\widetilde{S}}
\nc{\tilal}{\widetilde{\alpha}}
\nc{\tilbe}{\widetilde{\beta}}
\nc{\tilchi}{\widetilde{\chi}}
\nc{\barchi}{\chi'}
\nc{\rmCh}{\mathrm{Ch}}
\nc{\barrmCh}[2]{{\rmCh(#1, #2)}'}

\nc{\sfw}{\mathsf{w}}
\nc{\sfpar}{\mathsf{par}}
\nc{\ttg}{\mathtt{g}}
\nc{\sfcomp}{\mathsf{comp}}
\nc{\sfu}{\mathsf{u}}
\nc{\sfcyc}{\mathsf{cyc}}
\nc{\sfremov}{\mathsf{remov}}
\nc{\sfP}{\mathsf{P}}
\nc{\sfS}{\mathsf{S}}
\nc{\sfh}{\mathsf{h}}

\nc{\sfDS}{\mathsf{DS}}
\nc{\sfst}{\mathsf{st}}
\nc{\tDC}{\mathtt{DC}}
\nc{\scrC}{\mathscr{C}}
\nc{\scrD}{\mathscr{D}}
\nc{\scrL}{\mathscr{L}}
\nc{\SR}{\mathcal{SR}}
\nc{\sfposet}
{\mathsf{poset}}
\nc{\sftr}{\mathsf{tr}}
\nc{\rmcard}{\mathrm{card}}
\nc{\rmsor}{\mathrm{source}}
\nc{\ocard}{\overline{\mathrm{card}}}
\nc{\omin}{\overline{\mathrm{min}}}
\nc{\omax}{\overline{\mathrm{max}}}
\nc{\calL}{\mathcal{L}}
\nc{\calC}{\mathcal{C}}
\nc{\calK}{\mathcal{K}}
\nc{\calI}{\mathcal{I}}
\nc{\sfM}{\mathsf{M}}
\nc{\sfg}{\mathsf{g}}
\nc{\sfm}{\mathsf{m}}
\nc{\Ptab}{\mathtt{ins}}
\nc{\Qtab}{\mathtt{rec}}
\nc{\pt}{T^{\scalebox{0.5}{$\nearrow$}}}
\nc{\rmIm}{\mathrm{Im}}
\nc{\rmInv}{\mathrm{Inv}}
\nc{\rmInt}{\mathrm{Int}}
\nc{\Deq}{\overset{D}{\simeq}}
\nc{\Dleq}{\preceq_D}
\nc{\bdP}{\boldsymbol{P}}
\nc{\bdI}{\boldsymbol{I}}
\nc{\soc}{\mathrm{soc}}
\nc{\balproj}{\bal_{\mathrm{proj}}}
\nc{\balinj}{\bal_{\mathrm{inj}}}

\nc{\Larrow}{\leftarrow_\mathscr{L}}
\nc{\Lorder}{\le_\mathscr{L}}
\nc{\Lsim}{\sim_\mathscr{L}}

\nc{\rmfree}{\mathrm{free}}
\nc{\rmtor}{\mathrm{tor}}

\nc{\SG}{\mathfrak{S}}
\nc{\frakf}{\mathfrak{f}}
\nc{\frakR}{\mathfrak{R}}
\nc{\frakL}{\mathfrak{L}}
\nc{\frakJ}{\mathfrak{J}}
\nc{\frakC}{\mathfrak{C}}
\nc{\frakM}{\mathfrak{M}}

\nc{\PCT}{\mathrm{PCT}}
\nc{\SPCT}{\mathrm{SPCT}}
\nc{\RT}{\mathrm{RT}}
\nc{\SRT}{\mathrm{SRT}}
\nc{\SYT}{\mathrm{SYT}}
\nc{\SSYT}{\mathrm{SSYT}}
\nc{\LR}{\mathrm{LR}}
\nc{\sgn}{\mathrm{sgn}}
\nc{\RCT}{\mathrm{RCT}}
\nc{\SRCT}{\mathrm{SRCT}}
\nc{\SYRT}{\mathrm{SYRT}}
\nc{\SYCT}{\mathrm{SYCT}}
\nc{\SPYCT}{\mathrm{SPYCT}}
\nc{\tst}{\mathtt{st}}
\nc{\ttS}{\mathtt{S}}
\nc{\ttT}{\mathtt{T}}
\nc{\ttU}{\mathtt{U}}
\nc{\ttQ}{\mathtt{Q}}
\nc{\Span}{\mathrm{span}}
\nc{\comp}{\mathrm{comp}}
\nc{\rmst}{\mathrm{st}}
\nc{\bfDes}{\mathbf{Des}}
\nc{\Des}{\mathrm{Des}}
\nc{\set}{\mathrm{set}}
\nc{\wt}{\mathrm{wt}}
\nc{\ch}{\mathrm{ch}}
\nc{\id}{\mathrm{id}}
\nc{\Sym}{\mathrm{Sym}}
\nc{\Qsym}{\mathrm{QSym}}
\nc{\Nsym}{\mathrm{NSym}}
\nc{\sh}{\mathrm{sh}}
\nc{\bfC}{\mathbf{C}}
\nc{\bfD}{\mathbf{D}}
\nc{\bfS}{\mathbf{S}}
\nc{\bfm}{\mathbf{m}}
\nc{\hbfS}{\widehat{\mathbf{S}}}
\nc{\bfF}{\mathbf{F}}

\nc{\calA}{\mathcal{A}}
\nc{\calB}{\mathcal{B}}
\nc{\calG}{\mathcal{G}}
\nc{\calR}{\mathcal{R}}
\nc{\calS}{\mathcal{S}}
\nc{\calV}{\mathcal{V}}

\nc{\sfR}{\mathsf{R}}
\nc{\wcol}{w_{\mathrm{col}}}
\nc{\wrow}{w_{\mathrm{row}}}

\nc{\bfz}{\mathbf{z}}
\nc{\bfv}{\mathbf{z}}
\nc{\bfx}{\mathbf{x}}
\nc{\bfy}{\mathbf{y}}
\nc{\tal}{\lambda(\alpha)}
\nc{\tbe}{\widetilde{\beta}}
\nc{\opi}{\overline{\pi}}
\nc{\calP}{\mathcal{P}}
\nc{\rmtop}{\mathrm{top}}
\nc{\rad}{\mathrm{rad}}
\nc{\bfP}{\mathbf{P}}
\nc{\bfw}{\mathbf{w}}
\nc{\SET}{\mathrm{SET}}

\nc{\tcd}{\mathtt{cd}}
\nc{\tyd}{\mathtt{yd}}
\nc{\trd}{\mathtt{rd}}

\nc{\rmr}{\mathrm{r}}
\nc{\rmc}{\mathrm{c}}
\nc{\rmt}{\mathrm{t}}

\nc{\len}{\mathsf{len}}
\nc{\calF}{\mathcal{F}}
\nc{\calE}{\mathcal{E}}
\nc{\calT}{\mathscr{T}}
\nc{\sfT}{\mathsf{T}}
\nc{\calEsa}{\mathcal{E}^\upsig(\alpha)}
\nc{\tauC}{\tau_{\scalebox{0.5}{$C$}}}
\nc{\sytabC}{\sytab_{\scalebox{0.5}{$C$}}}
\nc{\pr}{\mathsf{pr}}
\nc{\Ups}{\Upsilon}
\nc{\pact}{\diamond}
\nc{\tauE}{\tau_{E}^{~}}
\nc{\tauF}{\tau_{\scalebox{0.5}{$F$}}}
\nc{\tauG}{\tau_{\scalebox{0.5}{$G$}}}
\nc{\rtE}{T_{\scalebox{0.5}{$E$}}}
\nc{\rtF}{T_{\scalebox{0.5}{$F$}}}
\nc{\rtG}{T_{\scalebox{0.5}{$G$}}}
\nc{\sytab}{\widehat{\tau}}
\nc{\hatE}{\widehat{E}}
\nc{\hJ}{\widehat{J}}
\nc{\hati}{\hat{i}}
\nc{\hcalE}{\widehat{\calE}}
\nc{\hatC}{\widehat{C}}
\nc{\bal}{{\boldsymbol{\upalpha}}}
\nc{\bbe}{{\boldsymbol{\upbeta}}}
\nc{\SPYRT}{\mathrm{SPYRT}}

\nc{\hfrakJ}{\widehat{\mathfrak{J}}}
\nc{\bfJ}{\boldsymbol{J}}
\nc{\hbfJ}{\boldsymbol{\widehat{J}}}
\nc{\bfM}{\boldsymbol{M}}

\nc{\bgam}{{\boldsymbol{\upgamma}}}
\nc{\bdel}{{\boldsymbol{\updelta}}}
\nc{\weakcon}{\odot}
\nc{\basisI}{I}

\newcommand{\tf}{\widetilde{f}}

\nc{\ldalpha}{\lambda(\alpha)}
\nc{\SRIT}{\mathrm{SRIT}}
\nc{\re}{\mathrm{rev}}
\nc{\otau}{\overline{\tau}}
\nc{\rtop}{{\rm top}}
\nc{\sfc}{\mathsf{c}}
\nc{\sfchar}{\mathsf{char}}
\nc{\sfcharbar}{\mathsf{char}'}
\nc{\sfr}{\mathsf{r}}
\nc{\sfSP}{\mathsf{SP}}
\nc{\sfRP}{\mathsf{RP}}
\nc{\sfRSP}{\mathsf{RSP}}
\nc{\tH}{\mathtt{H}}
\nc{\tV}{\mathtt{V}}
\nc{\rpi}{\mathring{\pi}}
\nc{\cpi}{\check{\pi}}
\nc{\frakm}{\mathfrak{m}}
\nc{\sfem}{\mathsf{em}}
\nc{\Hom}{\mathrm{Hom}}
\nc{\module}{\mathrm{mod}}
\nc{\rmread}{\mathsf{read}}
\nc{\tread}{\underline{\mathsf{read}}}
\nc{\ocalE}{\overline{\calE}}
\nc{\oE}{\overline{E}}

\nc{\SPCTsa}{\SPCT^\upsig(\alpha)}
\nc{\bfSsa}{\bfS_\alpha^\upsig}
\nc{\bfSsaC}{{\bfS}^\upsig_{\alpha,C}}
\nc{\hbfSsa}{\widehat{\bfS}_\alpha^\upsig}
\nc{\upineq}{\rotatebox{90}{$<$}}
\nc{\downineq}{\rotatebox{270}{$<$}}
\nc{\diagineq}{\rotatebox{135}{$<$}}
\nc{\sfB}{\mathsf{B}}
\nc{\hxi}{\widehat{\xi}}
\nc{\hxidwJ}{\hxi_{\scalebox{0.55}{$J$}}}
\nc{\hxiupJ}{\hxi^{\scalebox{0.55}{$J$}}}
\nc{\scrS}{\mathscr{S}}
\nc{\scrB}{\mathscr{B}}
\nc{\bfT}{\mathbf{T}}
\nc{\tshuffle}{\,\widetilde{\shuffle}\,}
\nc{\sJ}{\scalebox{0.55}{$J$}}
\nc{\sJo}{\scalebox{0.55}{$J_1$}}
\nc{\sJt}{\scalebox{0.55}{$J_2$}}
\nc{\Keq}{\overset{K}{\cong}}
\nc{\dKeq}{\overset{K^*}{\cong}}
\nc{\Rect}{\mathrm{Rect}}

\nc{\ra}{\rightarrow}

\nc{\matr}[2]{\left( \hspace{-1ex} \begin{array}{c} #1 \\ #2 \end{array} \hspace{-1ex} \right)}

\definecolor{wsgreen}{rgb}{0,0.5,0}

\nc{\DIRT}{\mathrm{DIRT}}
\nc{\hpi}{\pi}
\nc{\frakI}{\mathfrak{I}}
\nc{\hfrakI}{\widehat{\mathfrak{I}}}
\nc{\orho}{\overline{\rho}}
\nc{\autotheta}{\uptheta}
\nc{\calW}{\mathcal{W}}
\nc{\Endo}{\mathrm{End}}

\nc{\autophi}{\upphi}
\nc{\autochi}{\upchi}
\nc{\autoomega}{\upomega}
\nc{\hIM}{\widehat{\sfB}}
\nc{\bfpi}{\boldsymbol{\uppi}}
\nc{\bfopi}{\overline{\boldsymbol{\uppi}}}
\nc{\osfB}{\overline{\sfB}}
\nc{\rmw}{\mathrm{w}}
\nc{\conc}{\; {\bullet} \;}
\nc{\ostar}{\; \overline{\bullet} \;}
\nc{\rank}{\mathrm{rank}}
\nc{\fkp}{\mathfrak{p}}
\nc{\bfR}{\mathbf{R}}
\nc{\calD}{\mathcal{D}}

\nc{\upsig}{{\boldsymbol{\upsigma}}}
\nc{\bfSsaE}{{\bfS}^\upsig_{\alpha,E}}

\nc{\hfkp}{\widehat{\mathfrak{p}}}

\nc{\hautophi}{{\widehat{\autophi}}}
\nc{\hautotheta}{{\widehat{\autotheta}}}
\nc{\hautoomega}{{\widehat{\autoomega}}}
\nc{\rmperm}{\mathrm{perm}}
\nc{\Hnmod}{\text{$H_n(0)$-$\mathbf{mod}$}}
\nc{\modHn}{\text{$\mathbf{mod}$-$H_n(0)$}}

\nc{\laurH}{H_{\C[v^{\pm 1/2}];v}(\SG_n)}
\nc{\polyH}{H_{\C[v^{1/2}];v}(\SG_n)}
\nc{\genericH}{H_{\C(v^{1/2});v}(\SG_n)}
\nc{\compH}[1]{H_{\C;#1}(\SG_n)}

\nc{\bfsigJ}{\sigma_{\scalebox{0.55}{$J$}}}
\nc{\bfrhoJ}{\rho^{\scalebox{0.55}{$J$}}}

\nc{\teta}{\widetilde{\eta}}
\nc{\urmw}{\underline{\mathrm{w}}}
\nc{\sfcnt}{\mathsf{cnt}}

\nc{\pistar}[1]{\pi_{#1}^*}
\nc{\wfkp}{\widetilde{\mathfrak{p}}}
\nc{\bfpsi}{\boldsymbol{\uppsi}}

\nc{\boxn}{\setlength{\fboxsep}{1.5pt} \scriptsize \boxed{n}}

\nc{\yt}[1]{\todo[size=\tiny,color=blue!10]{#1 \\ \hfill --- Young-Tak}}
\nc{\YT}[1]{\todo[size=\tiny,inline,color=blue!10]{#1
		\\ \hfill --- Young-Tak}}
		
\nc{\yh}[1]{\todo[size=\tiny,color=cyan!10]{#1 \\ \hfill --- Young-Hun}}
\nc{\YH}[1]{\todo[size=\tiny,inline,color=cyan!10]{#1
		\\ \hfill --- Young-Hun}}
		
\nc{\sy}[1]{\todo[size=\tiny,color=magenta!10]{#1 \\ \hfill --- So-Yeon}}
\nc{\SY}[1]{\todo[size=\tiny,inline,color=magenta!10]{#1
		\\ \hfill --- So-Yeon}}

\nc{\nt}[1]{\todo[size=\tiny,color=green!10]{#1 \\ \hfill --- Note}}
\nc{\NT}[1]{\todo[size=\tiny,inline,color=green!10]{#1
		\\ \hfill --- Note}}

\definecolor{purple}{rgb}{0.44, 0.0, 1.0}

\definecolor{yhblue}{rgb}{0,0,0.6}

\newenvironment{red}{\relax\color{red}}{\hspace*{.5ex}\relax}
\newenvironment{blue}{\relax\color{blue}}{\hspace*{.5ex}\relax}
\newenvironment{green}{\relax\color{wsgreen}}{\hspace*{.5ex}\relax}
\newenvironment{magenta}{\relax\color{magenta}}{\hspace*{.5ex}\relax}
\newenvironment{purple}{\relax\color{purple}}{\hspace*{.5ex}\relax}

\nc{\ber}{\begin{red}}
\nc{\er}{\end{red}}
\nc{\beb}{\begin{blue}}
\nc{\eb}{\end{blue}}
\nc{\bema}{\begin{magenta}}
\nc{\ema}{\end{magenta}}
\nc{\begr}{\begin{green}}
\nc{\egr}{\end{green}}

\nc{\bepu}{\begin{purple}}
\nc{\epu}{\end{purple}}

\nc{\inn}[1]{\left\langle \, #1 \, \right\rangle}

\title[Relationship between generic Hecke algebra modules and $0$-Hecke modules]{Relationship between generic Hecke algebra modules and $0$-Hecke algebra modules}

\author[Y.-H. Kim]{Young-Hun Kim}
\address[Y.-H. Kim]{Department of Mathematics, Seoul Women’s University, Seoul 01797, Republic of Korea}
\email{yhkim@swu.ac.kr, ykim.math@gmail.com}

\author[S.-Y. Lee]{So-Yeon Lee}
\address[S.-Y. Lee]{Department of Mathematics, Seoul National University, Seoul 08826, Republic of Korea}
\email{sylsyl0814@gmail.com}

\author[Y.-T. Oh]{Young-Tak Oh}
\address[Y.-T. Oh]{Department of Mathematics / Institute for Mathematical and Data Sciences, Sogang University, Seoul 04107, Republic of Korea}
\email{ytoh@sogang.ac.kr}

\keywords{Hecke algebra,
$0$-Hecke algebra, Specialization functor,
skew Specht module,
Kazhdan--Lusztig cell module,
dual plactic class module}

\date{\today}
\subjclass[2020]{20C08, 05E10, 05E05}
\begin{document}

\begin{abstract}
Let $H_{R,a}(\mathfrak{S}_n)$ be the Hecke algebra over a principal ideal domain $R$ with parameter $a\in R$.
The representation theories of 
$H_{\mathbb{C}(v^{1/2});v}(\mathfrak{S}_n)$ and $H_{\mathbb{C};0}(\mathfrak{S}_n)$ differ substantially and have largely been studied independently.
In this paper, we relate the categories of finitely generated modules of these two algebras through specialization functors from the category of finitely generated $H_{\mathbb{C}[v^{1/2}];v}(\mathfrak{S}_n)$-modules.
We first investigate the fundamental properties of these functors concerning dimensions, essential surjectivity, and characters.
We then construct $H_{\mathbb{C}[v^{1/2}];v}(\mathfrak{S}_n)$-modules which, through these functors, provide connections between skew Specht modules of $H_{\mathbb{C}(v^{1/2});v}(\mathfrak{S}_n)$ and $0$-Hecke skew Specht modules of $H_{\mathbb{C};0}(\mathfrak{S}_n)$, and between Kazhdan--Lusztig cell modules of $H_{\mathbb{C}(v^{1/2});v}(\mathfrak{S}_n)$ and dual plactic class modules of $H_{\mathbb{C};0}(\mathfrak{S}_n)$.
Both families of $H_{\mathbb{C};0}(\mathfrak{S}_n)$-modules are closely related to regular Schur labeled skew shape posets.
We further explore how specialization functors relate the Frobenius characteristic to the quasisymmetric characteristic and connects the branching rules for skew Specht modules with those for $0$-Hecke skew Specht modules.
Finally, we use the established connection with Kazhdan--Lusztig cell modules to obtain a  result on the classification of dual plactic class modules up to isomorphism.
\end{abstract}

\maketitle

\tableofcontents

\section{Introduction}

Iwahori--Hecke algebras (Hecke algebras)  are  deformations of the group algebras of Coxeter groups and occupy a central place in representation theory.
Since the late 1970s, there have been 
many remarkable results on  representations of Hecke algebras (for instance, see~\cite{79KL, 86DJ, 96Ariki, 96GL, 06Geck, 09BK}).
In particular, the theory of the Kazhdan--Lusztig basis and the associated cell modules introduced in \cite{79KL} has been employed in revealing  beautiful connections between Hecke algebras and other topics such as Schubert varieties,  semisimple Lie algebras, and  categorifications  (see \cite{80KL,81BK,00Ariki,14EBW}).

For a principal ideal domain $R$ and $a\in R$, let $H_{R;a}(\SG_n)$ denote the Hecke algebra of the symmetric group $\SG_n$ over $R$ with parameter $a$ (see \cref{eq: def of Hecke alg} for its precise definition).
Even over $\mathbb C$, the representation theory of these algebras can vary significantly with the parameter.
Let $v^{1/2}$ be an indeterminate, and call $\genericH$ the \emph{generic Hecke algebra}.
For nonzero $q\in\mathbb C$, $\genericH$ and $\compH{q}$ share a cellular framework that allows their irreducible modules to be constructed and classified uniformly.
This common framework remains applicable when $q$ is a root of unity although the representation theory of $\compH{q}$ may then differ significantly from the generic case.
By contrast, the parameter of the \emph{$0$-Hecke algebra} $\compH{0}$ is not invertible, and thus the same cellular framework does not apply. 
The representation theory of $\compH{0}$ has largely developed independently of that of the generic Hecke algebra.
This contrast motivates us to investigate how these two representation theories can nevertheless be connected in a suitable way.

We now briefly review the representation theories of the generic Hecke algebra and the $0$-Hecke algebra.

The generic Hecke algebra $\genericH$ is isomorphic to the group algebra $\C(v^{1/2})[\SG_n]$.
Thus it is semisimple and its irreducible modules are parametrized by partitions of $n$.
We first recall two notable constructions of irreducible $\genericH$-modules.
One is due to  Kazhdan and Lusztig~\cite{79KL}.
Using a bar-invariant basis,  they constructed the Kazhdan--Lusztig cell modules
$\frakM_Q(\C(v^{1/2});v)$, where $Q$ ranges over the set $\SYT_n$ of standard Young tableaux of size $n$.
These modules form a complete list of  irreducible $\genericH$-modules.
For each $Q$,  the module $\frakM_Q(\C(v^{1/2});v)$ has a basis indexed by the dual plactic class $\Gamma_Q$ of $Q$, namely the set of permutations whose recording tableau under the Robinson--Schensted correspondence is $Q$.
The other is due to Dipper and James \cite{86DJ}.
Using an analogue of Young symmetrizers,  they constructed the Specht modules $S^{\lambda}(\C(v^{1/2});v)$, where $\lambda$ ranges over partitions of $n$.
These Specht modules form a complete list of pairwise nonisomorphic irreducible $\genericH$-modules.
For each $\lambda$, the module $S^{\lambda}(\C(v^{1/2});v)$ has a basis indexed by standard Young tableaux of shape $\lambda$.
We next recall the \emph{Frobenius characteristic}, denoted by $\ch^{\rm Frob}$.
It is an isomorphism between the Grothendieck ring of the tower of generic Hecke algebras $\genericH$ ($n \geq 0$) and the ring of symmetric functions.
It sends the isomorphism classes of irreducible modules to Schur functions (see \cref{subsubsec: Frob char}).

The $0$-Hecke algebra $\compH{0}$ is neither semisimple nor of finite representation type for $n>3$~\cite{11DY,02DHT}. 
Norton~\cite{79Norton} classified the irreducible modules and projective indecomposable modules of $\compH{0}$ up to isomorphism.
These modules are parametrized by compositions  of $n$. 
Let us introduce two families of $0$-Hecke modules in our concern.
To begin with, we recall the \emph{quasisymmetric characteristic}, denoted by $\ch$, introduced by Duchamp, Krob, Leclerc, and Thibon \cite{96DKLT}.
It is an isomorphism between the Grothendieck ring of the tower of $0$-Hecke algebras $\compH{0}$ ($n \geq 0$) and the ring of quasisymmetric functions.
It sends the isomorphism classes of irreducible modules to fundamental quasisymmetric functions.
Using this isomorphism, Duchamp, Hivert, and Thibon~\cite{02DHT} provided a representation theoretical interpretation of $P$-partition generating functions.
To each labeled poset $P$, they associated a $0$-Hecke module $M_P$ whose quasisymmetric characteristic is the $P$-partition generating function.
These modules have been dealt with in several papers (see \cite{03McNamara, 23CKO, 23KLO, 25Searles, 26KS}).
In particular, in \cite{23KLO}, the modules associated with regular Schur labeled skew shape posets were intensively studied in the context of Stanley's $P$-partitions conjecture.
This conjecture states that every labeled poset with a symmetric $P$-partition generating function is a Schur labeled skew shape poset.
In that paper, the authors classified the $\compH{0}$-modules associated with regular Schur labeled skew shape posets up to isomorphism and introduced pairwise nonisomorphic representatives $X_{\lambda/\mu}$, where $\lambda/\mu$ ranges over skew partitions of $n$.
These form the first family of modules considered here.
Each $X_{\lambda/\mu}$ has a basis indexed by standard Young tableaux of shape $\lambda/\mu$ and its images under $\ch$  is the skew Schur function $s_{\lambda/\mu}$.
It was also shown in~\cite{23KLO} that the modules $X_{\lambda/\mu}$ admit filtrations of a special type, called \emph{distinguished filtrations}.
To represent the successive quotients of these filtrations up to isomorphism, we introduce an $\compH{0}$-module $\scrD_Q$ for each $Q \in \SYT_n$.
These form the second family of modules considered here.
Each $\scrD_Q$ has $\Gamma_Q$ as a basis and
its images under $\ch$ is the Schur function $s_{\sh(Q)'}$, where $\sh(Q)'$ is the conjugate of the shape of $Q$.
 Hereafter, we call the modules $X_{\lambda/\mu}$ and $\scrD_Q$ the \emph{$0$-Hecke skew Specht module} attached to $\lambda/\mu$ and the \emph{dual plactic class module} attached to $Q$, respectively (see \cref{subsec: two kinds of notable 0 Hecke modules}).

The aim of this paper is to study the relationship between two families of $\genericH$-modules and those of $\compH{0}$-modules by means of suitable specialization functors.

In \cref{sec: specialization functor}, we  establish three fundamental properties of the specialization functors 
\begin{align*}
\calF_{\iota}&:\polyH\text{-}\module\to \genericH\text{-}\module\quad\text{and}\\
\calF_{\theta_q}&:\polyH\text{-}\module
\to H_{\C;\,q^2}(\SG_n)\text{-}\module \quad(q \in \C)
\end{align*}
induced by the natural embedding
$\iota:\C[v^{1/2}]\hookrightarrow\C(v^{1/2})$ and the evaluation homomorphism
$\theta_q:\C[v^{1/2}]\to\C$ at $v^{1/2} = q$, respectively.
Here, $H_{R;a}(\SG_n)\text{-}\module$ denotes the category of finitely generated $H_{R;a}(\SG_n)$-modules.
Let $M$ be a  finitely generated $\polyH$-module. 
We first determine the dimensions of $\calF_\iota(M)$ and $\calF_{\theta_q}(M)$.
More precisely, it holds that 
\[
\dim_{\C(v^{1/2})} \calF_\iota(M) = r_M
\quad \text{and} \quad
\dim_{\C} \calF_{\theta_q}(M) = r_M + j,
\]
where $r_M$ is the $\C[v^{1/2}]$-rank of $M$ and $j$ is the number of the nonzero invariant factors
whose specialization at $v^{1/2} = q$ is zero (\cref{thm: rank of tensored theta our cons}).
We next prove that $\calF_\iota$ and $\calF_{\theta_q}$ are essentially surjective (\cref{thm: essential surjectiveness}). 
Let $N$ be an object of the relevant target category. 
For $\calF_\iota$, we provide a preimage of $N$ that is free over $\C[v^{1/2}]$, and 
for $\calF_{\theta_q}$, we provide one that is $\C[v^{1/2}]$-torsion.
Finally, we deal with the characters $\chi_{\calF_\iota(M)}$ and $\chi_{\calF_{\theta_q}(M)}$.
In particular,
when $M$ is free over $\C[v^{1/2}]$,
we derive the identity
\[
\theta_q \circ (\chi_{\calF_\iota(M)} \circ \widetilde{\iota}\,) = \chi_{\calF_{\theta_q}(M)} \circ \widetilde{\theta_q},
\]
where $\widetilde{\iota}: \polyH \to \genericH$ and $\widetilde{\theta_q}: \polyH \to \compH{q^2}$ are the ring homomorphisms induced by $\iota$ and $\theta_q$, respectively
(\cref{thm: character and specialization}).

In \cref{sec: 0-Hecke modules from specializations}, we study the relationship between two families of generic Hecke algebra modules (skew Specht modules and Kazhdan--Lusztig cell modules) and two significant families of $0$-Hecke algebra modules.
In \cref{subsec: two kinds of notable 0 Hecke modules}, we introduce two  families of  $0$-Hecke modules in our consideration.
The first family consists of $0$-Hecke skew Specht modules $X_{\lambda/\mu}$, where $\lambda/\mu$ ranges over all skew partitions of $n$.
Each $X_{\lambda/\mu}$ has the set $\SYT(\lambda/\mu)$ of standard Young tableaux of shape $\lambda/\mu$ as a basis, and its image under $\ch$ is the skew Schur function $s_{\lambda/\mu}$.
The second family consists of the dual plactic class modules $\scrD_Q$, where  $Q$ ranges over $\SYT_n$, the set  of all standard Young tableaux of size $n$.
Each $\scrD_Q$ has the dual plactic class $\Gamma_Q$ of $Q$ as its basis, and its image under $\ch$ is the Schur function $s_{\sh(Q)'}$.

In \cref{sec: 0-Hecke modules from specializations}, we study the relationship between two families of generic Hecke algebra modules and two  families of $0$-Hecke algebra modules.
In \cref{subsec: two kinds of notable 0 Hecke modules}, we give a detailed introduction to two families of modules of the $0$-Hecke algebra: $0$-Hecke skew Specht modules and dual plactic class modules.
In \cref{subsec: 0-Hecke module from specialization 1}, we deal with the relationship between the skew Specht modules  and $0$-Hecke skew Specht modules.
For a principal ideal domain $R$ and $a \in R$, let  $S^{\lambda/\mu}(R;a)$ be the skew Specht module of $H_{R;a}(\SG_n)$ associated to $\lambda/\mu$.
We prove that 
\begin{align*}
\calF_{\iota}(S^{\lambda / \mu}(\C[v^{1/2}] ;v)) 
\cong S^{\lambda / \mu}(\C(v^{1/2}); v)
\quad \text{and} \quad
\calF_{\theta_0}(S^{\lambda / \mu}(\C[v^{1/2}] ;v))
\cong X_{\lambda/\mu}
\end{align*}
(\cref{thm: specialization of skew Specht module}).
The former equivalence follows  from \cref{cor: algebra homomorphism vs functor sub of regular}. 
The latter, however, cannot be deduced from this corollary because its hypothesis is not satisfied in this case (see \cref{rem: Specht and X}).
In fact, constructing an isomorphism between $\calF_{\theta_0}(S^{\lambda/\mu}(\C[v^{1/2}];v))$ and $X_{\lambda/\mu}$ is a nontrivial problem.
To solve this problem,  we introduce a subset $B(\lambda/\mu)$ of $S^{\lambda/\mu}(\C[v^{1/2}];v)$ parametrized by $\SYT(\lambda/\mu)$ and use an analogue of the Garnir relations to show that it is a $\C[v^{1/2}]$-basis (\cref{lem: skew specht basis}). 
We then analyze the $\polyH$-action on this basis to construct an explicit isomorphism.
In \cref{subsec: 0-Hecke module from specialization 2}, we deal with the relationship between the Kazhdan--Lusztig cell modules and the dual plactic class modules.
Let  $\frakM_Q(\C(v^{1/2});v)$ be the Kazhdan--Lusztig cell module associated with $Q \in \SYT_n$.
The purpose of this subsection is to find a module $M$ such that $\calF_\iota(M) \cong \frakM_Q(\C(v^{1/2});v)$ and $\calF_{\theta_0}(M) \cong \scrD_Q$.
As a first candidate, we canonically construct an $\polyH$-module $\bfM_Q^{(C)}$ using the Kazhdan--Lusztig basis $\{C_w \mid w \in \SG_n\}$ of $\laurH$ and the associated cell structure.
This module has a basis indexed by $\Gamma_Q$ (see \cref{eq: def of bfMCQ} for the definition).
This module satisfies that
 $\calF_{\iota}(\bfM^{(C)}_Q)
\cong \frakM_Q(\C(v^{1/2});v)$, however it is not satisfactory in the sense that  
\[
\calF_{\theta_0}(\bfM^{(C)}_Q)
\cong \bigoplus_{w \in \Gamma_Q} \bfF_{\comp(\Des_L(w))^\rmc} \not\cong \scrD_Q
\]
(\cref{thm: KL cell module C basis}).
This leads us to consider  a modification  of the Kazhdan--Lusztig basis.
More precisely, by using the basis $\{v^{\ell(w)/2} C_w \mid w \in \SG_n\}$ of $\laurH$ and the cell structure, 
 we construct an  $\polyH$-module $\bfM_Q^{(D)}$
  (for the definition, see \cref{def: def of MDQ}).
Then we show  
\begin{align*}
\calF_{\iota}(\bfM^{(D)}_Q)
\cong \frakM_Q(\C(v^{1/2});v) \quad \text{and} \quad
\calF_{\theta_0}(\bfM^{(D)}_Q)
\cong \scrD_Q
\end{align*}
(\cref{thm: obtain CGammaQ}).

In \cref{sec: properties of 0-Hecke modules specializatio}, we study characteristic maps.
To begin with, we review two character-theoretic criteria for $\genericH$-modules $M$ and $N$ to represent the same class in the Grothendieck group.
We then derive an analogue of these characterizations for $\compH{0}$-modules (\cref{prop: four equivalence relation on 0 Hecke}).
Using these characterizations, we demonstrate how the Frobenius characteristic and the quasisymmetric characteristic are related via specialization functors.
More precisely, we prove that for an $\polyH$-module $M$ that is free over $\C[v^{1/2}]$,
\[
\ch([\calF_{\theta_0}(M)])=\ch^{\rm Frob}([\calF_{\iota}(M)])
\]
(\cref{thm: characteristic pres prop spec}).
This relationship motivates the question of which $0$-Hecke modules have symmetric quasisymmetric characteristics.
We characterize these modules in terms of their characters (\cref{thm: char of symmetric 0-Hecke}).

In \cref{subsec:Relationship between restriction and specialization functors for skew Specht modules}, we study  the branching rules of the skew Specht modules $S^{\lambda/\mu}(\C(v^{1/2});v)$ and $0$-Hecke skew Specht modules $X_{\lambda/\mu}$, and connect them via specialization functors.
To do this, we introduce the  filtrations, called \emph{removable-corner filtrations},
\[
S_\bullet = (S_B)_{B \in \sfremov(\lambda/\mu)}, \quad \tilS_\bullet = (\tilS_B)_{B \in \sfremov(\lambda/\mu)}, \quad\text{and}\quad X_\bullet = (X_B)_{B \in \sfremov(\lambda/\mu)}
\]
of the restrictions of $S^{\lambda/\mu}(\C[v^{1/2}];v)$, $S^{\lambda/\mu}(\C(v^{1/2});v)$, and $X_{\lambda/\mu}$ to the relevant  subalgebras $H_{R;a}(\SG_{n-1})$, respectively.
Here, $\sfremov(\lambda/\mu)$ is the set of all removable boxes of the Young diagram of $\lambda/\mu$ with a suitable total order.
As the first result,
we decompose $\tilS_B$ and $X_B$ into direct sums of skew Specht modules and $0$-Hecke skew Specht modules, respectively, and prove that
\[
\calF_{\iota}(S_\bullet) \cong \tilS_\bullet
\quad \text{and} \quad
\calF_{\theta_0}(S_\bullet) \cong X_\bullet
\]
 (\cref{thm: skew Specht poly rest and specialization}).
The second result concerns the cyclic structure of $\tilS_B$ and $X_B$ when $\lambda/\mu$ is an ordinary partition (\cref{thm: generator of S and X as n-1}).
The cyclicity of $\tilS_B$ is immediate since it is isomorphic to a direct sum of pairwise nonisomorphic irreducible $H_{\C(v^{1/2});v}(\SG_{n-1})$-modules.
This decomposition does not, however, give an explicit cyclic generator. Here we construct explicit cyclic generators for $\tilS_B$ modeled on Young symmetrizers.
As for $X_B$, establishing its cyclicity is less immediate since its composition factors need not be multiplicity-free.  
We show the cyclicity by constructing an explicit cyclic generator.
It should be remarked that both $\tilS_B$ and $X_B$ can fail to be cyclic when $\lambda/\mu$ is not an ordinary partition (see \cref{rem: removable component is not cyclic}).

In \cref{subsec: Hom set Gamma Q and Gamma Qp}, we use the relationship established in \cref{subsec: 0-Hecke module from specialization 2} to obtain partial results on the Hom-spaces between dual plactic class modules and their classification up to isomorphism (\cref{thm: new thm for D Q iso to X}).
Let $Q_1,Q_2 \in \SYT_n$ with $\sh(Q_1) = \sh(Q_2)$.
Considering a nonzero homomorphism from $\frakM_{Q_1}(\C(v^{1/2});v)$ to $\frakM_{Q_2}(\C(v^{1/2});v)$, 
we construct a nonzero homomorphism $f^{\rm poly}_{Q_1, Q_2}$ from $\bfM^{(D)}_{Q_1}$ to $\bfM^{(D)}_{Q_2}$.
Then, by applying the specialization functor $\calF_{\theta_0}$ to $f^{\rm poly}_{Q_1, Q_2}$, we obtain a nonzero homomorphism $g_{Q_1, Q_2}$ from $\scrD_{Q_1}$ to $\scrD_{Q_2}$ (\cref{lem: new lem for g Q Q}).
With this preparation, we show that if at least one of $\Gamma_{Q_1}$ and $\Gamma_{Q_2}$ is a left weak Bruhat interval, then
\begin{align}\label{eq: Hom space is 1 dimnsional}
\Hom(\scrD_{Q_1}, \scrD_{Q_2}) = \C\{g_{Q_1, Q_2}\}.
\end{align}
We remark that if the assumption fails, then
  $ \dim\Hom(\scrD_{Q_1}, \scrD_{Q_2})$ can be greater than $1$
(for a concrete example, see \cref{rem: new remark for theorem on Hom}(2)).
For the classification, we endow every subset of $\SG_n$ with a colored digraph structure (see \cref{eq: colored digraph str def}).
Under the equality \cref{eq: Hom space is 1 dimnsional}, we derive the following classification theorem: 
\begin{align}\label{eq: classification in intro}
\scrD_{Q_1} \cong \scrD_{Q_2} \quad\text{if and only if}\quad \text{$\Gamma_{Q_1}$ and  $\Gamma_{Q_2}$ are isomorphic as colored digraphs.}
\end{align}
As a consequence, if exactly one of $\Gamma_{Q_1}$ and $\Gamma_{Q_2}$ is a left weak Bruhat interval, then $\scrD_{Q_1}$ and $\scrD_{Q_2}$ are not isomorphic (see \cref{rem: new remark for theorem on Hom}(1)).

\section{Preliminaries}
\label{sec: Prelim}

For integers $m$ and $n$, we set $[m,n]$ and $[n]$ to be the intervals $\{t\in \mathbb Z \mid m\le t \le n\}$ and $\{t\in \mathbb Z \mid 1\le t \le n\}$, respectively.
Throughout this paper, $n$ will denote a nonnegative integer unless otherwise stated. 
Every module will be assumed to be finitely generated.

\subsection{Compositions, Young diagrams, and  tableaux}\label{subsec: comp and diag}
A \emph{composition} $\alpha$ of $n$, denoted by $\alpha \models n$, is a finite ordered list of positive integers $(\alpha_1,\alpha_2,\ldots, \alpha_k)$ satisfying $\sum_{i=1}^k \alpha_i = n$.
We call $\alpha_i$ ($1 \le i \le k$) a \emph{part} of $\alpha$, $k =: \ell(\alpha)$ the \emph{length} of $\alpha$, and $n =:|\alpha|$ the \emph{size} of $\alpha$. 
And, we define the empty composition $\varnothing$ to be the unique composition of size and length $0$.
Whenever necessary, we set $\alpha_i = 0$ for all $i > \ell(\alpha)$.

Given $\alpha = (\alpha_1,\alpha_2,\ldots,\alpha_k) \models n$ and $I = \{i_1 < i_2 < \cdots < i_{l}\} \subseteq [n-1]$, 
let 
\begin{align*}
    \set(\alpha) &:= \{\alpha_1,\alpha_1+\alpha_2,\ldots, \alpha_1 + \alpha_2 + \cdots + \alpha_{k-1}\}, \text{ and} \\\comp(I) &:= (i_1,i_2 - i_1, i_3 - i_2, \ldots,n-i_{l}).
\end{align*}
The set of compositions of $n$ is in bijection with the set of subsets of $[n-1]$ under the correspondence $\alpha \mapsto \set(\alpha)$ (or $I \mapsto \comp(I)$).
The \emph{reverse composition $\alpha^\rmr$ of $\alpha$} is defined to be the composition $(\alpha_k, \alpha_{k-1}, \ldots, \alpha_1)$ and
the \emph{complement $\alpha^\rmc$ of $\alpha$} is defined to be the unique composition satisfying $\set(\alpha^c) = [n-1] \setminus \set(\alpha)$.

If a composition $\lambda = (\lambda_1, \lambda_2, \ldots, \lambda_k) \models n$ satisfies $\lambda_1 \ge \lambda_2 \ge \cdots \ge \lambda_k$, then it is called a \emph{partition} of $n$ and denoted as $\lambda \vdash n$.
The conjugate partition of $\lambda$, denoted by $\lambda'$, is the partition whose $i$th part is the number of parts of $\lambda$ that are greater than or equal to $i$.
Given two partitions $\lambda$ and $\mu$ with $\ell(\lambda) \ge \ell(\mu)$, we write $\lambda \supseteq \mu$ if $\lambda_i \geq  \mu_i$ for all $1 \le i \le \ell(\mu)$. 
A \emph{skew partition} $\lambda/\mu$ is a pair  $(\lambda, \mu)$ of partitions with $\lambda \supseteq \mu$.
We call $|\lambda/\mu| := |\lambda| - |\mu|$ the \emph{size} of $\lambda/\mu$.

Given a partition $\lambda$ of $n$, we define the \emph{Young diagram $\tyd(\lambda)$ of $\lambda$} to be the left-justified array of $n$ boxes, where the $i$th row from the top has $\lambda_i$ boxes for $1 \le i \le k$.
Similarly, given a skew partition $\lambda / \mu$, we define the \emph{Young diagram $\tyd(\lambda / \mu)$ of $\lambda / \mu$} to be the Young diagram $\tyd(\lambda)$ with all boxes belonging to $\tyd(\mu)$ removed.
We say that a box in $\tyd(\lambda/\mu)$ is \emph{in the $i$th row} if it is in the $i$th row from the top and \emph{in the $j$th column} if it is in the $j$th column from the left.
We denote by $(i,j)$ the box in the $i$th row and $j$th column.
A Young diagram is called \emph{connected}
if for each pair of consecutive rows, there are at least two boxes (one in each row) which have a common edge.
A skew partition is called \emph{connected} if the corresponding Young diagram is connected, and it is called \emph{basic} if the corresponding Young diagram contains neither empty rows nor empty columns.
In this paper, every skew partition is assumed to be  basic unless otherwise stated.

Given a skew partition $\lambda/\mu$ of  $n$, a \emph{bijective tableau} of shape $\lambda / \mu$ is a filling of $\tyd(\lambda / \mu)$ with distinct entries in $[n]$.
A bijective tableau is referred to as a \emph{standard Young tableau} if the elements in each row are arranged in increasing order from left to right, and the elements in each column are arranged in increasing order from top to bottom.
We denote by $\SYT(\lambda / \mu)$ the set of all standard Young tableaux of shape $\lambda / \mu$.

For later use, we introduce the following notation.

\begin{enumerate}[label = {\rm (\arabic*)}, leftmargin = 4ex, itemsep = 0.5ex]
\item 
For any filling $\ttT$ of  $\tyd(\lambda/\mu)$, let $\wcol(\ttT)$ be the word obtained by reading the entries of $\ttT$ from top to bottom starting with the leftmost column and $\wrow(\ttT)$ the word obtained by reading the entries of $\ttT$ from left to right starting with the topmost row.
When $\ttT$ is a bijective tableau, we regard $\wrow(\ttT)$ and $\wcol(\ttT)$ 
as permutations written in one-line notation.

\item 
Let $\ttT_{\lambda/\mu}$ be the filling of $\tyd(\lambda/\mu)$ obtained by filling $1, 2, \ldots, n$ from left to right starting with the topmost row and $\ttT'_{\lambda/\mu}$  the filling of $\tyd(\lambda/\mu)$ obtained by filling $1, 2, \ldots, n$ from top to bottom starting with the leftmost column.
\end{enumerate}

\subsection{Hecke algebras and their representations}

Let $\SG_n$ denote the symmetric group on $[n]$. Every permutation $\sigma \in \SG_n$ can be expressed as a product of simple transpositions $s_i := (i \ i+1)$ for $1 \leq i \leq n-1$.
A \emph{reduced expression for $\sigma$} is an expression that represents $\sigma$ in the shortest possible length, and 
the \emph{length $\ell(\sigma)$ of $\sigma$} is the number of simple transpositions in any reduced expression for $\sigma$.
The \emph{sign of $\sigma$}, denoted by $\sgn(\sigma)$, is defined by $(-1)^{\ell(\sigma)}$.
For $\alpha \models n$, define the \emph{Young subgroup of type $\alpha$} by
\[
\SG_\alpha := 
\left\langle s_i \; \middle| \; i \in \set(\alpha^\rmc) \right\rangle.
\]
We denote the longest element of $\SG_\alpha$ by $w_0(\alpha)$.

\subsubsection{Hecke algebras}\label{subsec: Hecke alg}
Let $R$ be a principal ideal domain and $a \in R$.
The \emph{Iwahori-Hecke algebra} $H_{R;a}(\SG_n)$ associated with the symmetric group $\SG_n$ is the $R$-algebra with identity generated by $T_1, T_2, \ldots, T_{n-1}$ subject to the following relations
\begin{alignat}{2}
\label{eq: def of Hecke alg}
T_i^2 &= (a-1)T_i + a
&\quad& \text{for $1\le i \le n-1$}, \nonumber \\ 
T_i T_{i+1} T_i &= T_{i+1} T_i T_{i+1}
&\quad& \text{for $1\le i \le n-2$}, \\ 
T_i T_j &= T_j T_i
&\quad& \text{if $|i-j| \ge 2$.}
\nonumber
\end{alignat}
For any reduced expression $s_{i_1} s_{i_2} \cdots s_{i_p}$ for $\sigma \in \SG_n$, let 
\[
T_{\sigma} := T_{i_1} T_{i_2} \cdots T_{i_p}.
\]
This element is independent of the choices of reduced expressions and $\{T_\sigma \mid \sigma \in \SG_n\}$ is an $R$-basis for $H_{R;a}(\SG_n)$.

Let us introduce the  elements of the Hecke algebra $H_{R;a}(\SG_n)$ that will be essential for our discussion.
For $S \subseteq \SG_n$, let 
\begin{align}\label{eq: def of bfx and bfy}
\bfx_S := \sum_{\sigma \in S} T_\sigma 
\quad \text{and} \quad 
\bfy_S := \sum_{\sigma \in S} \sgn(\sigma) a^{l_0 - \ell(\sigma)}T_\sigma,
\end{align}
where $l_0 := \max\{\ell(\sigma) \mid \sigma \in S \}$. And,
for $\alpha \models n$, let $\bfx_\alpha := \bfx_{\SG_\alpha}$
 and $\bfy_\alpha := \bfy_{\SG_\alpha}$.
For later use, we record the following identities from \cite[p.~843]{86Gyoja}:
for $\alpha \models n$ and $\sigma \in \SG_\alpha$,
\begin{align}\label{eq: Ti and bfy}
T_\sigma \bfx_\alpha = a^{\ell(\sigma)} \bfx_\alpha
= \bfx_\alpha T_\sigma 
\quad\text{and}\quad 
T_\sigma \bfy_\alpha 
= \sgn(\sigma) \bfy_\alpha
=
\bfy_\alpha T_\sigma.
\end{align}

In this paper, we study the representation theories of $\genericH$ and $\compH{0}$ through representations of $\polyH$, where $v^{1/2}$ is an indeterminate.
We refer to $\genericH$ as the \emph{generic Hecke algebra} and to $\compH{0}$ as the \emph{$0$-Hecke algebra}.

\subsubsection{Irreducible modules of Hecke algebras}
Here we introduce the irreducible modules of the generic Hecke algebra $\genericH$ and the $0$-Hecke algebra $\compH{0}$.

It is well known that 
the generic Hecke algebra $\genericH$ is isomorphic to the group algebra $\C(v^{1/2})[\SG_n]$ (see \cite{02Bourbaki, 81Lusztig}).
Therefore, it is semisimple and its irreducible modules are naturally parametrized by partitions of $n$.
For $\lambda \vdash n$, we define the \emph{Young symmetrizer associated to $\lambda$} in $\genericH$ by
\begin{align*}
\bfz_{\lambda}(\C(v^{1/2});v) := \bfy_{\lambda} T_{\wcol(\ttT_{\lambda})} \bfx_{\lambda'}.
\end{align*}
Then, we define the \emph{Specht module  of $\genericH$ associated to $\lambda$} by 
\begin{align*}
S^{\lambda}(\C(v^{1/2});v) := \genericH \bfz_{\lambda}(\C(v^{1/2});v).
\end{align*}
By the slight modification of \cite[Theorem ]{86DJ}, one can see that 
 $\{S^\lambda(\C(v^{1/2});v) \mid \lambda \vdash n \}$ is a complete list of pairwise nonisomorphic irreducible  $\genericH$-modules.

Unlike the generic Hecke algebra,
the $0$-Hecke algebra $\compH{0}$ is neither semisimple nor representation-finite for $n \ge 4$; it is of tame type when $n = 4$ and of wild type for $n > 4$ (see~\cite{11DY, 02DHT}). 
The classification of irreducible   $\compH{0}$-modules up to isomorphism was provided in  \cite{79Norton}.
There are  $2^{n-1}$ pairwise nonisomorphic irreducible $\compH{0}$-modules, which are naturally indexed by compositions of $n$.
For $\alpha \models n$, let $\bfF_{\alpha}:=\C v_{\alpha}$
and endow it with the $\compH{0}$-action given by
for each $1 \le i \le n-1$,
\begin{align}\label{eq: simple module 0-Hecke action}
T_i \cdot v_\alpha = \begin{cases}
0 & i \in \set(\alpha),\\
-v_\alpha & i \notin \set(\alpha).
\end{cases}
\end{align}
Then $\{\bfF_\alpha \mid \alpha \models n \}$ is a complete list of irreducible $\compH{0}$-modules.

\subsection{Characteristic maps}

Let $H_{R;a}(\SG_n)$-$\module$ be the category of finitely generated  $H_{R;a}(\SG_n)$-modules and $\calR(H_{R;a}(\SG_n)\text{-}\module)$ the $\Z$-span of the set of (representatives of) isomorphism classes of modules in $H_{R;a}(\SG_n)$-$\module$.
We denote by $[M]$ the isomorphism class corresponding to a finitely generated $H_{R;a}(\SG_n)$-module $M$. 
When $R$ is a field,  
we define the \emph{Grothendieck group $\calG_0(H_{R;a}(\SG_n)\text{-}\module)$ of $H_{R;a}(\SG_n)$-$\module$} by the quotient of $\calR(H_{R;a}(\SG_n)\text{-}\module)$ modulo the relations $[M] = [M'] + [M'']$ whenever there exists a short exact sequence $0 \ra M' \ra M \ra M'' \ra 0$. 
The equivalence classes of the irreducible $H_{R;a}(\SG_n)$-modules form a $\Z$-basis for $\calG_0(H_{R;a}(\SG_n)\text{-}\module)$.
Let
\[
\calG_0(H_{R;a}(\SG_\bullet)\text{-}\module) := \bigoplus_{n \ge 0} \calG_0(H_{R;a}(\SG_n)\text{-}\module).
\]
In the cases where $R$ is $\C$ or $\C(v^{1/2})$,
$\calG_0(H_{R;a}(\SG_\bullet)\text{-}\module)$ carries a Hopf algebra structure whose multiplication and comultiplication are induced by induction product and restriction, respectively (for instance, see \cite[Corollary 4.6 and Remark 4.7]{15SY}).

\subsubsection{The Frobenius characteristic}\label{subsubsec: Frob char}
A formal power series $f \in \Z[[x_1, x_2, \ldots]]$ is called a \emph{symmetric function} if it is of bounded degree and for every partition $(\lambda_1, \lambda_2, \ldots, \lambda_k)$, the coefficient of $x_{i_1}^{\lambda_1} x_{i_2}^{\lambda_2} \cdots x_{i_k}^{\lambda_k}$ in $f$ is equal to that of  $x_{1}^{\lambda_1} x_{2}^{\lambda_2} \cdots x_{k}^{\lambda_k}$ in $f$ for any sequence $(i_1, i_2, \cdots, i_k)$ of distinct positive integers.
Let $\Sym$ be the subring of $\Z[[x_1, x_2, \ldots]]$ consisting of symmetric functions.
Recall that $\Sym$ is a graded connected Hopf algebra
\[
\Sym = \bigoplus_{n \geq 0} \Sym_n,
\]
where $\Sym_n$ is the $\Z$-submodule of $\Sym$ consisting of homogeneous symmetric functions of degree $n$.
For each partition $\lambda$, let $s_{\lambda}$ be the \emph{Schur function} attached to $\lambda$. 
As $\lambda$ ranges over the partitions of $n$,  $s_{\lambda}$'s  form a basis of $\Sym_n$ (for example, see \cite{97Ful, 99Stanley}).

Consider the $\Z$-linear isomorphism
\begin{align}\label{eq: Frob char}
\ch^{\rm Frob}: \calG_0(H_{\C(v^{1/2});v}(\SG_\bullet)\text{-}\module) \ra \Sym
\end{align}
given by sending
 $[S^\lambda(\C(v^{1/2});v)]$ to $s_\lambda$ for every partition $\lambda$.
This map is known to be a Hopf algebra isomorphism (see \cite[Theorem 5.7]{14WW} and \cite[Proposition 2.9(ii)]{21Korff}).
We call this isomorphism as  the \emph{Frobenius characteristic}.

\subsubsection{The quasisymmetric characteristic}
A formal power series $f \in \Z[[x_1, x_2, \ldots]]$ is called a \emph{quasisymmetric function} if it is of bounded degree and for every composition $(\alpha_1, \alpha_2, \ldots, \alpha_k)$, the coefficient of $x_{i_1}^{\alpha_1} x_{i_2}^{\alpha_2} \cdots x_{i_k}^{\alpha_k}$ in $f$ is equal to that of  $x_{1}^{\alpha_1} x_{2}^{\alpha_2} \cdots x_{k}^{\alpha_k}$ in $f$ for any sequence $(i_1, i_2, \ldots i_k)$ with  $i_1 < i_2 < \cdots <i_k$.
Let $\Qsym$ be the subring of $\Z[[x_1, x_2, \ldots]]$ consisting of quasisymmetric functions.
Recall that $\Qsym$ is also a graded connected Hopf algebra
\[
\Qsym = \bigoplus_{n \geq 0} \Qsym_n,
\]
where $\Qsym_n$ is the $\Z$-submodule of $\Qsym$ consisting of homogeneous quasisymmetric functions of degree $n$.
For $\alpha \models n$ with $n \geq 1$, define 
\[
F_\alpha := \sum_{\substack{1 \le i_1 \le i_2 \le \cdots \le i_n \\ i_j < i_{j+1} \text{ if } j \in \set(\alpha)}} x_{i_1}x_{i_2} \cdots x_{i_n}
\]
and 
$F_{\varnothing} := 1$.
The function $F_\alpha$ is called \emph{the fundamental quasisymmetric function} attached to $\alpha$.
As $\alpha$ ranges over the compositions of $n$, $F_\alpha$'s form a basis of $\Qsym_n$.

Consider the $\Z$-linear isomorphism
\[
\ch: \calG_0(H_{\C;0}(\SG_\bullet)\text{-}\module) \ra \Qsym
\]
given by sending
 $[\bfF_\alpha]$ to $F_\alpha$ for every composition $\alpha$.
It was shown in \cite{96DKLT, 09BL} that this map is a  Hopf algebra isomorphism.
This isomorphism is called the 
\emph{quasisymmetric characteristic}.

\bigskip 
Before closing this section, we collect definitions and notation related to the weak Bruhat orders on $\SG_n$, which will be used repeatedly from \cref{sec: 0-Hecke modules from specializations} onward.
For $\sigma \in \SG_n$, let
\[
\Des_L(\sigma):= \{i \in [n-1] \mid \ell(s_i \sigma) < \ell(\sigma)\}
\ \  \text{and} \ \  
\Des_R(\sigma):= \{i \in [n-1] \mid \ell(\sigma s_i) < \ell(\sigma)\}.
\]
It is well known that if $\sigma = w_1 w_2 \cdots w_n$ in one-line notation, then
\begin{align*}
\begin{aligned}
\Des_L(\sigma) & = \{ i \in [n-1] \mid \text{$i$ is right of $i+1$ in $w_1 w_2 \cdots w_n$} \} \quad \text{and} \\
\Des_R(\sigma) & = \{ i \in [n-1] \mid w_i > w_{i+1}
\}.
\end{aligned} 
\end{align*}
The \emph{left weak Bruhat order} $\preceq_L$ 
(resp. \emph{right weak Bruhat order} $\preceq_R$) on $\SG_n$
is the partial order on $\SG_n$ whose covering relation $\preceq_L^{\rmc}$ (resp. $\preceq_R^\rmc$) is given as follows: 
\[
\sigma \preceq_L^\rmc s_i \sigma \text{ if }i \notin \Des_L(\sigma)
\quad \text{(resp. $\sigma \preceq_R^\rmc \sigma s_i$ if  $i \notin \Des_R(\sigma)$).}
\]
Although these two weak Bruhat orders are not identical,  
there exists a poset isomorphism $$(\SG_n, \preceq_L)\to (\SG_n, \preceq_R),\quad \sigma \mapsto \sigma^{-1}.$$
Given $\sigma, \rho \in \SG_n$, the \emph{left weak Bruhat interval $[\sigma, \rho]_L$} (resp. the \emph{right weak Bruhat interval $[\sigma, \rho]_R$}) denotes the closed interval $\{\gamma \in \SG_n \mid \sigma \preceq_L \gamma \preceq_L \rho \}$ (resp. $\{\gamma \in \SG_n \mid \sigma \preceq_R \gamma \preceq_R \rho \}$) with respect to the left weak Bruhat order (resp. the right weak Bruhat order).

\section{Specialization functors
\texorpdfstring{$\calF_\iota$ and $\calF_{\theta_q}$}
{F iota and F theta q}}
\label{sec: specialization functor}

\subsection{Definitions and notation}\label{subsec: def of specialization functor}
Let $R$ and $S$ be principal ideal domains (PIDs), and $a \in R$ and $b \in S$.
Suppose there exists a nonzero ring homomorphism $f: R \ra S$ such that $f(a)=b$.
Associated to the triple $(f,a,b)$, we define a functor $\calF_{f,a,b}$, which sends an  $H_{R;a}(\SG_n)$-module to an $H_{S;b}(\SG_n)$-module.
To do this, we first observe that $f$ induces a ring homomorphism
$$
\widetilde{f}: H_{R;a}(\SG_n) \ra H_{S;b}(\SG_n), \qquad \sum_{\sigma \in \SG_n}a_\sigma T_{\sigma,R} \mapsto \sum_{\sigma \in \SG_n} f(a_\sigma) T_{\sigma,S}.
$$
Via this homomorphism, 
$H_{S;b}(\SG_n)$ acquires the structure of an 
$(H_{S;b}(\SG_n), H_{R;a}(\SG_n))$-bimodule structure. Therefore,
we obtain a functor
\[
\calF_{f, a, b}: H_{R;a}(\SG_n)\text{-}\module \ra H_{S;b}(\SG_n)\text{-}\module
\]
such that for a finitely generated  $H_{R;a}(\SG_n)$-module $M$,
\[
\calF_{f, a, b}(M) = H_{S;b}(\SG_n) \otimes_{H_{R;a}(\SG_n)} M,
\]
and for an $H_{R;a}(\SG_n)$-module homomorphism $\varphi:M \to N$,
\[
\calF_{f,a, b}(\varphi) = {\rm id}\otimes \varphi: \calF_{f,a, b}(M) \ra \calF_{f,a, b}(N).
\]
We call the triple $(f, a, b)$ a \emph{specialization}, and $\calF_{f,a ,b}$ is the \emph{specialization functor} induced by $(f, a, b)$.
When $a$ and $b$ are clear from the context, we simply write $\calF_{f}$ for $\calF_{f,a ,b}$.
We note that $\calF_f(M)$ is isomorphic to $S \otimes_R M$ as an $S$-module.

Let us introduce the specialization functors in our consideration.
Let $v^{1/2}$ be an indeterminate and $q$ be an arbitrary complex number.
We define $\iota: \C[v^{1/2}] \ra \C(v^{1/2})$ to be the natural embedding and 
$\theta_q: \C[v^{1/2}] \ra \C$ to be the evaluation homomorphism at $v^{1/2} = q$.
In this paper, we focus on the specialization functors 
\[
\calF_\iota  \quad\text{and}\quad\calF_{\theta_q}
\]
induced by $(\iota, v, v)$ and $(\theta_q, v, q^2)$, respectively.
Given an $\polyH$-module $M$,
these specialization functors allow us to relate  $\calF_\iota(M)$ to $\calF_{\theta_q}(M)$.

\subsection{Properties of \texorpdfstring{$\calF_\iota$ and $\calF_{\theta_q}$}
{F iota and F theta q}}\label{subsec: properties of specialization functor}

In this subsection, we establish three noteworthy
properties of the specialization functors $\calF_\iota$ and  $\calF_{\theta_q}$ ($q \in \C$).
Before doing this, we introduce necessary notations.
Let $N$ be a finitely generated $\C[v^{1/2}]$-module.

We denote by $r_N$ the rank of $N$,  and by $N_{\rmtor}$ the torsion submodule of $N$, that is,
\[
N_{\rmtor} = \{x \in N \mid r x = 0 \;\text{for some $0 \neq r \in \C[v^{1/2}]$} \}.
\]
Let $M$ be a finitely generated $\polyH$-module.
We still use  the notation $r_M$ for its $\C[v^{1/2}]$-rank and $M_{\rmtor}$ for its $\C[v^{1/2}]$-torsion submodule.
One can easily check that $M_{\rmtor}$ is an $\polyH$-submodule of $M$.

The first property  we consider concerns the dimensions of $\calF_\iota(M)$ and $\calF_{\theta_q}(M)$.

\begin{theorem}\label{thm: rank of tensored theta our cons}
Let $M$ be  a finitely generated  $\polyH$-module.
\begin{enumerate}[label = {\rm (\arabic*)}]
\item The $\C(v^{1/2})$-dimension of $\calF_\iota(M)$ is $r_M$.
\item For $q \in \C$, the $\C$-dimension of $\calF_{\theta_q} (M)$ is $r_M + j$.
Here, $j$ is the number of the nonzero invariant factors $p(v^{1/2})$ of $M$ such that $p(v^{1/2})|_{v^{1/2} = q} = 0$.
\end{enumerate}
\end{theorem}

\begin{proof}
Let $(S,f) \in \{(\C(v^{1/2}), \iota)\} \cup \{ (\C, \theta_q) \mid q \in \C\}$.
Since $\C[v^{1/2}]$ is a PID, 
by the structure theorem for finitely generated modules over a PID, 
\[
M \cong \C[v^{1/2}]^{\oplus r_M} \oplus
\left( \bigoplus_{i=1}^k \C[v^{1/2}]\,/\,(p_i) \right),
\]
where $p_1, p_2, \ldots, p_k$ are the nonzero invariant factors 
of $M$.
Since $S$ is equipped with the $(S, \C[v^{1/2}])$-bimodule structure via $f$, 
applying the functor $S \otimes_{\C[v^{1/2}]} \text{-}$ to both sides, we have an $S$-module isomorphism
\begin{equation}\label{eq: S tensor including free part}
S \otimes_{\C[v^{1/2}]} M \cong 
S^{\oplus r_M} 
\oplus 
\left( \bigoplus_{i=1}^k S \otimes_{\C[v^{1/2}]} \C[v^{1/2}]\,/\,(p_i) \right).
\end{equation}

Let $1 \le i \le k$.
Consider the short exact sequence of $\C[v^{1/2}]$-modules
\[
\begin{tikzcd}
0 \arrow[r] & (p_i) \arrow[r, "\alpha"] & \C[v^{1/2}] \arrow[r, "\beta"] & \C[v^{1/2}] \, / \, (p_i) \arrow[r] & 0,
\end{tikzcd}
\]
where $\alpha$ is the natural embedding and $\beta$ is the canonical projection.
Since the tensor functor $S \otimes_{\C[v^{1/2}]} \text{-}$ is right exact, we have the right exact sequence
\[
\begin{tikzcd}
S \otimes_{\C[v^{1/2}]} (p_i) \arrow[r, "1 \otimes \alpha"] 
& S \otimes_{\C[v^{1/2}]} \C[v^{1/2}] \arrow[r, "1 \otimes \beta"] 
& S \otimes_{\C[v^{1/2}]} \C[v^{1/2}] \, / \, (p_i) \arrow[r] 
& 0.
\end{tikzcd}
\]
It follows that
\[
S \otimes_{\C[v^{1/2}]}\C[v^{1/2}] \, / \, (p_i)
\cong 
(S \otimes_{\C[v^{1/2}]} \C[v^{1/2}]) \, / \, \ker(1 \otimes \beta) = (S \otimes_{\C[v^{1/2}]} \C[v^{1/2}]) \, / \, \rmIm(1 \otimes \alpha).
\]
Moreover, since  $(1 \otimes \alpha) (1 \otimes p_i) = 1 \otimes p_i =  f(p_i) \otimes 1$ and $S$ is a field,  we have
\[
\rmIm(1 \otimes \alpha) \cong \begin{cases}
0  & \text{if $f(p_i) = 0$,}\\
S & \text{otherwise.}
\end{cases}
\]
Putting these together yields that
\begin{equation}\label{eq: dimension specialization functor to torsion}
S \otimes_{\C[v^{1/2}]} \C[v^{1/2}] \, / \, (p_i) 
\cong \begin{cases}
S  & \text{if $f(p_i) = 0$,}\\
0  & \text{otherwise.}
\end{cases}
\end{equation}
Now, the assertions (1) and (2) follow from 
\cref{eq: S tensor including free part} together with \cref{eq: dimension specialization functor to torsion}.
\end{proof}

\begin{remark}\label{rem: poly action and hecke action under free cond}
In \cref{thm: rank of tensored theta our cons}, consider the case where $M$ is free over $\C[v^{1/2}]$ with a basis $B$.
Let $(S,f)$ be either $(\C(v^{1/2}), \iota)$ or $(\C, \theta_q)$ for some $q \in \C$.
Then,  
$1 \otimes B$ is a basis for $\calF_f(M)$.
The $H_{S;f(v)}(\SG_n)$-action on $\calF_f(M)$ can be described as follows.
For $\sigma \in \SG_n$ and $b \in B$, let 
\begin{equation*}
T_\sigma \cdot b = \sum_{b' \in B}p^{\sigma}_{b, b'}(v^{1/2}) b'.
\end{equation*}
The $H_{S;f(v)}(\SG_n)$-action on $\calF_f(M)$ is given by 
\begin{equation}\label{eq: specialization and module action}
T_\sigma \cdot 1 \otimes b = \sum_{b' \in B} f\big(p^{\sigma}_{b, b'}(v^{1/2})\big)\, 1 \otimes b'.
\end{equation}
\end{remark}

\begin{corollary}\label{cor: algebra homomorphism vs functor sub of regular}
Let $M$ and $N$ be submodules of the regular representation of $\polyH$ such that $N \subseteq M$ and $M/N$ is free over $\C[v^{1/2}]$.
\begin{enumerate}[label={\rm (\arabic*)}]
\item 
We have an isomorphism of $\genericH$-modules
\[
\calF_{\iota}(M/N) \cong 
\inn{\widetilde{\iota}(M)}/ 
\inn{\widetilde{\iota}(N)},
\]
where 
$\inn{\widetilde{\iota}(M)}$ and
$\inn{\widetilde{\iota}(N)}$ are the submodules of
$H_{\C(v^{1/2}); v}(\SG_n)$ generated by
$\widetilde{\iota}(M)$ and $\widetilde{\iota}(N)$, respectively.

\item 
For $q \in \C$, we have $\dim_{\C} \widetilde{\theta_q}(M) \leq r_M$ and $\dim_{\C} \widetilde{\theta_q}(N) \leq r_N$.
If
\[
\dim_{\C}\bigl(\widetilde{\theta_q}(M)/\widetilde{\theta_q}(N)\bigr)
= r_{M/N},
\]
then we have an isomorphism of $\compH{q^2}$-modules
\[
\calF_{\theta_q}(M/N)
\cong
\widetilde{\theta_q}(M)/\widetilde{\theta_q}(N).
\]
\end{enumerate}
\end{corollary}
\begin{proof}
To prove (1) and (2) simultaneously, let $(S, f)$ be either $(\C(v^{1/2}), \iota)$ or $(\C, \theta_q)$, and set 
\[
M' = \inn{\widetilde{f}(M)} \quad\text{and}\quad N' = \inn{\widetilde{f}(N)}.
\] 
Since $M$, $N$, and $M/N$ are free over $\C[v^{1/2}]$, we can choose bases $B_N$ for $N$ and $B_{M}$ for $M$ such that $B_N \subseteq B_M$.
Since  $M'$ is spanned by $\widetilde{f}(B_M)$ over $S$, we have 
\[
\dim_S M' \leq |\widetilde{f}(B_M)| \leq |B_M| =  r_M.
\]
Also, since $N'$ is spanned by $\widetilde{f}(B_N)$ over $S$, we have
\[
\dim_S N' \leq |\widetilde{f}(B_N)| \leq |B_N| =  r_N.
\]

It remains to show that $\calF_f(M/N) \cong M'/N'$ under the assumption that $\dim_S M'/N' = r_{M/N}$ because this assumption always hold when $(S,f) = (\C(v^{1/2}), \iota)$.
Suppose that $\dim_S M' / N' = r_{M/N}$.
Let $B := B_M \setminus B_N$.
Then, the set $\{\widetilde{f}(b) + N' \mid b \in B \}$ spans  $M'/N'$. 
Since $|B| = r_{M/N}$, this set is an $S$-basis for $M'/N'$. 
On the other hand, by \cref{thm: rank of tensored theta our cons}, the set $\{1 \otimes (b + N) \mid b \in B \}$ is an $S$-basis for $\calF_{f}(M/N)$.

Let $\varphi: \calF_{f}(M/N) \ra M'/N'$
be the $S$-linear isomorphism sending $1 \otimes (b + N)$ to $\widetilde{f}(b)$ for $b \in B$.
For all $1 \le i \le n-1$ and $b \in B$,
\begin{align*}
\varphi(T_i \cdot (1\otimes (b + N))) 
&= \varphi(1\otimes (T_i  b + N))\\
&= \widetilde{f}(T_i b) + N'\\
&= T_i \cdot \widetilde{f}(b) + N'\\
&= T_i \cdot \varphi(1 \otimes (b+N)).
\end{align*}
Therefore, $\varphi$ is an $H_{S;f(v)}(\SG_n)$-module isomorphism.
\end{proof}

\begin{remark}
One can find examples such that the inequality in \cref{cor: algebra homomorphism vs functor sub of regular}(2) is strict. 
For example, if $M = \polyH (v-q) T_{\id}$, then $\calF_{\theta_q}(M) \cong \compH{q^2}$, whereas $\widetilde{\theta_q}(M) = 0$.
Another example can be found in \cref{rem: Specht and X}.
\end{remark}

The second property concerns the essential surjectivity of $\calF_\iota$ and $\calF_{\theta_q}$.
Recall that for categories $\calC$ and $\calD$, a functor $\calF: \mathcal C \ra \mathcal D$ is called \emph{essentially surjective} if every object $N$ of $\mathcal D$ is isomorphic to $\calF(M)$ for some object $M$ of $\mathcal C$.

\begin{theorem}\label{thm: essential surjectiveness}
\begin{enumerate}[label = {\rm (\arabic*)}]
    \item  The functor $\calF_\iota$ is essentially surjective.
    \item 
For $q \in \C$, the functor $\calF_{\theta_q}$ is essentially surjective.
\end{enumerate}
\end{theorem}

\begin{proof}
(1)
By \cref{cor: algebra homomorphism vs functor sub of regular}(1), for each $\lambda \vdash n$, we have 
\begin{equation*}
\calF_\iota\bigl(S^{\lambda}(\C[v^{1/2}];v)\bigr) \cong S^{\lambda}(\C(v^{1/2});v) 
\quad \text{as $\genericH$-modules.}
\end{equation*}
Since the category $\genericH$-$\module$ is semisimple and 
$\{ S^{\lambda}(\C(v^{1/2});v) \mid \lambda \vdash  n \}$
forms a complete set of non-isomorphic irreducible $\genericH$-modules, our assertion follows.

\medskip

(2) 
Let $N$ be a finitely generated  $\compH{q^2}$-module
and $\{x_1, x_2, \ldots, x_k\}$ a basis for $N$.
Let 
\begin{equation*}
M := \C[v^{1/2}]/(v - q^2) \otimes_{\C} N.
\end{equation*}
For all $1 \leq i \leq n-1$, define an action of $\mathring{T}_i$ on $M$ by
\begin{align}\label{eq: sfT actions}
\mathring{T}_i \cdot (1 \otimes x_j) := 1 \otimes (T_i \cdot x_j) \quad \text{for $1 \le j \le k$.}
\end{align}
Here, $\mathring{T}_i$'s denote the generators of $\polyH$ and $T_i$'s denote those of $\compH{q^2}$.
It is straightforward to see that $\mathring{T}_i$'s satisfy the braid relations.
They also satisfy the quadratic relations
$\mathring{T}_i^2 = (v - 1) \mathring{T}_i + v~(1 \le i \le n-1)$ since 
\begin{align*}
\mathring{T}_i^2 \cdot (1 \otimes x_j)
& = 1 \otimes T_i^2 \cdot x_j = 1 \otimes ((q^2-1)T_i + q^2)\cdot x_j
= ((v-1) \mathring{T}_i + v) \cdot (1 \otimes x_j).
\end{align*}
Thus, $M$ is an $\polyH$-module with the action given in 
\cref{eq: sfT actions}.

Note that
\begin{align*}
M & =  \C[v^{1/2}] \, / \, (v-q^2)  \otimes_{\C}  \left( \bigoplus_{j=1}^k \C x_j \right)
\cong \bigoplus_{j=1}^k \, \C[v^{1/2}] \, / \, (v-q^2)  \  \text{as $\C[v^{1/2}]$-modules.}
\end{align*}
This implies that $M$ is a $\C[v^{1/2}]$-torsion module with exactly $k$ nonzero invariant factors, all equal to $v-q^2$.
By \cref{thm: rank of tensored theta our cons}(2), the $\C$-dimension of $\calF_{\theta_q}(M)$ is equal to $k$.
Since
\begin{align*}
\calF_{\theta_q}(M)
& = \calF_{\theta_q}(\C[v^{1/2}] \, / \, (v - q^2) \otimes_{\C} N) 
= \C \, \otimes_{\C[v^{1/2}]} \, \C[v^{1/2}] \, / \, (v - q^2)\, \otimes_{\C} \, N,
\end{align*}
the set $\{1 \otimes 1 \otimes x_j \mid 1 \le j \le k\}$ spans $\calF_{\theta_q}(M)$, and thus it is a $\C$-basis for $\calF_{\theta_q}(M)$.
Moreover, for each $1\le j \le k$,
\[
T_i \cdot (1 \otimes 1 \otimes x_j)
= 1 \otimes \mathring{T}_i \cdot (1 \otimes x_j)
= 1 \otimes 1 \otimes T_i \cdot x_j \quad (1 \le i \le n-1).
\]
Consequently, the $\C$-linear map
\[
\calF_{\theta_q}(M) \ra N, \quad 1 \otimes 1 \otimes x_j \mapsto x_j  \quad (1 \le j \le k)
\]
is an $\compH{q^2}$-module isomorphism.
\end{proof}

\begin{remark}
In the proof of \cref{thm: essential surjectiveness}(2), the preimage $M$ of $N$ under $\calF_{\theta_q}$ that we found there is a $\C[v^{1/2}]$-torsion module.
Therefore,  \cref{thm: essential surjectiveness}(2) does not tell us the existence of a preimage of $N$ that is free over $\C[v^{1/2}]$.
\end{remark}

The third property  concerns the characters of $\calF_\iota(M)$ and $\calF_{\theta_q}(M)$. 
Let $R$ be a field and $a \in R$.
For any finitely generated $H_{R;a}(\SG_n)$-module $M$, the \emph{character} of $M$ is the map
\[
\chi_M : H_{R;a}(\SG_n) \ra R, \quad h \mapsto \sftr_{M}(h).
\]

\begin{theorem}\label{thm: character and specialization}
Let $M$ be a finitely generated $\polyH$-module.
\begin{enumerate}[label = {\rm (\arabic*)}]
\item 
Let $f$ be either $\iota$ or $\theta_q$ for some $q \in \C$.
Then
\[
\chi_{\calF_{f}(M)} = \chi_{\calF_{f}(M_{\rmtor})} + \chi_{\calF_{f}(M/M_{\rmtor})}.
\]

\item
Assume that $M$ is free over $\C[v^{1/2}]$.
Then, the image of $\chi_{\calF_{\iota}(M)} \circ \widetilde{\iota}$ lies in $\C[v^{1/2}]$ and 
the following diagram
\[
\begin{tikzcd}[column sep=3cm, row sep=1.5cm]
{\polyH} \arrow[d, "\widetilde{\theta_q}"'] \arrow[r, "{\chi_{\calF_\iota(M)} \; \circ \;\widetilde{\iota}}"] & {\C[v^{1/2}]} \arrow[d, "\theta_q"] \\
\compH{q^2} \arrow[r, "\chi_{\calF_{\theta_q}(M)}"] & \C 
\end{tikzcd}
\]
 commutes.
\end{enumerate}
\end{theorem}

\begin{proof}
(1) 
Let $M$ be a finitely generated $\polyH$-module.
Set $S = \C(v^{1/2})$ if $f = \iota$, and $S = \C$ otherwise.
For the assertion, it suffices to show that there is a short exact  sequence of $H_{S;f(v)}(\SG_n)$-modules 
\begin{equation*}
\begin{tikzcd}
0 \arrow[r] & \calF_f(M_{\rmtor}) \arrow[r] & \calF_f(M) \arrow[r] & \calF_f(M/M_{\rmtor}) \arrow[r] & 0
\end{tikzcd}.
\end{equation*}
Let \[
\begin{tikzcd}
0 \arrow[r] & M_{\rmtor} \arrow[r, "\alpha"] & M \arrow[r, "\beta"] & M/M_{\rmtor} \arrow[r] & 0
\end{tikzcd}
\]
be the canonical short exact sequence of $\polyH$-modules.
It splits as a sequence of $\C[v^{1/2}]$-modules, and thus
\[
\begin{tikzcd}
0 \arrow[r] & S \otimes_{\C[v^{1/2}]} M_{\rmtor} \arrow[r, "\id \otimes \alpha"] & S \otimes_{\C[v^{1/2}]} M \arrow[r, "\id \otimes \beta"] & S \otimes_{\C[v^{1/2}]} M/M_{\rmtor} \arrow[r] & 0
\end{tikzcd}
\]
is  exact as a sequence of $S$-modules.
By the definition of $\calF_f$, this sequence can be viewed as the following short exact sequence of $H_{S; f(v)}(\SG_n)$-modules:
\[
\begin{tikzcd}
0 \arrow[r] & \calF_f(M_{\rmtor}) \arrow[r, "\calF_f(\alpha)"] & \calF_f(M) \arrow[r, "\calF_f(\beta)"] & \calF_f(M/M_{\rmtor}) \arrow[r] & 0
\end{tikzcd}
\]
\medskip

(2) 
Let $B$ be a $\C[v^{1/2}]$-basis for $M$ and let $\sigma \in \SG_n$.
For each $b, b' \in B$, let $p^{\sigma}_{b,b'}(v^{1/2})$ be the  coefficient of $b'$ in the expansion of $T_\sigma \cdot b$ with respect to  $B$.
According to \cref{rem: poly action and hecke action under free cond}, we have 
\begin{align*}
\chi_{\calF_{\iota}(M)}(\widetilde{\iota}(T_\sigma)) = \sum_{b \in B} p^{\sigma}_{b,b}(v^{1/2})\quad\text{and}\quad \chi_{\calF_{\theta_q}(M)}(\widetilde{\theta_q}(T_\sigma)) = \sum_{b \in B} p^\sigma_{b,b}(v^{1/2})|_{v^{1/2} = q}.
\end{align*}
This proves the assertion.
\end{proof}

\section{Relationship between  generic Hecke algebra modules and 0-Hecke algebra modules}\label{sec: 0-Hecke modules from specializations}

In this section, we study the relationship between two families of generic Hecke algebra modules and  those of $0$-Hecke algebra modules whose images under $\ch$ are symmetric.
The former are skew Specht modules and Kazhdan--Lusztig cell modules, while the latter will be introduced in \cref{subsec: two kinds of notable 0 Hecke modules}.
To do this, we apply the specialization functors $\calF_\iota$ and $\calF_{\theta_0}$ to two families of $\polyH$-modules.
Every module $M$ in these families is free over $\C[v^{1/2}]$, thus $\calF_\iota(M)$ and $\calF_{\theta_0}(M)$ have the same dimension and their Hecke algebra actions are related as in \cref{rem: poly action and hecke action under free cond}.

Hereafter, for the $0$-Hecke algebra $\compH{0}$, we denote the generators $T_i$ ($1 \leq i \leq n-1$) and the elements $T_\sigma$ ($\sigma \in \SG_n$) by $\opi_i$ and $\opi_\sigma$, respectively.

\subsection{Two families of 0-Hecke algebra modules whose images under \texorpdfstring{$\ch$}{ch} are symmetric}\label{subsec: two kinds of notable 0 Hecke modules}

In \cite{72Stanley}, Stanley  introduced \emph{Schur labeled skew shape posets} and conjectured that they are precisely the posets $P$ whose $P$-partition generating function $K_P$ is symmetric.
On the other hand, Bj\"orner--Wachs~\cite{91BW} introduced \emph{regular posets}  in their study of convex subsets of $\SG_{n}$ under the weak Bruhat order. 
We call a poset that is both regular and Schur labeled skew shape a \emph{regular Schur labeled skew shape poset}.
For the precise definitions of Schur labeled skew shape posets, regular posets, and regular Schur labeled skew shape posets, see \cite[Section~2.3]{23KLO}.
In this subsection, we introduce two families of $0$-Hecke modules arising from regular Schur labeled skew shape posets.

Duchamp--Hivert--Thibon~\cite{02DHT} associated to each poset $P$ a right $0$-Hecke module whose image under quasisymmetric characteristic is the $P$-partition generating function $K_P$.
For technical reasons, we introduce its left module analogue.
Let $\sfP_n$ denote the set of all posets with underlying set $[n] := \{1,2,\ldots,n\}$.
Given $P \in \sfP_n$ with the partial order $\preceq_P$, 
let
\[
\Sigma_{L}(P) := \{\sigma \in \SG_n \mid \text{$\sigma(i) \leq \sigma(j)$ for all $i, j\in [n]$ with $i \preceq_P j$}\},
\]
which can be viewed as the set of  linear extensions of $P$ (see \cite[Subsection 2.3]{23KLO}).

\begin{definition}\label{def: poset modules}
Let $P \in \sfP_n$.
Define $M_P$ to be the $\compH{0}$-module with $\C \Sigma_{L}(P)$ as the underlying space and with the $\compH{0}$-action defined by
\begin{align*}
\opi_{i} \cdot \gamma :=
\begin{cases}
-\gamma & \text{if $i \in \Des_L(\gamma)$}, \\
0 & \text{if $i \notin \Des_L(\gamma)$ and $s_i\gamma \notin \Sigma_{L}(P)$,} \\
s_i \gamma & \text{if $i \notin \Des_L(\gamma)$ and $s_i\gamma \in \Sigma_{L}(P)$}.
\end{cases} 
\end{align*}
We call $M_P$ the \emph{poset module} associated with $P$.
\end{definition}

Let $\sfRSP_n$ denote the set of regular Schur labeled skew shape posets in $\sfP_n$.
The $\compH{0}$-modules $M_P$
for $P \in \sfRSP_n$
were intensively studied in \cite{23KLO}.
\footnote{It should be noted that $M_P$ is slightly different from the $\compH{0}$-module $\sfM_P$ defined in \cite[Definition 2.8]{23KLO}. 
More precisely, $M_P$ is the $\autotheta$-twist of $\sfM_P$, where $\autotheta$ is the automorphism of $\compH{0}$ given by $\pi_i \mapsto -\opi_i$.
For instance, see \cite[Theorem 3.6(a)]{23CKO}.
}
The families of $0$-Hecke modules considered here arise naturally from this context. 

The first family is obtained from the classification, up to isomorphism, of the $\compH{0}$-modules $M_P$ as $P$ ranges over $\sfRSP_n$, which was established in \cite[Theorem 5.5]{23KLO}. 
The authors constructed $\compH{0}$-modules admitting tableau bases that provide representatives for the isomorphism classes of $M_P$ ($P \in \sfRSP_n$) in the category $\compH{0}$-$\module$.
For a standard Young tableau $\ttS$ of skew shape, let 
\[
\bfDes(\ttS) := \{1 \leq i \leq n-1 \mid \text{$i$ is strictly left of $i+1$ in $\ttS$} \}.
\]
\begin{definition}{\rm (\cite[Appendix A]{23KLO})}\label{def: def of X lam mu} 
Let $\lambda/\mu$ be a skew partition of  $n$.
Define $X_{\lambda/\mu}$ to be the $\compH{0}$-module whose underlying vector space is $\C\SYT(\lambda/\mu)$ and whose action is given as follows: for $1 \leq i \leq n-1$ and $\ttS \in \SYT(\lambda/\mu)$,
\begin{align}\label{eq: left action for X}
\opi_i \cdot \ttS =
\begin{cases}
-\,\ttS
& \text{if $i \in \bfDes(\ttS)$,}\\[2pt]
~0
& \text{if $i \notin \bfDes(\ttS)$ and $s_i \cdot \ttS \notin \SYT(\lambda/\mu)$,}\\[2pt]
s_i \cdot \ttS
& \text{if $i \notin \bfDes(\ttS)$ and $s_i \cdot \ttS \in \SYT(\lambda/\mu)$.}
\end{cases}
\end{align} 
Here, $s_i \cdot \ttS$ denotes the tableau obtained from $\ttS$ by swapping $i$ and $i+1$.   
We call $X_{\lambda/\mu}$ the \emph{$0$-Hecke skew Specht module attached to $\lambda/\mu$}.
\end{definition}

It was shown in \cite{23KLO} that   $\{X_{\lambda/\mu} \mid \lambda/\mu \text{ is a skew partition of } n\}$ forms a complete set of representatives for the isomorphism classes of $M_P$ ($P \in \sfRSP_n$) and that
\begin{align}\label{eq: ch of X}
\ch([X_{\lambda/\mu}]) = s_{\lambda/\mu}.
\end{align}

The second family arises from filtrations of $M_P$'s for $P \in \sfRSP_n$ giving a representation-theoretic interpretation of the Schur-positive expansion of skew Schur functions.
To describe these filtrations, we recall several notions related to the Robinson--Schensted correspondence.
For $\sigma \in \SG_n$, we denote by $(\Ptab(\sigma), \Qtab(\sigma))$ the image of $\sigma$ under the \emph{Robinson–Schensted correspondence}, which is  a one-to-one correspondence from $\SG_n$ to $\bigcup_{\lambda \vdash n} \SYT(\lambda) \times \SYT(\lambda)$.
We call $\Ptab(\sigma)$ and $\Qtab(\sigma)$ the \emph{insertion tableau} and \emph{recording tableau} of $\sigma$, respectively.
The \emph{dual Knuth equivalence} $\dKeq$ is the equivalence relation on $\SG_n$ defined by
\[
\sigma \dKeq \rho \quad\text{if and only if}\quad \Qtab(\sigma) = \Qtab(\rho).
\]
Let $\SYT_n := \bigcup_{\lambda \vdash n} \SYT(\lambda)$.
For each $Q \in \SYT_n$, define
\[
\Gamma_Q := \{\sigma \in \SG_n \mid \Qtab(\sigma) = Q \}
\]
and call it the \emph{dual plactic class} of $Q$.
For more information, see \cite{05BB, 91Sag, 99Stanley}.

Let $P \in \sfRSP_n$ and 
let \begin{equation}\label{eq: dist filt for M P}
0 = M_0 \subset M_1 \subset M_2 \subset \cdots \subset M_l = M_P
\end{equation}
be a filtration of $M_P$ constructed as in the proof of \cite[Theorem~6.7]{23KLO}. 
According to the construction, each successive quotient  $M_i/M_{i-1}$ has a basis $\{\gamma + M_{i-1} \mid \gamma \in \Gamma_{Q^{(i)}}\}$ for a unique $Q^{(i)} \in \SYT_n$.
Moreover,
\begin{align}\label{eq: character Mi Mi-1 Schur}
\ch([M_i/M_{i-1}]) = s_{\sh(Q^{(i)})'}.
\end{align}
To study these successive quotients, we introduce the following $\compH{0}$-modules.

\begin{definition}\label{def: def of X Q}
For each $Q \in \SYT_n$, we define $\scrD_Q$ to be the $\compH{0}$-module whose underlying vector
space is $\C\Gamma_Q$ and whose action\footnote{
The well-definedness of this action follows from the proof of \cite[Theorem 6.7]{23KLO}. }
is given as follows: for $1 \leq i \leq n-1$ and $\gamma \in \Gamma_Q$,
\begin{align}\label{eq: action on DQ}
\opi_i \cdot \gamma =
\begin{cases}
- \gamma & \text{if $i \in \Des_L(\gamma)$,}\\
0 & \text{if $i  \notin \Des_L(\gamma)$ and $s_i \gamma \notin \Gamma_Q$,}\\
s_i \gamma & \text{if $i \notin \Des_L(\gamma)$ and $s_i \gamma \in \Gamma_Q$.}
\end{cases}
\end{align}
We call  $\scrD_Q$ the  \emph{dual plactic class module attached to $Q$}.
\end{definition}

Every successive quotient of the filtration in \cref{eq: dist filt for M P} is isomorphic to a dual plactic class module.
To be precise, the $\C$-linear map 
\begin{align}\label{eq: isom btw filt and scrD}
M_i/M_{i-1} \ra \scrD_{Q^{(i)}}, \quad\gamma + M_{i-1} \mapsto \gamma \quad\text{for $\gamma \in \Gamma_{Q^{(i)}}$}
\end{align}
is an $\compH{0}$-module isomorphism.

Conversely, every dual plactic class module of $\compH{0}$ appears as a successive quotient in a filtration constructed in the proof of \cite[Theorem~6.7]{23KLO}. Indeed, a particularly simple choice of $P \in \sfRSP_n$ covers all these modules at once.
Let $P_o \in \sfP_n$ be the poset in which no two distinct elements are comparable.
Then $P_o$ belongs to $\sfRSP_n$ and 
$M_{P_o} \cong \compH{0}$.
Following the proof of \cite[Theorem~6.7]{23KLO}, we construct a filtration of $M_{P_o}$.
Let $\trianglelefteq$ be the dominance order on the set of partitions of  $n$, and write $\SYT_n= \{Q_1, Q_2, \ldots, Q_l\}$ 
in an order that
if $\sh(Q_i) \trianglelefteq \sh(Q_j)$ and $\sh(Q_i) \ne \sh(Q_j)$, then $i \le j$.
For $1 \le i \le l$, set $M_i$ to be the $\C$-span of $\bigcup_{j=1}^i\Gamma_{Q_j}$.
Then $M_i$ closed under the $\compH{0}$-action, thus we have a filtration
\begin{align}\label{eq: regular rep of 0 Hecke and its dist filt}
0=M_0 \subset M_1 \subset M_2 \subset \cdots \subset M_l= M_{P_o}.
\end{align}
By \cref{eq: isom btw filt and scrD},
for each $Q \in \SYT_n$, $\scrD_{Q} \cong M_i/M_{i-1}$ if $Q = Q_i$.

Let $\lambda \vdash n$.
Combining the above isomorphism with the property \cref{eq: character Mi Mi-1 Schur}, we derive that for all $Q \in \SYT(\lambda')$,
\begin{align}\label{eq: ch of D Q}
\ch([\scrD_{Q}]) = s_{\lambda}.
\end{align}

From the discussion above  it follows that 
\begin{align}\label{eq: MP set in DQ set}
\{M_P \mid P \in \sfRSP_n \text{ such that $\ch([M_P]) = s_{\lambda}$} \} \subseteq
\{\scrD_Q \mid Q \in \SYT(\lambda') \}. 
\end{align}
By \cite[Theorem 2.2]{05W} together with \cref{eq: ch of X}, up to isomorphism, $X_{\lambda}$ and $X_{\lambda^\circ}$ are the only $0$-Hecke skew Specht modules in   $\{\scrD_Q \mid Q \in \SYT(\lambda')\}$.
Indeed, there are many $\scrD_Q$'s that are not isomorphic to $X_{\lambda}$ or $X_{\lambda^\circ}$ (see \cref{rem: new remark for theorem on Hom}).

\subsection{Skew Specht modules and $0$-Hecke skew Specht modules}\label{subsec: 0-Hecke module from specialization 1}

In the case of symmetric groups, skew analogues of Young symmetrizers and Specht modules were well studied in \cite{79PJ, 89GW}. However, to the best of our knowledge, such analogues have not been developed for Hecke algebras. Therefore, we begin by introducing them for the Hecke algebra $H_{R;a}(\SG_n)$, where $R$ is a PID and $a \in R$.

Let $\lambda/\mu$ be a skew partition of  $n$.
Define $\alpha(\lambda/\mu)$ by the composition 
\[
\alpha(\lambda/\mu) := (\alpha_1, \ldots, \alpha_{\ell(\lambda)}),
\]
where $\alpha_i$ is the number of boxes in the $i$th row of $\tyd(\lambda/\mu)$.

\begin{definition}\label{def: skew Young sym and skew Specht}
Let $\lambda/\mu$ be a skew partition of  $n$.
\begin{enumerate}
\item The \emph{Young symmetrizer 
associated to $\lambda/\mu$} in $H_{R;a}(\SG_n)$ is defined by 
\begin{align*}
\bfz_{\lambda/\mu}(R;a) := \bfy_{\alpha(\lambda/\mu)} T_{\wcol(\ttT_{\lambda/\mu})} \bfx_{\alpha(\lambda'/\mu')}.
\end{align*}
\item
The \emph{skew Specht module  of $H_{R;a}(\SG_n)$ associated to $\lambda/\mu$} is defined by 
\begin{align*}
S^{\lambda/\mu}(R;a) := H_{R;a}(\SG_n)\bfz_{\lambda/\mu}(R;a).
\end{align*}
\end{enumerate}
\end{definition}

For the notation in \cref{def: skew Young sym and skew Specht}, refer to \cref{sec: Prelim}.
When $R$ and $a$ are clear from the context, we simply write $\bfz_{\lambda/\mu}$ and $S^{\lambda/\mu}$ for $\bfz_{\lambda/\mu}(R;a)$ and $S^{\lambda/\mu}(R;a)$, respectively.

In this subsection, we apply the specialization functors  $\calF_\iota$, $\calF_{\theta_q}$, and $\calF_{\theta_0}$ to 
$$S^{\lambda/\mu}(\C[v^{1/2}]; v),$$ and provide a relationship among the Hecke algebra modules $S^{\lambda/\mu}(\C(v^{1/2});v)$, $S^{\lambda/\mu}(\C; q^2)$, and $X_{\lambda/\mu}$.
Note that  $S^{\lambda/\mu}(\C[v^{1/2}]; v)$ is free over $\C[v^{1/2}]$ since it is a submodule of the regular representation  of $H_{\C[v^{1/2}];v}(\SG_n)$.

The following theorem is the main result of this subsection.
\begin{theorem}\label{thm: specialization of skew Specht module}
Let $\lambda/\mu$ be a skew partition of $n$.
\begin{enumerate}[label = {\rm (\arabic*)}]
\item $\calF_\iota(S^{\lambda/\mu}(\C[v^{1/2}];v)) \cong S^{\lambda/\mu}(\C(v^{1/2});v)$ as $\genericH$-modules.
\item 
For nonzero $q \in \C$, $\calF_{\theta_q}(S^{\lambda/\mu}(\C[v^{1/2}];v)) \cong S^{\lambda/\mu}(\C;q^2)$  
as $\compH{q^2}$-modules.
\item 
$\calF_{\theta_0}(S^{\lambda/\mu}(\C[v^{1/2}];v)) \cong X_{\lambda/\mu}$ as $\compH{0}$-modules.
\end{enumerate}  
\end{theorem}

\begin{remark}\label{rem: Specht and X}
Unlike \cref{thm: specialization of skew Specht module}(1) and (2), $\calF_{\theta_0}(S^{\lambda/\mu}(\C[v^{1/2}];v))$ is not necessarily isomorphic to $S^{\lambda/\mu}(\C;0)$.  
For instance, let $\lambda = (3,2)$.
Then 
\[
\dim_{\C} \calF_{\theta_0}(S^{\lambda/\mu}(\C[v^{1/2}];v)) =
\dim_{\C} X_{\lambda} = 5. 
\]
On the other hand, we have $\dim S^{\lambda}(\C;0) = 2$.
To see this, note that 
\[
S^\lambda(\C;0) = 
H_{\C;0}(\SG_5)
\opi_{32154}\opi_{14253}
\pi_{1} \pi_3,
\]
where $\pi_i = \opi_i + 1$ for $1 \le i \le 4$.
By \cite[Theorem 1(2)]{22JKLO}, we have an $H_{\C;0}(\SG_5)$-module isomorphism
\[
H_{\C;0}(\SG_5)
\opi_{32154}\opi_{14253}\pi_{1} \pi_3 \cong \overline{\sfB}(35241,45231),
\]
where $\overline{\sfB}(35241,45231)$ is a negative weak Bruhat interval module introduced in \cite[Definition 1]{22JKLO}.
It is straightforward that $\dim_{\C} \overline{\sfB}(35241,45231) = 2$.
Therefore, 
\[
\calF_{\theta_0}(S^{\lambda}(\C[v^{1/2}];v)) \not\cong S^\lambda(\C;0).
\]
\end{remark}

Before proving this theorem, we  introduce definitions, notation, and lemmas that are necessary for the proof.

First, we introduce a subset of $S^{\lambda/\mu}(\C[v^{1/2}];v)$, which will later be shown to form a basis for $S^{\lambda/\mu}(\C[v^{1/2}];v)$.
A bijective tableau $\ttT$ is called \emph{row-standard} if the entries in each row of $\ttT$ increase from left to right.
Denote by $\RT(\lambda/\mu)$ the set of all row-standard  tableaux of shape $\lambda/\mu$.
For $\ttT \in \RT(\lambda/\mu)$, define 
$\sfh_{\ttT} \in \polyH$  by
\[
\sfh_{\ttT} := T_{\wrow(\ttT)} \bfz_{\lambda/\mu}.
\]
Then, the desired subset is defined by 
\[
B(\lambda/\mu) := \{\sfh_{\ttS} \mid \ttS \in \SYT(\lambda/\mu) \}.
\]

In \cref{lem: triangularity of elt in  B}, we analyze the expansion of each element in $B(\lambda/\mu)$ in the standard basis $\{T_\sigma \mid \sigma \in \SG_n\}$.
The following lemma is immediate from this analysis.

\begin{lemma}\label{lem: lin indep of B(lam/mu)}
Let $\lambda/\mu$ be a skew partition of $n$.
Then, $B(\lambda/\mu)$ is $\C[v^{1/2}]$-linearly independent.
\end{lemma}

Second, we present two crucial lemmas for the proof of the main theorem.
To this end, we begin with introducing  the Garnir element corresponding to a cell of $\tyd(\lambda/\mu)$ and establish its property.
Let $(i,j)\in \tyd(\lambda/\mu)$ be a box such that $(i+1,j)\in \tyd(\lambda/\mu)$.
We define $\SG^{\lambda/\mu}_{i,j}$ to be the set of permutations $\wrow(\ttT)$, where $\ttT$ ranges over all row-standard tableaux of shape $\lambda/\mu$ that agree with $\ttT_{\lambda/\mu}$ at every box $(r,c)\in \tyd(\lambda/\mu)$ except when
\begin{align}\label{eq: Garnir belt}
\text{(i) } r=i \text{ and } c\ge j,
\qquad
\text{(ii) } r=i+1 \text{ and } c\le j.
\end{align}
Equivalently, 
\begin{align}\label{eq: Garnir permutations alternating def}
\SG^{\lambda/\mu}_{i,j}
= 
\calD^{\{\ttT_{\lambda / \mu}(i,\lambda_i)\}}
\cap 
\left\langle s_i \mid \ttT_{\lambda / \mu}(i,j) \le i < \ttT_{\lambda / \mu}(i+1,j) \right\rangle,
\end{align}
where $\calD^{\{\ttT_{\lambda / \mu}(i,\lambda_i)\}}$ denotes the set of permutations $\sigma \in \SG_n$ with $\Des_R(\sigma) \subseteq \{\ttT_{\lambda / \mu}(i,\lambda_i)\}$.
We then define the \emph{Garnir element corresponding to $(i,j)$} to be the element $G^{\lambda/\mu}_{i,j}$ in $\polyH$ given by
\[
G^{\lambda/\mu}_{i,j} := \bfy_{\SG^{\lambda/\mu}_{i,j}}.
\]
\begin{example}
Let $\lambda / \mu = (3,3,2) / (1)$ and $(i,j) = (1,2)$.
Then
\[
\ttT_{(3,3,2) / (1)} =
 \begin{array}{l}
\begin{ytableau}
\none & *(yellow!50) 1 & *(yellow!50) 2 \\
*(yellow!50) 3 & *(yellow!50) 4 & 5  \\
6 & 7 
\end{ytableau}
\end{array}.
\]
Here, the boxes $(r,c)$ highlighted in yellow are those satisfying either
(i) $r=i$ and $c\ge j$, or
(ii) $r=i+1$ and $c\le j$.
We see that $\SG^{(3,3,2)/(1)}_{1,2}$ consists of $\wrow(\ttT)$, where $\ttT$ ranges over the following tableaux:
\[
\begin{array}{l}
\begin{ytableau}
\none & *(yellow!50) 1  & *(yellow!50) 2 \\
*(yellow!50) 3 & *(yellow!50) 4 & 5  \\
6 & 7 
\end{ytableau}
\end{array}
\quad
 \begin{array}{l}
\begin{ytableau}
\none & *(yellow!50) 1 & *(yellow!50) 3 \\
*(yellow!50) 2 & *(yellow!50) 4 & 5  \\
6 & 7 
\end{ytableau}
\end{array}
\quad
\begin{array}{l}
\begin{ytableau}
\none & *(yellow!50) 1 & *(yellow!50) 4 \\
*(yellow!50) 2 & *(yellow!50) 3 & 5  \\
6 & 7 
\end{ytableau}
\end{array}
\quad
\begin{array}{l}
\begin{ytableau}
\none & *(yellow!50) 2 & *(yellow!50) 3 \\
*(yellow!50) 1 & *(yellow!50) 4 &  5  \\
6 & 7 
\end{ytableau}
\end{array}
\begin{array}{l}
\begin{ytableau}
\none & *(yellow!50) 2 & *(yellow!50) 4 \\
*(yellow!50) 1 & *(yellow!50) 3 & 5  \\
6 & 7 
\end{ytableau}
\end{array}
\quad
\begin{array}{l}
\begin{ytableau}
\none & *(yellow!50) 3 & *(yellow!50) 4 \\
*(yellow!50) 1 & *(yellow!50) 2 & 5  \\
6 & 7 
\end{ytableau}
\end{array}
\]
Therefore,
\begin{align*}
\SG^{(3,3,2) / (1)}_{1,2} 
& = \left\{
1234567, 1324567, 1423567, 2314567, 2413567, 3412567
\right\}, 
\end{align*}
and the Garnir element $G^{(3,3,2) / (1)}_{1,2}$ is given by 
\begin{align*} 
v^4 T_{\id} - v^3 T_{s_2} + v^2 T_{s_3s_2} + v^2 T_{s_1s_2} - v T_{s_1s_3s_2} + T_{s_2s_1s_3s_2}.
\end{align*}
\end{example}

The following lemma shows that the  Garnir element corresponding to a cell of $\tyd(\lambda/\mu)$  annihilates the Young symmetrizer $\bfz_{\lambda/\mu}$.
To avoid interrupting the flow of the main argument, we defer its proof to \cref{lem: Garnir rel}.

\begin{lemma}\label{lem: Garnir rel sec 4}
Let $\lambda/\mu$ be a skew partition of $n$ and $(i, j) \in \tyd(\lambda/\mu)$.
Then,
\begin{equation*}
G^{\lambda/\mu}_{i, j} \bfz_{\lambda/\mu}(\C[v^{1/2}];v) = 0.
\end{equation*}
\end{lemma}

\begin{remark}
A relation analogous to \cref{lem: Garnir rel sec 4} is stated in \cite{92KW} for ordinary partitions. However, that statement is incorrect. Indeed, let $\lambda=(2,2)$. Using the notation of \cite{92KW}, a direct computation gives
\begin{align*}
Y^{\lambda}(q)
&= (1 + T_{212} + T_{232} + T_{2132})
   \cdot (q^2 - qT_1 - qT_3 + T_{13}), \\
G^\lambda_{2,2}(q)
&= q^2 - qT_2 + T_{23}.
\end{align*}
It then follows by direct computation that
\[
Y^{\lambda}(q)G^\lambda_{2,2}(q)\neq 0.
\]
Thus, the relation asserted in \cite{92KW} fails even for the ordinary partition $\lambda=(2,2)$. To the best of our knowledge, no analogue of \cref{lem: Garnir rel sec 4} for general skew shapes has appeared in the literature.
\end{remark}

For $B \in \tyd(\lambda/\mu)$, let 
\begin{equation}\label{def of SYTB}
\SYT(\lambda/\mu)^B_{\boxn} := 
\{\ttS \in \SYT(\lambda/\mu) \mid \text{$n$ lies weakly below $B$ in $\ttS$}  \}.
\end{equation}
We are now ready to introduce two crucial lemmas.

\begin{lemma}\label{lem: expansion of row-standard into standard basis}
Let $\ttT \in \RT(\lambda/\mu)$, and set $B := \ttT^{-1}(n)$.
Then $\sfh_{\ttT}$ admits a unique $\Z[v]$-linear expansion 
\begin{align}\label{eq: row standard h expansion}
\sfh_{\ttT}
=
\sum_{\ttS \in \SYT(\lambda/\mu)^{B}_{\boxn}}
f_{\ttT,\ttS}(v)\,\sfh_{\ttS}
\end{align}
in 
$\{\sfh_{\ttS} \mid \ttS \in \SYT(\lambda/\mu)^{B}_{\boxn}\}.$
Moreover, if $\ttT \notin \SYT(\lambda/\mu)$, then 
each coefficient in this expansion is divisible by $v$.
\end{lemma}

\begin{proof}
We prove the assertion by using induction on $\ell(\wrow(\ttT))$.

If $\ell(\wrow(\ttT))=0$, then $\ttT = \ttT_{\lambda/\mu}$, and so $\sfh_{\ttT} = \sfh_{\ttT_{\lambda/\mu}}$.
Due to \cref{lem: lin indep of B(lam/mu)}, 
this is the unique expansion of the form \cref{eq: row standard h expansion}.

Now, suppose that $\ell(\wrow(\ttT)) > 0$ and that the assertion holds for $\ttT' \in \RT(\lambda/\mu)$ with $\ell(\wrow(\ttT')) < \ell(\wrow(\ttT))$.
If $\ttT \in \SYT(\lambda/\mu)$, then the assertion can be shown as in the case where $\ell(\wrow(\ttT)) = 0$.
Consider the case where $\ttT \notin \SYT(\lambda/\mu)$.
Since $\ttT \in \RT(\lambda/\mu) \setminus \SYT(\lambda/\mu)$,
we can choose $(i, j) \in \tyd(\lambda/\mu)$ such that 
\[
(i+1, j) \in \tyd(\lambda/\mu) \quad \text{and} \quad \ttT(i, j) > \ttT(i+1, j).
\]
Let $D$ be the set of boxes $(r,c) \in \tyd(\lambda/\mu)$ satisfying \cref{eq: Garnir belt}.
Enumerate the elements of $\{ \ttT(r,c) \mid (r,c) \in D \}$ in increasing order as
\[
x_1 < x_2 < \cdots < x_{\lambda_i - \mu_{i+1} + 1}.
\]
For the reader's convenience, we illustrate below the $i$th and $(i+1)$st rows of $\ttT$, with the boxes in $D$ highlighted in red and their entries indicated.
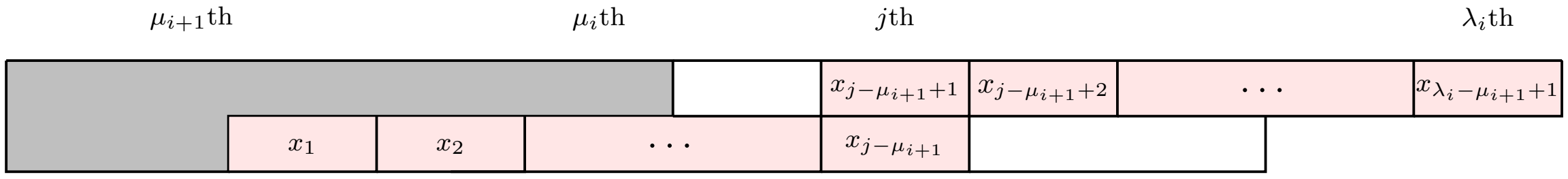
\[
\scalebox{0.9}{$
\def \hd {1.6}
\def \vd {0.6}
\begin{tikzpicture}

\draw[thick, fill=gray!50] (-0.5*\hd, 2*\vd) -- (4*\hd, 2*\vd) -- (4*\hd, 1*\vd) -- (1*\hd, 1*\vd) -- (1*\hd, 0*\vd) -- (-0.5*\hd, 0*\vd) -- (-0.5*\hd, 2*\vd);

\node at (0.75*\hd, 2.75*\vd) {\tiny $\mu_{i+1}$th};
\node at (3.5*\hd, 2.75*\vd) {\tiny $\mu_i$th};
\node at (5.5*\hd, 2.75*\vd) {\tiny $j$th};
\node at (9.5*\hd, 2.75*\vd) {\tiny $\lambda_{i}$th};

\draw[fill = red!10] (1*\hd, 0*\vd) -- (6*\hd, 0*\vd) -- (6*\hd, 1*\vd) -- (10*\hd, 1*\vd) -- (10*\hd, 2*\vd) -- (5*\hd, 2*\vd) -- (5*\hd, 1*\vd) -- (1*\hd, 1*\vd) -- (1*\hd, 0*\vd);

\draw[thick] (1*\hd, 0*\vd) -- (2*\hd, 0*\vd) -- (2*\hd, 1*\vd);
\node at (1.5*\hd, 0.45*\vd) {\tiny $x_{1}$};

\draw[thick] (2*\hd, 0*\vd) -- (3*\hd, 0*\vd) -- (3*\hd, 1*\vd);
\node at (2.5*\hd, 0.45*\vd) {\tiny $x_{2}$};

\draw[thick] (2.5*\hd, 0*\vd) -- (5*\hd, 0*\vd) -- (5*\hd, 1*\vd);
\node at (4*\hd, 0.45*\vd) {$\cdots$};

\draw[thick] (5*\hd, 0*\vd) -- (6*\hd, 0*\vd) -- (6*\hd, 1*\vd);
\node at (5.5*\hd, 0.45*\vd) {\tiny $x_{j - \mu_{i+1}}$};

\draw[thick] (6*\hd, 0*\vd) -- (8*\hd, 0*\vd) -- (8*\hd, 1*\vd);

%%%%%%%%%%%%
\draw[thick] (4*\hd, 2*\vd) -- (10*\hd, 2*\vd);

\draw[thick] (4*\hd, 1*\vd) -- (5*\hd, 1*\vd) -- (5*\hd, 2*\vd);
% \node at (5.5*\hd, 1.45*\vd) {$*$};

\draw[thick] (5*\hd, 1*\vd) -- (6*\hd, 1*\vd) -- (6*\hd, 2*\vd);
\node at (5.5*\hd, 1.45*\vd) {\tiny $x_{j-\mu_{i+1}+1}$};

\draw[thick] (6*\hd, 1*\vd) -- (7*\hd, 1*\vd) -- (7*\hd, 2*\vd);
\node at (6.5*\hd, 1.45*\vd) {\scriptsize $x_{j-\mu_{i+1}+2}$};

\draw[thick] (7*\hd, 1*\vd) -- (9*\hd, 1*\vd) -- (9*\hd, 2*\vd);
\node at (8*\hd, 1.45*\vd) {$\cdots$};

\draw[thick] (9*\hd, 1*\vd) -- (10*\hd, 1*\vd) -- (10*\hd, 2*\vd);
\node at (9.5*\hd, 1.45*\vd) 
{\tiny $x_{\lambda_i-\mu_{i+1}+1}$};

\end{tikzpicture}
$}
\]

Let $\mathtt{U}_0$ be the bijective tableau of shape $\lambda/\mu$ that agrees with $\ttT$ on every box in $\tyd(\lambda/\mu)\setminus D$ and whose boxes in $D$ are filled with
$x_1, x_2, \ldots, x_{\lambda_i - \mu_{i+1} + 1}$
from left to right, first along the $i$th row and then along the $(i+1)$st row.
We illustrate below the $i$th and $(i+1)$st rows of $\mathtt{U}_0$, with the boxes in $D$ highlighted in red and their entries indicated.
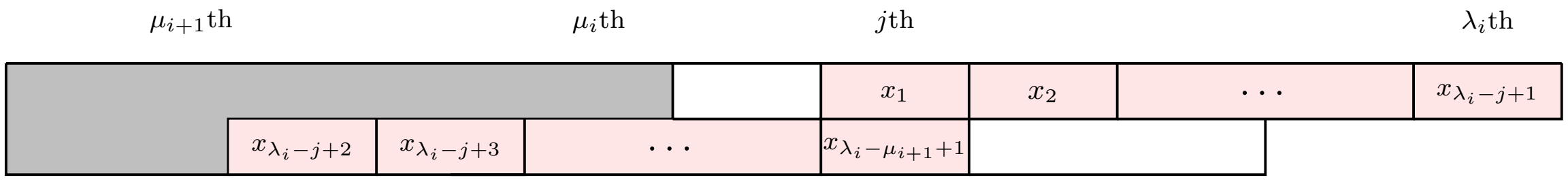
\[
\scalebox{0.9}{$
\def \hd {1.6}
\def \vd {0.6}
\begin{tikzpicture}

\draw[thick, fill=gray!50] (-0.5*\hd, 2*\vd) -- (4*\hd, 2*\vd) -- (4*\hd, 1*\vd) -- (1*\hd, 1*\vd) -- (1*\hd, 0*\vd) -- (-0.5*\hd, 0*\vd) -- (-0.5*\hd, 2*\vd);
% \draw[thick, fill=gray!50] (0*\hd, 1*\vd) -- (1.5*\hd, 1*\vd)

% \node at (-0.5*\hd, 1.5*\vd) {\tiny $i$th};
% \node at (-0.5*\hd, 0.5*\vd) {\tiny $(i+1)$st};

\node at (0.75*\hd, 2.75*\vd) {\tiny $\mu_{i+1}$th};
\node at (3.5*\hd, 2.75*\vd) {\tiny $\mu_i$th};
\node at (5.5*\hd, 2.75*\vd) {\tiny $j$th};
\node at (9.5*\hd, 2.75*\vd) {\tiny $\lambda_{i}$th};

\draw[fill = red!10] (1*\hd, 0*\vd) -- (6*\hd, 0*\vd) -- (6*\hd, 1*\vd) -- (10*\hd, 1*\vd) -- (10*\hd, 2*\vd) -- (5*\hd, 2*\vd) -- (5*\hd, 1*\vd) -- (1*\hd, 1*\vd) -- (1*\hd, 0*\vd);

\draw[thick] (1*\hd, 0*\vd) -- (2*\hd, 0*\vd) -- (2*\hd, 1*\vd);
\node at (1.5*\hd, 0.45*\vd) {\tiny $x_{\lambda_i - j + 2}$};

\draw[thick] (2*\hd, 0*\vd) -- (3*\hd, 0*\vd) -- (3*\hd, 1*\vd);
\node at (2.5*\hd, 0.45*\vd) {\tiny $x_{\lambda_i - j + 3}$};

% \draw[thick] (2.5*\hd, 0*\vd) -- (3*\hd, 0*\vd) -- (3*\hd, 1*\vd);
% \node at (2.75*\hd, 0.45*\vd) {\tiny $x_2$};

\draw[thick] (2.5*\hd, 0*\vd) -- (5*\hd, 0*\vd) -- (5*\hd, 1*\vd);
\node at (4*\hd, 0.45*\vd) {$\cdots$};

\draw[thick] (5*\hd, 0*\vd) -- (6*\hd, 0*\vd) -- (6*\hd, 1*\vd);
\node at (5.5*\hd, 0.45*\vd) {\tiny $x_{\lambda_i - \mu_{i+1} +1}$};

\draw[thick] (6*\hd, 0*\vd) -- (8*\hd, 0*\vd) -- (8*\hd, 1*\vd);

%%%%%%%%%%%%
\draw[thick] (4*\hd, 2*\vd) -- (10*\hd, 2*\vd);

\draw[thick] (4*\hd, 1*\vd) -- (5*\hd, 1*\vd) -- (5*\hd, 2*\vd);
% \node at (5.5*\hd, 1.45*\vd) {$*$};

\draw[thick] (5*\hd, 1*\vd) -- (6*\hd, 1*\vd) -- (6*\hd, 2*\vd);
\node at (5.5*\hd, 1.45*\vd) {\tiny $x_{1}$};

\draw[thick] (6*\hd, 1*\vd) -- (7*\hd, 1*\vd) -- (7*\hd, 2*\vd);
\node at (6.5*\hd, 1.45*\vd) {\scriptsize $x_{2}$};

\draw[thick] (7*\hd, 1*\vd) -- (9*\hd, 1*\vd) -- (9*\hd, 2*\vd);
\node at (8*\hd, 1.45*\vd) {$\cdots$};

\draw[thick] (9*\hd, 1*\vd) -- (10*\hd, 1*\vd) -- (10*\hd, 2*\vd);
\node at (9.5*\hd, 1.45*\vd) 
{\tiny $x_{\lambda_i - j +1}$};

\end{tikzpicture}
$}
\]
For each $\sigma \in \SG^{\lambda/\mu}_{i,j}$, let $\mathtt{U}_\sigma$ be the bijective tableau of shape $\lambda/\mu$ satisfying 
\begin{align*}
\wrow(\mathtt{U}_\sigma) = \wrow(\mathtt{U}_0)\sigma.
\end{align*}

Note that $\SG^{\lambda/\mu}_{i,j}$ is a left weak Bruhat interval since 
the sets on the right hand side of \cref{eq: Garnir permutations alternating def} are left weak Bruhat intervals.
Let $\sigma_0$ be the unique element of maximal length in $\SG^{\lambda/\mu}_{i,j}$.
One can easily see that
\[
\sigma_0(r) = \begin{cases}
r & \text{if $r < t_1$ or $r > t_3$,} \\
r + j - \mu_{i+1} & \text{if $t_1 \le r < t_2$,} \\
r - t_2 + t_1  & \text{if $t_2 \le k \le t_3$,}
\end{cases}
\]
where
\[
t_1 = \sum_{1 \leq k \leq i-1} (\lambda_k - \mu_k) + (j - \mu_i), \quad
t_2 = \sum_{1 \leq k \leq i} (\lambda_k - \mu_k) + 1, \quad \text{and} \quad
t_3 = \sum_{1 \leq k \leq i} (\lambda_k - \mu_k) + (j - \mu_{i+1}).
\]
From this, we deduce that $\wrow(\ttT) = \wrow(\mathtt{U}_0) \sigma_0$, and so
\begin{align}\label{eq: R = U sigma 0}
\ttT = \mathtt{U}_{\sigma_0}.
\end{align}
Moreover, for each $\sigma \in \SG^{\lambda/\mu}_{i,j}$,
combining the inequalities
$\wrow(\mathtt{U}_0)(t_1) < \wrow(U_0)(t_1 + 1) < \cdots < \wrow(\mathtt{U}_0)(t_3)$
with the fact that $\sigma$ permutes only the integers in $[t_1,\, t_3]$, we have
\begin{align}\label{eq: length preserving cond on Rd}
\ell(\wrow(\mathtt{U}_\sigma)) = \ell(\wrow(\mathtt{U}_0)\sigma)  = \ell(\wrow(\mathtt{U}_0)) + \ell(\sigma).
\end{align}
Hence
\begin{align}\label{eq: expansion of TG}
\begin{aligned}
T_{\wrow(\mathtt{U}_0)} G^{\lambda/\mu}_{i,j}
& = \sum_{\sigma \in \SG^{\lambda/\mu}_{i,j}} \sgn(\sigma) v^{\ell(\sigma_0) - \ell(\sigma)}
T_{\wrow(\mathtt{U}_0)} T_\sigma
\\
& \overset{\cref{eq: length preserving cond on Rd}}{=}  \sum_{\sigma \in \SG^{\lambda/\mu}_{i,j}} \sgn(\sigma) v^{\ell(\sigma_0) - \ell(\sigma)} T_{\wrow(\mathtt{U}_\sigma)} \\
& \overset{\cref{eq: R = U sigma 0}}{=} \sgn(\sigma_0) T_{\wrow(\ttT)} + 
\sum_{\sigma \in \SG^{\lambda/\mu}_{i,j} \setminus \{\sigma_0\}} \sgn(\sigma) v^{\ell(\sigma_0) - \ell(\sigma)} T_{\wrow(\mathtt{U}_\sigma)}.
\end{aligned}
\end{align}
By~\cref{lem: Garnir rel sec 4}, we have $T_{\wrow(\mathtt{U}_0)} G^{\lambda/\mu}_{i,j} \bfz_{\lambda/\mu} = 0$.
Thus, multiplying the first and last equalities of \cref{eq: expansion of TG} on the right by $\bfz_{\lambda/\mu}$ yields
\begin{align}\label{eq: expansion of hR}
0 = \sgn(\sigma_0)\sfh_{\ttT} + 
\sum_{\sigma \in \SG^{\lambda/\mu}_{i,j} \setminus \{\sigma_0\}} \sgn(\sigma) v^{\ell(\sigma_0) - \ell(\sigma)} T_{\wrow(\mathtt{U}_\sigma)} \bfz_{\lambda/\mu}.
\end{align}

For each $\sigma \in \SG^{\lambda/\mu}_{i,j}$, let $\ttT_\sigma$ be the row-standard tableau of shape $\lambda/\mu$ whose $r$th row contains the same entries as the $r$th row of $\mathtt{U}_\sigma$ for all $1 \le r \le \ell(\lambda)$.
By definition, there exists $\xi_\sigma \in \SG_{\alpha(\lambda/\mu)}$ such that
\[
\wrow(\mathtt{U}_\sigma) = \wrow(\ttT_\sigma)\xi_\sigma.
\]
Moreover, since $\ttT_\sigma$ is row-standard, we have 
\[
\ell(\wrow(\mathtt{U}_\sigma)) = \ell(\wrow(\ttT_\sigma)) + \ell(\xi_\sigma).
\]
It follows that for each $\sigma \in \SG^{\lambda/\mu}_{i,j}$,
\begin{align*}
T_{\wrow(\mathtt{U}_\sigma)} \bfz_{\lambda/\mu} 
& = T_{\wrow(\ttT_\sigma)} T_{\xi_\sigma} \; \bfy_{\alpha(\lambda/\mu)} T_{\wcol(\ttT_{\lambda/\mu})} \bfx_{\alpha(\lambda'/\mu')} \\
& \overset{\text{\cref{eq: Ti and bfy}}}{=} \sgn(\xi_\sigma) T_{\wrow(\ttT_\sigma)} \bfy_{\alpha(\lambda/\mu)} T_{\wcol(\ttT_{\lambda/\mu})} \bfx_{\alpha(\lambda'/\mu')} \\
& = \sgn(\xi_\sigma) T_{\wrow(\ttT_\sigma)} \bfz_{\lambda/\mu} \\
& = \sgn(\xi_\sigma) \sfh_{\ttT_\sigma}.
\end{align*}
Combining this with \cref{eq: expansion of hR}, we have 
\begin{align}\label{eq: expansion of hR to hRsigma}
\sfh_{\ttT} 
& =
-\sgn(\sigma_0) \sum_{\sigma \in \SG^{\lambda/\mu}_{i,j} \setminus \{\sigma_0\}} 
\sgn(\sigma)
\sgn(\xi_\sigma)
v^{\ell(\sigma_0) - \ell(\sigma)}  \sfh_{\ttT_\sigma}.
\end{align}

Let $\sigma \in G^{\lambda/\mu}_{i, j} \setminus \{\sigma_0\}$.
From the inequalities $\ell(\wrow(\ttT_\sigma)) \leq  \ell(\wrow(\mathtt{U}_\sigma))$ and  $ \ell(\sigma_0) - \ell(\sigma) > 0$, we derive that 
\begin{equation}\label{eq: lenth cond on row standard}
\ell(\wrow(\ttT_\sigma)) 
\leq 
\ell(\wrow(\mathtt{U}_0)) + \ell(\sigma) 
< \ell(\wrow(\mathtt{U}_0)) + \ell(\sigma_0)
\overset{\cref{eq: R = U sigma 0}}{=} \ell(\wrow(\ttT)).
\end{equation}
By the induction hypothesis, $\sfh_{\ttT_\sigma}$ is a $\Z[v]$-linear combination of $\{\sfh_{\ttS} \mid \ttS \in \SYT(\lambda/\mu)^{\ttT^{-1}_\sigma (n)}_{\boxn}\}$.
If $B \notin D$, then $\ttT_{\sigma}^{-1}(n) = B$. 
Otherwise, $B = (i,\lambda_i)$ and $\ttT_{\sigma}^{-1}(n)$ lies in either the $i$th row or $(i+1)$st row.
It follows that  $\ttT_{\sigma}^{-1}(n)$ is weakly below $B$.
Thus, by \cref{eq: expansion of hR to hRsigma}, $\sfh_{\ttT}$ is a $\Z[v^{1/2}]$-linear expansion of $\{\sfh_{\ttS} \mid \ttS \in \SYT(\lambda/\mu)^B_{\boxn}\}$.
The uniqueness of this expansion follows from \cref{lem: lin indep of B(lam/mu)}.
Finally, every coefficient in this expansion is divisible by $v$ since $\ell(\sigma_0) - \ell(\sigma) \geq 1$ on the right hand side of \cref{eq: expansion of hR to hRsigma}.
\end{proof}

\begin{lemma}\label{lem: skew specht basis}
Let $\lambda/\mu$ be a skew partition of  $n$.
Then, $B(\lambda/\mu)$ is a $\C[v^{1/2}]$-basis for $S^{\lambda/\mu}(\C[v^{1/2}];v)$.
\end{lemma}

\begin{proof}
By \cref{lem: lin indep of B(lam/mu)}, it remains to prove that
$B(\lambda/\mu)$ spans
$S^{\lambda/\mu}(\C[v^{1/2}];v)$ over $\C[v^{1/2}]$.
Since $S^{\lambda/\mu}(\C[v^{1/2}];v)$ is 
spanned by $\{T_\sigma\bfz_{\lambda/\mu} \mid \sigma\in\SG_n\}$,
it suffices to show that each element of this set belongs to the
$\C[v^{1/2}]$-span of $B(\lambda/\mu)$.

Let $\sigma \in \SG_n$.
Note that the set of minimal coset representatives for
$\SG_n/\SG_{\alpha(\lambda/\mu)}$ is $\{\wrow(\ttT) \mid \ttT \in \RT(\lambda/\mu) \}$.
It follows that there exist $\ttT_\sigma \in \RT(\lambda/\mu)$ and $\sigma' \in \SG_{\alpha(\lambda/\mu)}$ such that $\sigma = \wrow(\ttT_\sigma)  \sigma'$ and $\ell(\sigma) = \ell(\wrow(\ttT_\sigma)) +  \ell(\sigma')$.
Using this decomposition, we have
\begin{align*}
T_\sigma \bfz_{\lambda/\mu} & = T_{\wrow(\ttT_\sigma)} T_{\sigma'} \bfz_{\lambda/\mu} =T_{\wrow(\ttT_\sigma)} \sgn(\sigma') \bfz_{\lambda/\mu} = \sgn(\sigma') \sfh_{\ttT_\sigma},
\end{align*}
where the second equality follows from \cref{eq: Ti and bfy}.
By \cref{lem: expansion of row-standard into standard basis},
$\sfh_{\ttT_\sigma}$ belongs to the $\C[v^{1/2}]$-span of
$B(\lambda/\mu)$, and hence so does $T_\sigma\bfz_{\lambda/\mu}$.
\end{proof}

Now we are ready to prove  \cref{thm: specialization of skew Specht module}.

\begin{proof}[Proof of \protect\cref{thm: specialization of skew Specht module}.]

(1) By \cref{cor: algebra homomorphism vs functor sub of regular}(1), 
\[
\calF_\iota(S^{\lambda/\mu}(\C[v^{1/2}];v)) \cong S^{\lambda/\mu}(\C(v^{1/2});v) \quad\text{as $\genericH$-modules.}
\]

(2) By \cref{lem: skew specht basis},  $\widetilde{\theta_q}(B(\lambda/\mu))$ spans  $S^{\lambda/\mu}(\C;q^2)$, which is the image of $S^{\lambda/\mu}(\C[v^{1/2}];v)$ under $\widetilde{\theta_q}$.
Furthermore, since $q \neq 0$, using \cref{lem: triangularity of elt in  B}, we see that $\widetilde{\theta_q}(B(\lambda/\mu))$ is linearly independent over $\C$. 
It follows that $\dim_{\C}S^{\lambda/\mu}(\C;q^2) = |B(\lambda/\mu)|$, and thus by \cref{cor: algebra homomorphism vs functor sub of regular}(2),
\[
\calF_{\theta_q}(S^{\lambda/\mu}(\C[v^{1/2}];v)) \cong S^{\lambda/\mu}(\C;q^2) \quad\text{as $\compH{q^2}$-modules.}
\]

(3)  
Combining \cref{thm: rank of tensored theta our cons} with \cref{lem: skew specht basis}, we see that $1 \otimes B(\lambda/\mu)$ is a $\C$-basis for $\calF_{\theta_0}(S^{\lambda/\mu}(\C[v^{1/2}];v))$.
Moreover, by construction, $\SYT(\lambda/\mu)$ is a $\C$-basis for $X_{\lambda/\mu}$.
Define a $\C$-linear isomorphism $\psi: \calF_{\theta_0}(S^{\lambda/\mu}(\C[v^{1/2}];v)) \ra X_{\lambda/\mu}$ by letting
\[
\psi(1 \otimes \sfh_{\ttS}) = \ttS \quad\text{for $\ttS \in \SYT(\lambda/\mu)$}
\]
and extending linearly.
To prove our assertion, it suffices to show that $\psi$ preserves the $\compH{0}$-action.
Let $\ttS \in \SYT(\lambda/\mu)$ and $1 \leq k \leq n-1$.
We have three cases.
\smallskip

\textbf{Case 1:} $k$ and $k+1$ are in the same row of $\ttS$.  

In this case, $k$ lies immediately to the left of $k+1$ in $\ttS$.
Set $i := \wrow(\ttS)^{-1}(k)$.
Then
\[
k \notin \Des_L(\wrow(\ttS)),
\quad
s_k \wrow(\ttS) = \wrow(\ttS) s_i,
\quad \text{and} \quad 
s_i \in \SG_{\alpha(\lambda/\mu)}.
\]
It follows that 
\begin{align*}
T_k \sfh_{\ttS} 
& = T_k T_{\wrow(\ttS)} \bfz_{\lambda/\mu}
= T_{s_k\wrow(\ttS)} \bfz_{\lambda/\mu}
= T_{\wrow(\ttS)s_i} \bfz_{\lambda/\mu}
= T_{\wrow(\ttS)} T_i \bfz_{\lambda/\mu} \\
& = T_{\wrow(\ttS)} T_i \bfy_{\alpha(\lambda/\mu)} T_{\wcol(\ttT_{\lambda/\mu})} \bfx_{\alpha(\lambda'/\mu')} \\
& = - T_{\wrow(\ttS)} \bfy_{\alpha(\lambda/\mu)} T_{\wcol(\ttT_{\lambda/\mu})} \bfx_{\alpha(\lambda'/\mu')} \\
& = - \sfh_{\ttS}.
\end{align*}
Thus, $\opi_k \cdot 1 \otimes \sfh_{\ttS} 
= 1 \otimes T_k  \sfh_{\ttS}
= - 1 \otimes \sfh_{\ttS}$.
Moreover, since $\opi_k \cdot \ttS = -\ttS$, we have that
\[
\psi(\opi_k \cdot 1 \otimes \sfh_{\ttS}) 
= \psi(- 1 \otimes \sfh_{\ttS}) 
= - \ttS = \opi_k \cdot \ttS 
= \opi_k \cdot \psi(1 \otimes \sfh_{\ttS}).
\]

\textbf{Case 2:} $k$ and $k+1$ are in the same column of $\ttS$.

In this case, $k \notin \Des_L(\wrow(\ttS))$ and $s_k \cdot \ttS \in \RT(\lambda/\mu) \setminus \SYT(\lambda/\mu)$.
By \cref{lem: expansion of row-standard into standard basis},  $\sfh_{s_k \cdot \ttS}$ is the $v \Z[v]$-linear combination of $B(\lambda/\mu)$.
Combining this with $T_k \sfh_{\ttS} = \sfh_{s_k \cdot \ttS}$ yields  
\[
\opi_k \cdot 1 \otimes \sfh_{\ttS} = 1 \otimes T_k \sfh_{\ttS} = 1 \otimes \sfh_{s_k \cdot \ttS} = 0.
\]
Moreover, since  $\opi_k \cdot \ttS = 0$, we have that  
\[
\psi(\opi_k \cdot 1 \otimes \sfh_{\ttS})  
= 0 = \opi_k \cdot \ttS 
= \opi_k \cdot \psi(1 \otimes \sfh_{\ttS}).
\]

\textbf{Case 3:} $k$ and $k+1$ are in the distinct rows and  columns of $\ttS$.  

In this case, $s_k \cdot \ttS \in \SYT(\lambda/\mu)$.
Suppose that $k+1$ is below $k$. 
Then $\opi_k \cdot \ttS = s_k \cdot \ttS$ and $k \notin \Des_L(\wrow(\ttS))$.
The latter implies that
$\opi_k \cdot 1 \otimes \sfh_{\ttS} 
= 1 \otimes T_k \sfh_{\ttS} 
= 1 \otimes \sfh_{s_k \cdot \ttS}$.
Thus,
\[
\psi(\opi_k \cdot 1 \otimes \sfh_{\ttS}) 
= \psi(1 \otimes \sfh_{s_k \cdot \ttS}) 
= s_k \cdot \ttS 
= \opi_k \cdot \ttS 
= \opi_k \cdot \psi(1 \otimes \sfh_{\ttS}).
\]
Suppose that $k+1$ is above $k$.
Then $\opi_k \cdot \ttS = -\ttS$ and $k \in \Des_L(\wrow(\ttS))$.
The latter implies that
\[
T_k T_{\wrow(\ttS)} = T_k^2 T_{\wrow(s_k \cdot \ttS)}    = (v-1) T_{\wrow(\ttS)} + v T_{\wrow(s_k \cdot \ttS)}.
\]
Thus,
$\opi_k \cdot 1 \otimes \sfh_{\ttS} 
= 1 \otimes T_k  \sfh_{\ttS} 
= 1 \otimes ((v-1)  \sfh_{\ttS} + v \sfh_{s_k \cdot \ttS})
= - 1 \otimes \sfh_{\ttS}$.
Consequently,
\[
\psi(\opi_k \cdot 1 \otimes \sfh_{\ttS}) = \psi(- 1 \otimes \sfh_{\ttS}) 
= -\ttS = \opi_k \cdot \ttS
= \opi_k \cdot \psi(1 \otimes \sfh_{\ttS}).
\]

Therefore, $\psi$ is an $\compH{0}$-module isomorphism, as required.
\end{proof}

For readers' understanding, 
we close this subsection by providing two tables on structural properties of the $\genericH$-module $S^{\lambda/\mu}(\C(v^{1/2});v)$ and the $\compH{0}$-module $X_{\lambda/\mu}$. 
We say that $\lambda/\mu$ is  \emph{connected} if its Young diagram is connected, and is a \emph{ribbon} if its Young diagram does not contain any $2 \times 2$ square.
The first table concerns their indecomposability (see \cref{table: indecomposability of X and S}).
The indecomposability of $S^{\lambda/\mu}(\C(v^{1/2});v)$ is obtained by considering the Schur expansion of its image under Frobenius characteristic.  
For the indecomposability of $X_{\lambda/\mu}$, see \cite[Proposition 7.3]{23KLO}.

The second table is about the projectivity of $S^{\lambda/\mu}(\C(v^{1/2});v)$ and $X_{\lambda/\mu}$ (see  \cref{table: projectivity comparing}).
The content on  $X_{\lambda/\mu}$
follows from the subsequent proposition.
\begin{proposition}\label{prop: proj and ribbon} 
Let $\lambda/\mu$ be a  skew partition  of  $n$. Then, 
$X_{\lambda/\mu}$
is projective if and only if $\lambda/\mu$ is a ribbon.
\end{proposition}
\begin{proof}
The ``if'' part follows from   \cite[Remark A.1(2)]{23KLO}.
For the ``only if'' part, suppose that $X_{\lambda/\mu}$ is projective. 
By  \cite[Lemma 5.4]{23KLO}  together with \cite[Remark A.1(2)]{23KLO}, the projective cover of $X_{\lambda/\mu}$ is  $X_{\nu/\kappa}$
for some ribbon $\nu/\kappa$. 
Since $X_{\lambda/\mu}$ is projective, it is isomorphic to its projective cover $X_{\nu/\kappa}$.
Therefore, by \cite[Theorem 5.5]{23KLO}, $\lambda/\mu=\nu/\kappa$,
and thus $\lambda/\mu$ is a ribbon.
\end{proof}

\begin{table}[t]
\centering
\begin{tabular}{c|c!{\vrule width 1.5pt}c|c}
\hline
 \multicolumn{2}{c!{\vrule width 1.5pt}}{$\lambda/\mu$} & $S^{\lambda/\mu}(\C(v^{1/2});v)$ 
& $X_{\lambda/\mu}$ 
 \\
\Xhline{1.5pt}
\multirow{2}{*}{connected} & \makecell{of the form $\nu$ or $\nu^\circ$}
& irreducible & indecomposable \\
\cline{2-4}
& \makecell{otherwise}
&  not irreducible
& indecomposable\\
\hline
\multirow{2}{*}{disconnected} & 
contains a disconnected ribbon&not irreducible  & not indecomposable\\
\cline{2-4}
& 
otherwise & not irreducible & unknown\\
\hline
\end{tabular}
\captionsetup{width=\textwidth}
\caption{Indecomposability of $S^{\lambda/\mu}(\C(v^{1/2});v)$ and 
$X_{\lambda/\mu}$}
\label{table: indecomposability of X and S}
\end{table}
\begin{table}[t]
\centering
\begin{tabular}{c!{\vrule width 1.5pt}c|c}
\hline
$\lambda/\mu$ & $S^{\lambda/\mu}(\C(v^{1/2});v)$
& $X_{\lambda/\mu}$ \\\Xhline{1.5pt}
 ribbon
& projective &  projective\\
\hline
not a ribbon &  projective & not projective \\
\hline
\end{tabular}
\captionsetup{width=\textwidth}
\caption{Projectivity of $S^{\lambda/\mu}(\C(v^{1/2});v)$ and 
$X_{\lambda/\mu}$}
\label{table: projectivity comparing}
\end{table}

\subsection{Kazhdan--Lusztig cell modules and dual plactic class modules}\label{subsec: 0-Hecke module from specialization 2}

In this subsection, we consider $\polyH$-modules arising from the
Kazhdan--Lusztig basis and from a rescaling of this basis, together with
their subquotients. We describe their images under the specialization
functors $\calF_\iota$ and $\calF_{\theta_q}$, where $q \in \C$.
In this way, we provide a relationship between
Kazhdan--Lusztig cell modules and dual plactic class modules.

We begin by reviewing the Kazhdan--Lusztig basis and Kazhdan--Lusztig cell modules.
For $f(v^{1/2}) \in \C[v^{\pm 1/2}]$, let $\overline{f(v^{1/2})} = f(v^{-1/2})$.
Using this, Kazhdan and Lusztig~\cite{79KL} introduced the bar involution  $\bar{~}\bar{~}$ of $\laurH$ by defining
\[
\overline{\sum_{w\in\SG_n} a_w T_w} = \sum_{w\in\SG_n} \overline{a_w} \, T_{w^{-1}}^{-1}  \quad \text{for all $a_w \in \C[v^{\pm1/2}]$ and $w \in \SG_n$}.
\]
Recall that the \emph{$($strong$)$ Bruhat order} $\leq$ is the partial order on $\SG_n$ generated by the following covering relations:
\[
\sigma \leq^{\rmc} \rho \quad\text{if $\ell(\rho) = \ell(\sigma)+1$ and $\sigma^{-1} \rho$ is a transposition.}
\]

\begin{theorem}{\rm (\cite[Theorem 1.1]{79KL})}\label{thm: KL basis thm}
For $w \in \SG_n$, 
there exists a unique $C_w \in H_{\C[v^{\pm 1/2}];v}(\SG_{n})$ such that 
\begin{enumerate}[label = {\rm (\arabic*)}]
\item 
$\overline{C_w} = C_w$, and
\item
$C_w = \sum_{y \le w} \sgn(y)\sgn(w) v^{\ell(w)/2} v^{-\ell(y)} \overline{P_{y,w}} T_y$,
where $P_{y,w} \in \Z[v]$ is of degree at most $(\ell(w) - \ell(y) - 1)/2$ for $y < w$, and $P_{w,w} = 1$.
\end{enumerate}
\end{theorem}
By \cref{thm: KL basis thm}(2), we have 
\begin{align}\label{eq: Kazhdan Lusztig basis elt into standard basis}
C_w = v^{-\ell(w)/2} T_w + \sum_{y < w} d_{y, w} T_y \quad (d_{y,w} \in \Z[v^{\pm 1/2}]).
\end{align}
Hence, $\{C_w \mid w \in \SG_n\}$  forms a basis for $\laurH$, called the \emph{Kazhdan--Lusztig basis}.
The polynomials $P_{y,w}~(y \le w)$ are called the \emph{Kazhdan--Lusztig polynomials}.
For $y,w \in \SG_n$ with $y \not\le w$, we define $P_{y,w} = 0$ for convenience.

To describe the action of $T_i$'s on each $C_w$, we need another order on $\SG_n$ called the \emph{Kazhdan--Lusztig order}.
It is defined as follows:
for $y, w \in \SG_n$, we write $y \prec w$ if $y < w$ and $P_{y,w}$ has degree $(\ell(w) - \ell(y) - 1)/2$.
In addition, for $y, w \in \SG_n$, define $\mu(y,w) \in \Z$ to be the coefficient of $v^{(\ell(w) - \ell(y) - 1)/2}$ in $P_{y,w}$.
Note that $\mu(y,w) \neq 0$ if and only if $y \prec w$.
For $i \in [n-1]$ and $w \in \SG_n$, the action of $T_i$ on $C_w$ is described in~\cite[p.171]{79KL} as follows:
\begin{align}\label{eq: Ti action on KL elements}
T_i \cdot C_w = \begin{cases}
-C_w & \text{if $i \in \Des_L(w)$}, \\[2ex]
v^{1/2} C_{s_iw} + v C_w + v^{1/2} \displaystyle
\sum_{\substack{z \prec w \\ s_i z < z}} \mu(z,w) C_z & \text{if $i \notin \Des_L(w)$}.
\end{cases}
\end{align}

For $\sigma, \rho \in \SG_n$, we write $\sigma \Larrow  \rho$ if $C_\sigma$ appears with nonzero coefficient in $h C_\rho$ for some $h \in \laurH$.
We note that $\Larrow$ is reflexive.
The \emph{left cell preorder} $\Lorder'$ is the transitive closure of $\Larrow$.
Let $\Lsim$ be the relation on $\SG_n$ given by $\sigma \Lsim \rho$ if $\sigma \Lorder' \rho$ and $\rho \Lorder' \sigma$.
By definition, $\Lsim$ is an equivalence relation.
The equivalence classes of $\Lsim$ are called the \emph{left cells} of $\SG_n$.
It is well known that for $w \in \SG_n$, the left cell containing $w$ is given by
\begin{align}\label{eq: left cell and Q tableau}
\{u \in \SG_n \mid \Qtab(u) = \Qtab(w)\}
\end{align}
(see, for example, \cite{00Ariki, 88GM, 79KL}).
Hence, left cells are parametrized by $\SYT_n$, and the left cell corresponding to $Q \in \SYT_n$ is exactly equal to $\Gamma_Q$.
The preorder $\Lorder'$ induces the partial order $\Lorder$ on $\{\Gamma_Q \mid Q \in \SYT_n \}$ defined by 
$\Gamma_Q \Lorder \Gamma_{Q'}$ if there exist $u \in \Gamma_Q$ and $w \in \Gamma_{Q'}$ such that $u \Lorder' w$.
We write $\Gamma_Q <_{\scrL} \Gamma_{Q'}$ if $\Gamma_Q \Lorder \Gamma_{Q'}$ and $\Gamma_Q \neq \Gamma_{Q'}$.

To give a uniform definition of Kazhdan--Lusztig cell modules for various Hecke algebras, we let $R$ be a PID with an invertible element $a$.
Suppose that there exists a ring homomorphism $f:\C[v^{\pm 1/2}]\ra R$ satisfying $f(v) = a$.
Recall in \cref{subsec: def of specialization functor} that this homomorphism  induces the ring homomorphism $\widetilde{f}: \laurH \ra H_{R;a}(\SG_n)$.
Since $f(v^{1/2})$ is invertible, it follows from \cref{eq: Kazhdan Lusztig basis elt into standard basis} that the set $\{\widetilde{f}(C_w) \mid w \in \SG_n \}$ forms an $R$-basis for $H_{R;a}(\SG_n)$.

For $Q \in \SYT_n$, let 
\begin{align*}
\frakJ_Q(R;a) 
& := R\text{-span of }  \{\widetilde{f}(C_w) \mid  \text{$w \in \Gamma_{Q'}$ for some $Q' \in \SYT_n$ with $\Gamma_{Q'} \Lorder \Gamma_Q$} \}\ \  \text{and} \\
\hfrakJ_Q(R;a)
& := R\text{-span of } \{\widetilde{f}(C_w) \mid \text{$w \in \Gamma_{Q'}$ for some $Q' \in \SYT_n$ with $\Gamma_{Q'} <_{\scrL} \Gamma_Q$}  \}.
\end{align*}

We remark that there may be more than one choice of a ring homomorphism $f:\C[v^{\pm 1/2}]\ra R$ satisfying $f(v)=a$. Nevertheless, the definitions of $\frakJ_Q(R;a)$ and $\hfrakJ_Q(R;a)$ do not depend on these choices.
To see this, let $g$ be another such ring homomorphism.
Since $g(v) = a = f(v)$, by \cref{thm: KL basis thm}(2),  $g(C_w) = \pm f(C_w)$ for all $w \in \SG_n$.

By the definition of the preorder $\Lorder'$, it is clear that $\frakJ_Q(R;a)$ and $\hfrakJ_Q(R;a)$ are left ideals of $H_{R;a}(\SG_n)$.
We define $\frakM_Q(R;a)$ by the left $H_{R;a^2}(\SG_{n})$-module 
\[
\frakM_Q(R;a) := \frakJ_Q(R;a)  \,/\, \hfrakJ_Q(R;a).
\]
This module is called the \emph{Kazhdan--Lusztig $($left$)$ cell module of $H_{R;a}(\SG_{n})$} associated with the left cell $\Gamma_Q$.
By definition, $\left\{ \widetilde{f}(C_w) + \hfrakJ_Q(R;a)  \mid  w \in \Gamma_Q \right\}$ forms a basis for $\frakM_Q(R;a)$.
For further background and details on Kazhdan--Lusztig $($left$)$ cell modules, see \cite{79KL, 96GL, 99Mat, 03Geordie}.

We next introduce the homomorphisms needed for the arguments in this subsection.
Let $\iota_0:\C[v^{1/2}]\hookrightarrow \C[v^{\pm 1/2}]$ and $\iota':\C[v^{\pm 1/2}]\hookrightarrow \C(v^{1/2})$
be the natural embeddings. 
For a nonzero $q\in\C$, let $\theta_q':\C[v^{\pm 1/2}]\to\C$
be the evaluation homomorphism at $v^{\pm 1/2} = q^{\pm 1}$.

Now, we introduce two kinds of $\polyH$-modules arising from the Kazhdan--Lusztig basis and describe their images under the specialization functions.
We first consider the $\C[v^{1/2}]$-span of $\{C_w \mid w\in\SG_n\}$, which is closed under the $\polyH$-action by \cref{eq: Ti action on KL elements}.
We denote the resulting $\polyH$-module by $\bfM^{(C)}_n$.
The following proposition describes its images under the specialization functors.
\begin{proposition}\label{lem: Cn isom to regular and simples}
\begin{enumerate}[label = {\rm (\arabic*)}]
    \item $\calF_{\iota} (\bfM^{(C)}_n) \cong H_{\C(v^{1/2});v}(\SG_n)$  as $\genericH$-modules.
    \item For a nonzero $q \in \C$, $\calF_{\theta_q} (\bfM^{(C)}_n) \cong H_{\C; q^2}(\SG_n)$  as $\compH{q^2}$-modules.
    \item $\calF_{\theta_0}(\bfM^{(C)}_n) \cong \bigoplus_{w \in \SG_n} \bfF_{\comp(\Des_L(w))}$ as $\compH{0}$-modules.
\end{enumerate}
\end{proposition}
\begin{proof}
(1) Since $\iota = \iota' \circ \iota_0$, we have
\[
\calF_{\iota}(\bfM^{(C)}_n) \cong \calF_{\iota'}(\calF_{\iota_0}(\bfM^{(C)}_n))\cong \calF_{\iota'}(\laurH) \cong \genericH,
\]
as required.

(2) Since $\theta_q = \theta'_q \circ \iota_0$, we have 
\[
\calF_{\theta_q}(\bfM^{(C)}_n) \cong \calF_{\theta'_q}(\calF_{\iota_0}(\bfM^{(C)}_n))\cong \calF_{\theta'_q}(\laurH) \cong \compH{q^2},
\]
as required.

(3) Note that $\bfM^{(C)}_n$ is free over $\C[v^{1/2}]$.
By \cref{rem: poly action and hecke action under free cond}, $\{1 \otimes C_w \mid w \in \SG_n \}$ forms a basis for $\calF_{\theta_0}(\bfM^{(C)}_n)$.
Moreover, it follows from \cref{eq: specialization and module action} and \cref{eq: Ti action on KL elements}  that the $\compH{0}$-action on this basis  is given as follows: for $w \in \SG_n$  and for $1 \leq i \leq n-1$,
\begin{align}\label{eq: F_o Mcn action}
\opi_i \cdot 1 \otimes C_w = \begin{cases}
- 1 \otimes C_w & \text{if $i \in \Des_L(w)$},\\
0 & \text{if $i \not\in \Des_L(w)$}.
\end{cases}
\end{align}
This shows that for each $w \in \SG_n$, $\C\{1 \otimes C_w\}$ is isomorphic to $\bfF_{\comp(\Des_L(w))^\rmc}$, and thus we obtain the desired isomorphism.
\end{proof}

We then introduce  subquotients of $\bfM^{(C)}_n$, which can be viewed as an analogue of the Kazhdan--Lusztig cell modules.
For $Q \in \SYT_n$, let
\begin{align*}
\bfJ^{(C)}_Q
&:= \C[v^{1/2}]\text{-span of } \, \{C_w \mid  \text{$w \in \Gamma_{Q'}$ for some $Q' \in \SYT_n$ with $\Gamma_{Q'} \Lorder \Gamma_Q$} \}\ \  \text{and}  \\
\hbfJ^{(C)}_Q
&:= \C[v^{1/2}]\text{-span of } \,  \{C_w \mid  \text{$w \in \Gamma_{Q'}$ for some $Q' \in \SYT_n$ with $\Gamma_{Q'} <_{\scrL} \Gamma_Q$} \}.
\end{align*}
From the definition of $\Lorder$ together with \cref{eq: analogue of KL elts}, one can verify that the $\bfJ^{(C)}_Q$ and $\hbfJ^{(C)}_Q$ are stable under the action of $\polyH$.
We define the $\polyH$-module $\bfM^{(C)}_Q$ by
\begin{align}\label{eq: def of bfMCQ}
\bfM^{(C)}_Q := \bfJ^{(C)}_Q / \hbfJ^{(C)}_Q.
\end{align}
The following theorem describes the images of $\bfM^{(C)}_Q$ under the specialization functors.

\begin{theorem}\label{thm: KL cell module C basis}
Let  $Q \in \SYT_n$.
\begin{enumerate}[label = {\rm (\arabic*)}]
\item 
$\calF_{\iota}(\bfM^{(C)}_Q) \cong \frakM_Q(\C(v^{1/2});v)$
as $\genericH$-modules.
\item 
For nonzero $q \in \C$, $\calF_{\theta_q}(\bfM^{(C)}_Q) \cong \frakM_Q(\C;q)$ as $\compH{q^2}$-modules.
\item 
$\calF_{\theta_0}(\bfM^{(C)}_Q) \cong \bigoplus_{w \in \Gamma_Q}  \bfF_{\comp(\Des_L(w))^\rmc}$ as $\compH{0}$-modules.
\end{enumerate} 
\end{theorem}
\begin{proof}
(1)  We have an $\laurH$-module isomorphism
\begin{align*}
\calF_{\iota_0}(\bfM^{(C)}_Q) \cong \frakM_Q(\C[v^{\pm 1/2}];v).
\end{align*}
Since $\iota = \iota' \circ \iota_0$, we have 
\[
\calF_{\iota}(\bfM^{(C)}_Q) \cong \calF_{\iota'}(\frakM_Q(\C[v^{\pm 1/2}];v)) \cong \frakM_Q(\C(v^{1/2});v).
\]

(2) The assertion can be proved in the same way as in (1).

(3) Note that  $\{C_w + \hbfJ^{(C)}_Q \mid w \in \Gamma_Q \}$ is a $\C[v^{1/2}]$-basis for $\bfM^{(C)}_Q$.
And, in view of \cref{rem: poly action and hecke action under free cond}, $\{1 \otimes (C_w + \hbfJ^{(C)}_Q) \mid w \in \Gamma_Q \}$ is a $\C$-basis for $\calF_{\theta_0}(\bfM^{(C)}_Q)$.
For each $w \in \Gamma_Q$, it follows from \cref{eq: F_o Mcn action} that  $\C \{1 \otimes (C_w + \hbfJ^{(C)}_Q) \}$ is an $\compH{0}$-submodule of $\calF_{\theta_0}(\bfM^{(C)}_Q)$ isomorphic to $\bfF_{\comp(\Des_L(w))^\rmc}$.
Thus, we obtain the desired isomorphism.
\end{proof}

Next, we define two kinds of $\polyH$-modules arising from the rescaled
Kazhdan--Lusztig basis $\{D_w \mid w\in \SG_n\}$, where
\[
D_w := v^{\ell(w)/2} C_w.
\]
We first consider the $\C[v^{1/2}]$-span of $\{D_w \mid w \in \SG_n \}$.
From \cref{eq: Ti action on KL elements} it follows that for $i \in [n-1]$ and $w \in \SG_n$,
\begin{align}\label{eq: analogue of KL elts}
T_i \cdot D_w = \begin{cases}
- D_w & \text{if $i \in \Des_L(w)$,} \\[2ex]
\displaystyle D_{s_i w} + \sum_{z \le w} f_z(v^{1/2}) D_z & \text{if $i \notin \Des_L(w)$,}
\end{cases}
\end{align}
where $f_z(v^{1/2}) \in v^{1/2}\Z[v^{1/2}]$.
Therefore, this  space is closed under the $\polyH$-action. 
We denote the resulting $\polyH$-module by $\bfM^{(D)}_n$.

\begin{lemma}\label{lem: Dw in polyH}
\begin{enumerate}[label = {\rm (\arabic*)}]
\item 
For $w \in \SG_n$, $D_w$ is an element of $\polyH$. 
More precisely,
\begin{align*}
D_w = T_w + \sum_{y < w} g_{y, w} T_y \quad(g_{y, w} \in v^{1/2} \Z[v^{1/2}]).
\end{align*}
\item 
$\bfM^{(D)}_n$ is the regular representation of $\polyH$. 
\end{enumerate}
\end{lemma}
\begin{proof}
(1) Let $w \in \SG_n$.
By \cref{thm: KL basis thm}, we have
\begin{align*}
D_w =  \sum_{y \le w} \sgn(y)\sgn(w) v^{\ell(w) - \ell(y)} \overline{P_{y,w}} T_y = T_w + \sum_{y < w} \sgn(y)\sgn(w) v^{\ell(w) - \ell(y)} \overline{P_{y,w}} T_y. 
\end{align*}
Let  $y \in \SG_n$ with $y < w$.
It suffices to show that
\begin{align}\label{eq: y < w Dw Ty coeff} 
v^{\ell(w)-\ell(y)}\overline{P_{y,w}} \in v^{1/2}\Z[v^{1/2}].
\end{align} 
Again, by \cref{thm: KL basis thm},
\[ 
\deg P_{y,w} \leq \frac{\ell(w)-\ell(y)-1}{2}. 
\]
It follows that
\begin{align*}
\ell(w) - \ell(y) - \deg P_{y, w} 
&\geq \ell(w) - \ell(y) - \frac{\ell(w) - \ell(y) - 1}{2} 
= \frac{\ell(w) - \ell(y) + 1}{2} 
> 0.
\end{align*}
Combining this with the fact that the degree of every monomial appearing in $\overline{P_{y,w}}$ is greater than or equal to $-\deg P_{y,w}$, we have \cref{eq: y < w Dw Ty coeff}.

(2) By (1), $\{ D_w \mid w \in \SG_n \}$ is a basis for $\polyH$, and the transition matrix from $\{ D_w \mid w \in \SG_n \}$ to $\{ T_w \mid w \in \SG_n \}$ is unitriangular.
Therefore, $\bfM^{(D)}_n = \polyH$.
\end{proof}

The following proposition, which describes the images of $\bfM^{(D)}_n$ under the specialization functors,
 is an immediate consequence of 
\cref{lem: Dw in polyH}(2). 

\begin{proposition}\label{lem: Dn isom to regular}
\begin{enumerate}[label = {\rm (\arabic*)}]
    \item $\calF_{\iota} (\bfM^{(D)}_n) \cong H_{\C(v^{1/2});v}(\SG_n)$  as $\genericH$-modules.
    \item For a nonzero $q \in \C$, $\calF_{\theta_q} (\bfM^{(D)}_n) \cong \compH{q^2}$  as $\compH{q^2}$-modules.
    \item $\calF_{\theta_0}(\bfM^{(D)}_n) \cong \compH{0}$ as $\compH{0}$-modules.
\end{enumerate}
\end{proposition}

We next introduce the most notable family 
of $\polyH$-modules in this subsection.
It consists of the
subquotients of $\bfM^{(D)}_n$, which can be viewed as another analogue of the Kazhdan--Lusztig cell modules.
For $Q \in \SYT_n$, let
\begin{align*}
\bfJ^{(D)}_Q
&:= \C[v^{1/2}]\text{-span of } \, \{D_w \mid  \text{$w \in \Gamma_{Q'}$ for some $Q' \in \SYT_n$ with $\Gamma_{Q'} \Lorder \Gamma_Q$} \}\ \  \text{and}  \\
\hbfJ^{(D)}_Q
&:= \C[v^{1/2}]\text{-span of } \, \{D_w \mid  \text{$w \in \Gamma_{Q'}$ for some $Q' \in \SYT_n$ with $\Gamma_{Q'} <_{\scrL} \Gamma_Q$}  \}.
\end{align*}
From the definition of $\Lorder$ together with \cref{eq: analogue of KL elts}, we derive that  $\bfJ^{(D)}_Q$ and $\hbfJ^{(D)}_Q$ are stable under the action of $\polyH$.
We define the $\polyH$-module $\bfM^{(D)}_Q$ by
\begin{align}\label{def: def of MDQ}
\bfM^{(D)}_Q := \bfJ^{(D)}_Q / \hbfJ^{(D)}_Q.
\end{align}
The following theorem describes the images of $\bfM^{(D)}_Q$ under the specialization functors.

\begin{theorem}\label{thm: obtain CGammaQ}
Let  $Q \in \SYT_n$.
\begin{enumerate}[label = {\rm (\arabic*)}]
\item $\calF_{\iota}(\bfM^{(D)}_Q) \cong \frakM_Q(\C(v^{1/2}); v)$ as $\genericH$-modules. 
\item For nonzero $q \in \C$, $\calF_{\theta_q}(\bfM^{(D)}_Q) \cong \frakM_Q(\C;q)$ as $\compH{q^2}$-modules.
\item 
$\calF_{\theta_0}(\bfM^{(D)}_Q) \cong \scrD_Q$ as $\compH{0}$-modules.
\end{enumerate} 
\end{theorem}
\begin{proof}
Note that by \cref{lem: Dw in polyH}(2), $\bfJ^{(D)}_Q$ and  $\hbfJ^{(D)}_Q$ are submodules of the regular representation of $\polyH$.
In addition, $\bfM^{(D)}_Q = \bfJ^{(D)}_Q / \hbfJ^{(D)}_Q$ is free over $\C[v^{1/2}]$.

(1) Since 
\[
\inn{\widetilde{\iota}(\bfJ^{(D)}_Q)} = \frakJ_Q(\C(v^{1/2});v) \quad\text{and}\quad
\inn{\widetilde{\iota}(\hbfJ^{(D)}_Q)} = \hfrakJ_Q(\C(v^{1/2});v),
\]
it follows from  \cref{cor: algebra homomorphism vs functor sub of regular}(1) that  
\[
\calF_\iota(\bfM^{(D)}_Q) \cong \inn{\widetilde{\iota}(\bfJ^{(D)}_Q)}/\inn{\widetilde{\iota}(\bfJ^{(D)}_Q)} = \frakJ_Q(\C(v^{1/2});v)/\hfrakJ_Q(\C(v^{1/2});v) = \frakM_Q(\C(v^{1/2});v).
\]

(2) Note that 
\[
\widetilde{\theta_q}(\bfJ^{(D)}_Q) = \frakJ_Q(\C;q^2) \quad\text{and}\quad
\widetilde{\theta_q}(\hbfJ^{(D)}_Q) = \hfrakJ_Q(\C;q^2).
\]
Since $\dim_{\C} \frakM_Q(\C;q^2) = |\Gamma_Q| = r_{\bfM^{(D)}_Q}$, 
it follows from  \cref{cor: algebra homomorphism vs functor sub of regular}(2) that  
\[
\calF_{\theta_q}(\bfM^{(D)}_Q) \cong \widetilde{\theta_q}(\bfJ^{(D)}_Q)/\widetilde{\theta_q}(\hbfJ^{(D)}_Q) = \frakJ_Q(\C;q^2)/\hfrakJ_Q(\C;q^2) = \frakM_Q(\C;q^2).
\]

(3) By \cref{lem: Dw in polyH}(1), we have
\begin{align*}
\widetilde{\theta_0}(D_w) = \opi_w \quad \text{for $w \in \SG_n$}, 
\end{align*}
and thus, 
\begin{align*}
\widetilde{\theta_0}(\bfJ^{(D)}_Q)
& = 
\{\opi_w \in \SG_n \mid  \text{$w \in \Gamma_{Q'}$ for some $Q' \in \SYT_n$ with $\Gamma_{Q'} \Lorder \Gamma_Q$} \} \ \ \text{and}
\\
\widetilde{\theta_0}(\hbfJ^{(D)}_Q)
& =
\{\opi_w \in \SG_n \mid  \text{$w \in \Gamma_{Q'}$ for some $Q' \in \SYT_n$ with $\Gamma_{Q'} <_{\scrL} \Gamma_Q$} \}.
\end{align*}
Then, in view of \cref{eq: action on DQ}, we have an $\compH{0}$-module isomorphism 
\[
\widetilde{\theta_0}(\bfJ^{(D)}_Q)/
\widetilde{\theta_0}(\hbfJ^{(D)}_Q) \cong \scrD_Q.
\]
Since $\dim_{\C} \scrD_Q = |\Gamma_Q| = r_{\bfM^{(D)}_Q}$, it follows from  \cref{cor: algebra homomorphism vs functor sub of regular}(2) that   
\[
\calF_{\theta_0}(\bfM^{(D)}_Q) \cong \widetilde{\theta_0}(\bfJ^{(D)}_Q)/ \widetilde{\theta_0}(\hbfJ^{(D)}_Q) \cong \scrD_Q.
\qedhere
\]
\end{proof}

The $\polyH$-module  $\bfM^{(D)}_n$ admits a filtration that behaves well under the specialization functors  $\calF_{\iota}$ and $\calF_{\theta_0}$. 
To see this, write $\SYT_n$ as $=\{Q_1,Q_2,\ldots,Q_l\}$, where the  elements are enumerated in such a way that $i < j$ whenever $\Gamma_{Q_i} <_{\scrL} \Gamma_{Q_j}$.
Let 
\[
M^{(D)}_i:=\sum_{1 \leq t \leq i} \bfJ^{(D)}_{Q_t} \quad(1 \leq i \leq l) \quad\text{and}\quad M^{(D)}_0 := 0.
\]
Since each $\bfJ^{(D)}_{Q_t}$ is a left ideal of $\polyH$,  $M^{(D)}_i$ is an $\polyH$-module. 
Thus, we have a filtration 
\begin{align}\label{eq: filt of M D}
0=M^{(D)}_0 \subset M^{(D)}_1 \subset M^{(D)}_2 \subset \cdots \subset M^{(D)}_l=\bfM^{(D)}_n.
\end{align}

\begin{proposition}\label{prop: filt of MnD}
In the filtration \cref{eq: filt of M D}, we have the following.
\begin{enumerate}[label = {\rm (\arabic*)}]
\item  $\calF_{\iota}(M^{(D)}_i) / \calF_{\iota}(M^{(D)}_{i-1}) \cong \frakM_{Q_i}(\C(v^{1/2});v)$ for all $1 \leq i \leq l$.
\item $\calF_{\theta_0}(M^{(D)}_i) / 
\calF_{\theta_0}(M^{(D)}_{i-1}) \cong \scrD_{Q_i}$ for all $1 \leq i \leq l$.
\end{enumerate}
\end{proposition}
\begin{proof}
(1) Note that  the $\C(v^{1/2})$-linear map
\[
\calF_\iota(M^{(D)}_i) \ra \sum_{1 \leq t \leq i} \frakJ_{Q_t}(\C(v^{1/2});v), \quad 1 \otimes D_w \mapsto D_w
\quad
\text{for $w \in \bigcup_{1 \leq t \leq i} \Gamma_{Q_t}$}
\]
is an $\genericH$-module isomorphism.
Therefore,  
\begin{align*}
\calF_\iota(M^{(D)}_i) / \calF_{\iota}(M^{(D)}_{i-1}) 
& \cong \sum_{1 \leq t \leq i} \frakJ_{Q_t}(\C(v^{1/2});v) \, \big/  \sum_{1 \leq t < i} \frakJ_{Q_t}(\C(v^{1/2});v) \\
& \cong \frakJ_{Q_i}(\C(v^{1/2});v) \big/ \hfrakJ_{Q_i}(\C(v^{1/2});v) 
\\& = \frakM_{Q_i}(\C(v^{1/2});v).
\end{align*}
Here, the second isomorphism follows from the identity
\[
\frakJ_{Q_i}(\C(v^{1/2});v) \cap \sum_{1 \leq t < i} \frakJ_{Q_t} = \hfrakJ_{Q_i}(\C(v^{1/2});v).
\]

(2) By  \cref{lem: Dw in polyH}, $M^{(D)}_i$ is a submodule of the regular representation of $\polyH$ and  
\begin{align}\label{eq: theta 0  til M D i}
 \widetilde{\theta_0}(M^{(D)}_i) =  \C\{\opi_w \mid w \in \Gamma_{Q_t}  \text{ for some $1 \leq t \leq i$} \}.
\end{align}
By \cref{cor: algebra homomorphism vs functor sub of regular}(2),  $\calF_{\theta_0}(M^{(D)}_i)$ is isomorphic to $\widetilde{\theta_0}(M^{(D)}_i)$.
The explicit isomorphism is given as follows:
\[
\calF_{\theta_0}(M^{(D)}_i) \ra \widetilde{\theta_0}(M^{(D)}_i), \quad 1 \otimes D_w \mapsto \opi_w
\quad
\text{for $w \in \bigcup_{1 \leq t \leq i} \Gamma_{Q_t}$}
\]
(see the proof of \cref{cor: algebra homomorphism vs functor sub of regular}(2)).
Therefore, 
\[
\calF_{\theta_0}(M^{(D)}_i) / 
\calF_{\theta_0}(M^{(D)}_{i-1}) \cong \widetilde{\theta_0}(M^{(D)}_i)/\widetilde{\theta_0}(M^{(D)}_{i-1}) \cong \scrD_{Q_i}.
\]
Here, the third isomorphism follows from  \cref{eq: action on DQ} together with \cref{eq: theta 0  til M D i}.
\end{proof}

We summarize the results of this subsection in \cref{tab: M CD specialization}.
\begin{table}[h]
\centering
\begin{tabular}{c!{\vrule width 1.5pt}c|c|c|c}
 & $\calF_\iota(\bfM^{(X)}_n)$
& $\calF_{\theta_0}(\bfM^{(X)}_n)$
& $\calF_\iota(\bfM^{(X)}_Q)$
& $\calF_{\theta_0}(\bfM^{(X)}_Q)$ \\[1ex]
\Xhline{1.5pt}
&&&&\\[-2ex]
$X = C$
& $\genericH$
& $\displaystyle \bigoplus_{w \in \SG_n}
   \bfF_{\comp(\Des_L(w))^\rmc}$
& $\frakM_Q$
& $\displaystyle \bigoplus_{\ttT \in \SYT(\lambda)}
   \bfF_{\comp(\bfDes(\ttT))^\rmc}$ \\[3ex]
\hline
&&&&\\[-2ex]
$X = D$
& $\genericH$
& $\compH{0}$
& $\frakM_Q$
& $\scrD_Q$
\end{tabular}
\captionsetup{width=\textwidth}
\caption{The images of $\bfM^{(X)}_n$ and $\bfM^{(X)}_Q$, up to isomorphism, under $\calF_\iota$ and $\calF_{\theta_0}$ for $X = C, D$}
\label{tab: M CD specialization}
\end{table}

As seen in this table,  $\calF_{\theta_0}(\bfM^{(C)}_n)$ and $\calF_{\theta_0}(\bfM^{(C)}_Q)$ $(Q \in \SYT_n)$ are semisimple, and $\calF_{\theta_0}(\bfM^{(D)}_n)$ is a direct sum of projective indecomposable modules.  
The decomposability of $\calF_{\theta_0}(\bfM^{(D)}_Q)$ is not   known in general. 
However, if  $\scrD_Q$ is isomorphic to $X_{\lambda}$ or $X_{\lambda^\circ}$,  then $\calF_{\theta_0}(\bfM^{(D)}_Q)$ is indecomposable.
Indeed, for every $\lambda \vdash n$, there exists such a tableau $Q \in \SYT_n$ (for details, see the discussion immediately following \cref{eq: MP set in DQ set}).
As a byproduct,  we deduce the following interesting properties: 
\begin{enumerate}[label = $\bullet$, itemsep = 0.5ex]
\item
$\bfM^{(C)}_n \not\cong \bfM^{(D)}_n$ for all $n \geq 3$.
\item
$\bfM^{(C)}_Q \not\cong \bfM^{(D)}_Q$ if $\scrD_Q$ is not semisimple.
\end{enumerate}

\section{The Frobenius characteristic and the quasisymmetric characteristic}\label{sec: properties of 0-Hecke modules specializatio}

In this section,  we describe the  relationship between the Frobenius characteristic and the quasisymmetric characteristic via specialization functors. 
And, we characterize $0$-Hecke algebra modules whose images under the quasisymmetric characteristic are symmetric in terms of their characters.

To begin with, we introduce notation and basic properties concerning characters.
Let $R$ be a field and let $a \in R$.
Recall that for an $H_{R;a}(\SG_n)$-module $M$, the character $\chi_M$ of $M$ is defined by 
\[
\chi_M : H_{R;a}(\SG_n) \ra R, \quad h \mapsto \sftr_{M}(h).
\]

\begin{lemma}\label{lem: char and Grothendieck}
Let $R$ be a field with characteristic zero and let $a \in R$.
For $H_{R;a}(\SG_n)$-modules $M$ and $N$,
\begin{align*}
\chi_M = \chi_N \quad\text{if and only if} \quad \text{$[M] = [N]$ in $\calG_0(H_{R;a}(\SG_n)\text{-}\module)$.}
\end{align*}
\end{lemma}

\begin{proof}
The ``if'' part is clear from the definition of characters.

Let us prove the ``only if'' part.
Since $R$ has characteristic zero and the classes of irreducible $H_{R;a}(\SG_n)$-modules form a $\Z$-basis of the Grothendieck group, it suffices to show that the  characters of pairwise nonisomorphic irreducible $H_{R;a}(\SG_n)$-modules are linearly independent over $R$.
Every irreducible $H_{R;a}(\SG_n)$-module is absolutely irreducible.
If $a\neq 0$, this follows from \cite[Theorem~3.4]{96GL} together with the fact that $H_{R;a}(\SG_n)$ is a cellular algebra; otherwise, it follows from the fact that every irreducible $H_{R;0}(\SG_n)$-module is one-dimensional.
By \cite[Theorem~2.3.13]{89NT}, the characters of pairwise nonisomorphic  irreducible $H_{R;a}(\SG_n)$-modules are linearly independent over $R$.
Therefore, the assertion follows.
\end{proof}

For $\alpha \models n$ with $\set(\alpha) = \{i_1 < i_2 < \cdots < i_k \}$, let
\[
\sfcyc(\alpha) := (1 \, 2 \, \cdots \, i_1) (i_1+1 \; i_1 + 2 \; \cdots \; i_2) \cdots (i_k + 1 \; i_k + 2 \; \cdots \; n).
\]
If we write $[n-1] \setminus \set(\alpha)$ as $\{j_1 <  j_2 < \cdots < j_{n-1-k}\}$, then
\begin{align}\label{eq: reduced exp of cyc}
s_{j_1} s_{j_2}  \cdots s_{j_{n-1-k}}
\end{align}
is a reduced expression for $\sfcyc(\alpha)$.
In the case where $(R,a) = (\C(v^{1/2}), v)$, it was shown in \cite[Theorem~5.1]{91Ram} that for $\sigma \in \SG_n$, there exist $a_{\sigma, \lambda} \in \Z[v]$ ($\lambda \vdash n$) such that 
\begin{align*}
\chi_M(T_{\sigma}) 
= \chi_M\left( \sum_{\lambda \vdash n} a_{\sigma, \lambda} T_{\sfcyc(\lambda)}\right)\quad\text{for all $\genericH$-modules $M$.}
\end{align*}
Combining this with \cref{lem: char and Grothendieck} shows that for $\genericH$-modules $M$ and $N$, the following conditions are equivalent:
\begin{enumerate}[label = {\rm (\roman*)}]
\item 
$\chi_M = \chi_N$.
\item 
$\chi_M(T_{\sfcyc(\lambda)}) = \chi_N(T_{\sfcyc(\lambda)})$ for all $\lambda \vdash n$.
\item 
$[M] = [N]$ in $\calG_0(\genericH\text{-}\module)$.
\end{enumerate}

The following proposition is a $0$-Hecke analogue of the above equivalences.

\begin{proposition}\label{prop: four equivalence relation on 0 Hecke}
For $\compH{0}$-modules $M$ and $N$,  the following are equivalent.
\begin{enumerate}[label = {\rm (\arabic*)}]
\item $\chi_M = \chi_N$.
\item $\chi_M(\opi_{\sfcyc(\alpha)}) = \chi_N(\opi_{\sfcyc(\alpha)})$ for all $\alpha \models n$.
\item $[M] = [N]$ in $\calG_0(\compH{0}\text{-}\module)$.
\end{enumerate}
\end{proposition}
\begin{proof}
By \cref{lem: char and Grothendieck}, (1) is equivalent to (3).
It is clear that (1) implies (2), thus it remains to show that (2) implies (1).
For $\sigma \in \SG_n$, we set  
\[
\upalpha_{\sigma} := \comp([n-1] \setminus \{j_1, j_2, \ldots, j_k\})
\]
for any reduced expression $s_{j_1} s_{j_2} \cdots s_{j_k}$ for  $\sigma$.
This composition is well-defined by 
Matsumoto's theorem.
For the required implication, it suffices to show that for every $\compH{0}$-module $M$ and $\sigma \in \SG_n$,
\begin{align}\label{eq: character on opi sig and opi cyc comp}
\chi_{M}(\opi_{\sigma}) =  \sgn(\sigma)  \sgn(\sfcyc(\upalpha_{\sigma}))
\chi_{M}(\opi_{\sfcyc(\upalpha_{\sigma})}).
\end{align}

Let $\alpha \models n$.
By \cref{eq: simple module 0-Hecke action}, for any $\sigma \in \SG_n$,
\begin{align}\label{eq: simple module character}
&\chi_{\bfF_{\alpha}}(\opi_{\sigma}) = \begin{cases}
\sgn(\sigma) & \text{if $\set(\alpha) \subseteq \set(\upalpha_{\sigma})$,}\\
0 & \text{otherwise.}
\end{cases}
\end{align}
Also, for any $\beta \models n$, using the reduced expression for $\sfcyc(\beta)$ in \cref{eq: reduced exp of cyc}, we see that $\upalpha_{\sfcyc(\beta)} = \beta$.
Combining this equality with \cref{eq: simple module character}, we derive that for any $\sigma \in \SG_n$,
\begin{align}\label{eq: chi F and chi F cyc}
\chi_{\bfF_{\alpha}}(\opi_{\sigma}) 
= \sgn(\sigma)  \sgn(\sfcyc(\upalpha_{\sigma}))
\chi_{\bfF_{\alpha}}(\opi_{\sfcyc(\upalpha_{\sigma})}).
\end{align}
Now, the equality \cref{eq: character on opi sig and opi cyc comp} follows from \cref{eq: chi F and chi F cyc} because $\chi_M = \sum_{\beta \models n} c_{\beta} \chi_{\bfF_{\beta}}$ for any $\compH{0}$-module $M$ with $[M] = \sum_{\beta \models n} c_\beta [\bfF_\beta]$ in $\calG_0(\compH{0}\text{-}\module)$.
\end{proof}

We state the first main result of this section.

\begin{theorem}\label{thm: characteristic pres prop spec}
Let $M$ be an $\polyH$-module that is free over $\C[v^{1/2}]$.
Then 
\[
\ch([\calF_{\theta_0}(M)]) = \ch^{\rm Frob}([\calF_{\iota}(M)]).
\]
\end{theorem}
\begin{proof}
Let $c_\lambda \in \Z_{\geq 0}~(\lambda \vdash n)$ such that
\begin{align}\label{eq: calF iota M cong Specht}
\calF_\iota(M) \cong \bigoplus_{\lambda \vdash n} S^{\lambda}(\C(v^{1/2});v)^{\oplus c_\lambda}
\end{align}
and
\[
N := \bigoplus_{\lambda \vdash n} S^{\lambda}(\C[v^{1/2}];v)^{\oplus c_\lambda}.
\]
By \cref{thm: specialization of skew Specht module}(1), we have 
\[
\calF_\iota(M) \cong \calF_\iota(N) \quad\text{as $\genericH$-modules,}
\]
and hence $\chi_{\calF_\iota(M)} = \chi_{\calF_\iota(N)}$.
By \cref{thm: character and specialization}(2), it follows that 
\[
\chi_{\calF_{\theta_0}(M)} =  \chi_{\calF_{\theta_0}(N)}.
\]
Combining this with \cref{prop: four equivalence relation on 0 Hecke}, we have 
\begin{equation}\label{eq: same characteristic M and M'}
    \ch([\calF_{\theta_0}(M)]) = \ch([\calF_{\theta_0}(N)]).
\end{equation}
Moreover, by \cref{thm: specialization of skew Specht module}(3), 
\[
\calF_{\theta_0}(N) \cong  \bigoplus_{\lambda \vdash n} X_\lambda^{\oplus c_\lambda} \quad\text{as $\compH{0}$-modules}.
\]
This, together with \cref{eq: ch of X} and \cref{eq: same characteristic M and M'}, implies that 
\[
\ch([\calF_{\theta_0}(M)]) 
= \ch([\calF_{\theta_0}(N)]) 
= \sum_{\lambda\vdash n} c_\lambda s_\lambda.
\]
From \cref{eq: calF iota M cong Specht} it follows that $\ch^{\rm Frob}([\calF_\iota(M)]) = \sum_{\lambda\vdash n} c_\lambda s_\lambda$, thus the assertion follows.
\end{proof}

Next, we state the second main result of this section.

\begin{theorem}\label{thm: char of symmetric 0-Hecke}
Let $M$ be an $\compH{0}$-module.
Then  $\ch([M])$ is symmetric if and only if
$\chi_M(\opi_{\sfcyc(\alpha)}) =  \chi_M(\opi_{\sfcyc(\beta)})$ for all $\alpha, \beta \models n$ with $\lambda(\alpha) = \lambda(\beta)$.
\end{theorem}

\begin{proof} 
Let $V$ be the $\Z$-span of $\{ \chi_{\bfF_\alpha} \mid \alpha \models n\}$.
Define a $\Z$-linear map
\[
\sfchar:
\calG_0(\compH{0}\text{-}\module)
\to
V,
\quad [\bfF_\alpha] \mapsto \chi_{\bfF_\alpha} 
\quad \text{for $\alpha \models n$}.
\]
By \cref{lem: char and Grothendieck}, $\sfchar$ is a $\Z$-linear isomorphism and $\sfchar([N]) = \chi_N$ for every $\compH{0}$-module $N$.
Let
\[
W := \{f \in V \mid f(\opi_{\sfcyc(\alpha)}) = f(\opi_{\sfcyc(\beta)}) \text{ for all $\alpha, \beta \models n$ with $\lambda(\alpha) = \lambda(\beta)$}\}.
\]
Since $\ch$ and $\sfchar$ are $\Z$-linear isomorphisms, the assertion will follow once we establish the equality
\begin{align}\label{eq: ch inv and char inv}
\ch^{-1}(\Sym_n) = \sfchar^{-1}(W).
\end{align}
We prove this in the following two steps.
\smallskip

{\it Step 1.} We show that $\ch^{-1}(\Sym_n) \subseteq \sfchar^{-1}(W)$.
Since the elements $\ch([X_\lambda]) = s_\lambda$ $(\lambda \vdash n)$ form a $\Z$-basis for $\Sym_n$, it suffices to show that $\sfchar([X_\lambda]) \in W$ for all $\lambda \vdash n$, that is,
\[
\chi_{X_\lambda}(\opi_{\sfcyc(\alpha)}) = \chi_{X_\lambda}(\opi_{\sfcyc(\beta)})
\quad \text{for all $\alpha, \beta \models n$ with $\lambda(\alpha) = \lambda(\beta)$}.
\]
Fix $\lambda \vdash n$ and $\alpha, \beta \models n$ such that $\lambda(\alpha) = \lambda(\beta)$.
By \cref{thm: specialization of skew Specht module},
\[
\calF_{\iota}(S^{\lambda}(\C[v^{1/2}]; v)) \cong S^{\lambda}(\C(v^{1/2}); v) \quad\text{and}\quad
\calF_{\theta_0}(S^{\lambda}(\C[v^{1/2}]; v)) \cong X_{\lambda}.
\]
Combining this with \cref{thm: character and specialization}(2) yields
\begin{align}\label{eq: W and Symn}
\chi_{X_{\lambda}}(\opi_\sigma) 
= \theta_0(\chi_{S^{\lambda}(\C(v^{1/2});v)}(T_\sigma)) \quad\text{for all $\sigma \in \SG_n$}.
\end{align}
Note that $\sfcyc(\alpha)$ and $\sfcyc(\beta)$ are minimal-length permutations in the same conjugacy class.
By \cite[(1) on p.~261]{00GP}, $T_{\sfcyc(\alpha)}-T_{\sfcyc(\beta)}$ is a finite sum of commutators in $\genericH$.
Since the trace of every commutator is zero, it follows that
\begin{align}\label{eq: same cycle type has the same trace}
\chi_{S^{\lambda}(\C(v^{1/2});v)}(T_{\sfcyc(\alpha)}) 
= \chi_{S^{\lambda}(\C(v^{1/2});v)}(T_{\sfcyc(\beta)}).
\end{align}
Using \cref{eq: W and Symn} and
\cref{eq: same cycle type has the same trace}, we have
\begin{align*}
\chi_{X_{\lambda}}(\opi_{\sfcyc(\alpha)}) 
= \theta_0(\chi_{S^{\lambda}(\C(v^{1/2});v)}(T_{\sfcyc(
\alpha)}))
= \theta_0(\chi_{S^{\lambda}(\C(v^{1/2});v)}(T_{\sfcyc(
\beta)}))
= \chi_{X_\lambda}(\opi_{\sfcyc(\beta)}),
\end{align*}
as required.
\smallskip

{\it Step 2.} We show that $\sfchar^{-1}(W) / \ch^{-1}(\Sym_n)$ is zero.
Note that $\Qsym_n/\Sym_n$ is torsion-free over $\Z$.
Since $\sfchar^{-1}(W)/\ch^{-1}(\Sym_n)$ is isomorphic to a submodule of $\Qsym_n/\Sym_n$ via $\ch$, it is also torsion-free over $\Z$.
Thus it suffices to show that $\sfchar^{-1}(W)$ and $\ch^{-1}(\Sym_n)$ have the same $\Z$-rank, or equivalently, that $W$ and $\Sym_n$ have the same $\Z$-rank.

Since $\ch$ and $\sfchar$ preserve $\Z$-rank, {\it Step~1} shows that the $\Z$-rank of $W$ is greater than or equal to the $\Z$-rank of $\Sym_n$, which is the number of partitions of $n$.
By \cref{prop: four equivalence relation on 0 Hecke}, every $f \in W$ is uniquely determined by the integer values $f(\opi_{\sfcyc(\alpha)})$ for $\alpha \models n$.
By the definition of $W$, these values depend only on $\lambda(\alpha)$, so $f$ is uniquely determined by
$f(\opi_{\sfcyc(\lambda)})$ for $\lambda \vdash n$.
It follows that the $\Z$-rank of $W$ is less than or equal to the number of partitions of $n$.
Thus $W$ and $\Sym_n$ have the same $\Z$-rank.
\end{proof}

\section{Removable-corner filtrations of restrictions of skew Specht modules}\label{subsec:Relationship between restriction and specialization functors for skew Specht modules}

\subsection{The removable-corner filtration of $S^{\lambda/\mu}(\C[v^{1/2}];v)\restarrow$ and its specializations}

We begin by introducing some necessary notation. Throughout this section, let
$\lambda/\mu$ be a skew partition of  $n$. If $M$ is an
$H_{R;a}(\SG_n)$-module, where $R$ is a PID and $a \in R$, we write
\[
M\restarrow
=
M\restarrow^{H_{R;a}(\SG_n)}_{H_{R;a}(\SG_{n-1})}.
\]
A box $B \in \tyd(\lambda/\mu)$ is called \emph{removable} if
$\tyd(\lambda/\mu) \setminus \{B\}$ is a Young diagram. We denote by
$\sfremov(\lambda/\mu)$ the set of removable boxes of
$\tyd(\lambda/\mu)$.

In this section, we introduce a special  filtration of $S^{\lambda/\mu}(\C[v^{1/2}];v)\restarrow$ and investigate their images under the specialization functors $\calF_{\iota}$ and $\calF_{\theta_0}$.

\begin{lemma}\label{lem: skew Specht submodule n-1}
Let $B \in \sfremov(\lambda/\mu)$.  Then the $\C[v^{1/2}]$-submodule of $S^{\lambda/\mu}(\C[v^{1/2}];v)$ spanned by $\{\sfh_{\ttS} \mid \ttS \in \SYT(\lambda/\mu)^B_{\boxn} \}$
is closed under the $H_{\C[v^{1/2}];v}(\SG_{n-1})$-action. 
\end{lemma}

\begin{proof}
For simplicity, we denote by $\calV$  the $\C[v^{1/2}]$-span of $\{\sfh_{\ttS} \mid \ttS \in \SYT(\lambda/\mu)^B_{\boxn} \}$.
Let $1 \le k < n-1$ and $\ttS \in \SYT(\lambda/\mu)^B_{\boxn}$.
For the assertion, we have to prove that $T_k \sfh_{\ttS} \in \calV$. 
To do this, we use several facts established in the proof of \cref{thm: specialization of skew Specht module}(3).
We have  three cases.

{\bf Case 1:} $k$ and $k+1$ are in the same row of $\ttS$.

In this case, we have proven that $T_k \sfh_{\ttS} = - \sfh_{\ttS}$.
Thus, $T_k \sfh_{\ttS} \in \calV$.

{\bf Case 2:} $k$ and $k+1$ are in the same column of $\ttS$.

In this case, we have proven that $T_k \sfh_{\ttS} = \sfh_{s_k \cdot \ttS}$.
By \cref{lem: expansion of row-standard into standard basis},
\[
\sfh_{s_k \cdot \ttS}
= \sum_{\mathtt{U} \in \SYT(\lambda/\mu)^{\ttS^{-1}(n)}_{\boxn}}
f_{s_k \cdot \ttS, \mathtt{U}}(v) \sfh_{\mathtt{U}},
\]
where every $f_{s_k \cdot \ttS, \mathtt{U}}(v)$ is in $\Z[v]$.
Since $\ttS^{-1}(n)$ is weakly below $B$, we have $T_k \sfh_{\ttS} \in \calV$.

{\bf Case 3:} $k$ and $k+1$ are in the distinct rows and columns of $\ttS$.

In this case, we have shown that
\[
T_k \sfh_{\ttS} =
\begin{cases}
\sfh_{s_k \cdot \ttS} 
& \text{if $k+1$ is below $k$,} \\
(v-1)\sfh_\ttS + v \sfh_{s_k \cdot \ttS}
& \text{if $k+1$ is above $k$.}
\end{cases}
\]
Since $1 \le k < n-1$, $\ttS^{-1}(n) = (s_k \cdot \ttS)^{-1}(n)$, and this box is weakly below $B$.
Thus, $T_k \sfh_{\ttS} \in \calV$.
\end{proof}

For each $B \in \sfremov(\lambda/\mu)$, we denote the
$H_{\C[v^{1/2}];v}(\SG_{n-1})$-module obtained in
\cref{lem: skew Specht submodule n-1} by
$S_B$.

Define a total order $\boxleq$ on $\sfremov(\lambda/\mu)$ by
\begin{align}\label{eq: total order condition}
B \boxleq B' \quad \text{if $B$ is weakly below $B'$.}
\end{align}
We enumerate the elements of $\sfremov(\lambda/\mu)$ as $B_1,B_2,\ldots,B_k$ in such a way that
\[
B_1 \boxl B_2 \boxl \cdots \boxl B_k .
\]
By the definition of $\SYT(\lambda/\mu)^B_{\boxn}$ given in \cref{def of SYTB}, we obtain the following filtration:
\begin{align}\label{eq: first remov filt}
S_{\bullet} \; := \;\left(
0 =:S_{B_0} \subset S_{B_1} \subset \cdots 
\subset S_{B_k} = S^{\lambda/\mu}(\C[v^{1/2}];v)\restarrow\right)
\end{align}
We refer to this as the \emph{removable-corner filtration} of $S^{\lambda/\mu}(\C[v^{1/2}];v)\restarrow$.

In the following, we describe 
the filtrations 
$\calF_\iota(S_{\bullet})$ and $\calF_{\theta_0}(S_{\bullet})$.
Before doing this, we introduce
certain filtrations  of $S^{\lambda/\mu}(\C(v^{1/2});v) \restarrow$ and $X_{\lambda/\mu} \restarrow$, constructed analogously to $S_{\bullet}$. 
For each $B \in \sfremov(\lambda/\mu)$, consider the $\C(v^{1/2})$-subspace of $S^{\lambda/\mu}(\C(v^{1/2});v)$ spanned by $\{\sfh_{\ttS} \mid \ttS \in \SYT(\lambda/\mu)^B_{\boxn} \}$.
Following the proof of \cref{lem: skew Specht submodule n-1}, we can see that this subspace  is closed under the $H_{\C(v^{1/2});v}(\SG_{n-1})$-action.
We denote the  resulting $H_{\C(v^{1/2});v}(\SG_{n-1})$-module
 by $\widetilde{S}_B$.
Then, we obtain the filtration
\begin{align}\label{eq: first remov filt}
\tilS_{\bullet} \; := \;\left(
0 =:\tilS_{B_0} \subset \tilS_{B_1} \subset \cdots 
\subset \tilS_{B_k} = S^{\lambda/\mu}(\C(v^{1/2});v)\restarrow\right),
\end{align}
called the
\emph{removable-corner filtration} of $S^{\lambda/\mu}(\C(v^{1/2});v)\restarrow$.

Next, for each $B \in \sfremov(\lambda/\mu)$,  consider the $\C$-subspace of $X_{\lambda/\mu}$ spanned by $\SYT(\lambda/\mu)^{B}_{\boxn}$.
By the definition of $\compH{0}$-action on $X_{\lambda/\mu}$ in \cref{eq: left action for X}, 
this space is closed under the $H_{\C;0}(\SG_{n-1})$-action.
We denote the  resulting $H_{\C;0}(\SG_{n-1})$-module
 by $X_B$.
Then, we obtain the filtration
\begin{align}\label{eq: first remov filt}
X_{\bullet} \; := \;\left(
0 =:X_{B_0} \subset X_{B_1} \subset \cdots 
\subset X_{B_k} = X_{\lambda/\mu}\restarrow\right),
\end{align}
called the
\emph{removable-corner filtration} of $X_{\lambda/\mu}\restarrow$.

Given two $H_{R;a}(\SG_{n-1})$-module filtrations 
\begin{align*}
0 \subsetneq M_1 \subsetneq \cdots 
\subsetneq M_r = M 
\quad\text{and}\quad 
0  \subsetneq N_1 \subsetneq \cdots 
\subsetneq N_s = N,    
\end{align*}
we say that 
these filtrations are  \emph{isomorphic} if 
(i) $r = s$ and (ii) there exists an $H_{R;a}(\SG_{n-1})$-module isomorphism $f: M \ra N$ such that $f(M_i) = N_i$ for all $1 \leq i \leq r$.

We are now ready to state the first main result. 
\begin{theorem}\label{thm: skew Specht poly rest and specialization}
\begin{enumerate}[label = {\rm (\arabic*)}]
\item $\calF_{\iota}(S_{\bullet})$ is isomorphic to $\tilS_{\bullet}$.
\item $\calF_{\theta_0}(S_{\bullet})$ is isomorphic to $X_{\bullet}$.
\item For $1 \leq i \leq k$, we have an $H_{\C(v^{1/2});v}(\SG_{n-1})$-module isomorphism
\[ 
\tilS_{B_i} \cong \bigoplus_{1\le j \le i}
S^{(\lambda \setminus \{B_j\})/\mu}(\C(v^{1/2});v).
\]

\item For $1 \leq i \leq k$, we have an $H_{\C;0}(\SG_{n-1})$-module isomorphism
\[
X_{B_i} \cong 
\bigoplus_{1 \le j \le i}
X_{(\lambda\setminus\{B_j\})/\mu}.
\]
\end{enumerate}
\end{theorem}

\begin{proof}
(1) Consider the $\genericH$-module isomorphism $\varphi: \calF_{\iota}(S^{\lambda/\mu}(\C[v^{1/2}];v)) \ra S^{\lambda/\mu}(\C(v^{1/2});v)$ given in the proof of \cref{thm: specialization of skew Specht module}(1).
Since $\varphi(1 \otimes \sfh_{\ttS}) = \sfh_{\ttS}$ for each  $\ttS \in \SYT(\lambda/\mu)$, it is an $H_{\C(v^{1/2});v}(\SG_{n-1})$-module isomorphism such that
\[
\varphi(\calF_\iota(S_{B_i})) = \tilS_{B_i} \quad\text{for all $1 \leq i \leq k$}.
\]
Therefore, $\calF_{\iota}(S_{\bullet})$ is isomorphic to $\tilS_{\bullet}$.

(2) Consider the $\compH{0}$-module isomorphism $\psi: \calF_{\theta_0}(S^{\lambda/\mu}(\C[v^{1/2}];v)) \ra X_{\lambda/\mu}$ given in the proof of \cref{thm: specialization of skew Specht module}(3).
Since $\psi(1 \otimes \sfh_{\ttS}) = \ttS$ for each $\ttS \in \SYT(\lambda/\mu)$, it is an $H_{\C;0}(\SG_{n-1})$-module isomorphism such that  
\[
\psi(\calF_{\theta_0}(S_{B_i})) = X_{B_i} \quad\text{for all $1 \leq i \leq k$.}
\]
Therefore, $\calF_{\theta_0}(S_{\bullet})$ is isomorphic to $X_{\bullet}$.

(4) For $1 \leq j \leq k$, let $M_j$ be the $\C$-subspace of
$X_{\lambda/\mu}\restarrow$ spanned by
$\{\ttS \in \SYT(\lambda/\mu) \mid \ttS(B_j)=n\}$.
From \cref{eq: left action for X}, we deduce that each $M_j$ is an
$H_{\C;0}(\SG_{n-1})$-submodule of
$X_{\lambda/\mu}\restarrow$ and is isomorphic to
$X_{(\lambda\setminus\{B_j\})/\mu}$.
Therefore, for $1 \leq i \leq k$, 
\[
X_{B_i}
=
\bigoplus_{1 \leq j \leq i} M_j
\cong
\bigoplus_{1 \leq j \leq i}
X_{(\lambda\setminus\{B_j\})/\mu}
\]
as $H_{\C;0}(\SG_{n-1})$-modules.

(3)
Let $1 \leq i \leq k$.
Since $\genericH$-$\module$ is semisimple, it suffices to show that 
\[
\ch^{\rm Frob}([\tilS_{B_i}]) = \sum_{1 \leq j\leq i} s_{(\lambda \setminus \{B_j\})/\mu}.
\]
Applying (1), \cref{thm: characteristic pres prop spec}, and (2) successively gives 
\begin{align}\label{eq: frob - frob equality}
\ch^{\rm Frob}([\tilS_{B_i}]) = 
\ch^{\rm Frob}([\calF_{\iota}(S_{B_i})]) = \ch([\calF_{\theta_0}(S_{B_i})]) = \ch([X_{B_i}]).
\end{align}
Therefore,
\[
\ch^{\rm Frob}([\tilS_{B_{i}}]) 
\overset{\cref{eq: frob - frob equality}}{=} \ch([X_{B_i}]) 
\overset{(4)}{=} \sum_{1 \leq j \leq i} 
\ch([X_{(\lambda \setminus \{B_j\})/\mu}])
= \sum_{1 \leq j \leq i} s_{(\lambda \setminus \{B_j\})/\mu}. \qedhere
\]
\end{proof}

From now on, we assume that $\lambda/\mu$ is an ordinary partition, that is, $\mu=\emptyset$, and let $B\in\sfremov(\lambda)$. 
In what follows, we investigate the cyclicity of $\tilS_B$ and $X_B$, together with their cyclic generators.
By \cref{thm: skew Specht poly rest and specialization}(3), we have
\[
\tilS_{B}  \cong \bigoplus_{B' \leq_{\lambda} B} S^{\lambda \setminus \{B'\}}(\C(v^{1/2});v).
\]
Since 
$\tilS_{B}$ is isomorphic to a direct sum of pairwise nonisomorphic irreducible modules, it is a cyclic $H_{\C(v^{1/2});v}(\SG_{n-1})$-module.
However, as the isomorphism is not explicit, identifying a cyclic generator within $\tilS_{B}$ is nontrivial.
In contrast to $\tilS_B$, the situation for $X_B$ is somewhat different. By \cref{thm: skew Specht poly rest and specialization}(4), we have
\[
X_B \cong \bigoplus_{B' \leq_\lambda B} X_{\lambda\setminus\{B'\}}.
\]
Unlike the decomposition of $\tilS_B$ into pairwise nonisomorphic irreducible modules, the right-hand side is a direct sum of indecomposable $H_{\C;0}(\SG_{n-1})$-modules. Hence, this decomposition alone does not imply that 
$X_B$ is cyclic. 

Our aim is to construct cyclic generators for $\tilS_B$ and $X_B$, respectively.
We first introduce notation for describing these generators.

\begin{enumerate}[label = {\rm (\roman*)}, leftmargin = 4ex, itemsep = 1ex]
\item 
Let $\alpha=(\alpha_1,\ldots,\alpha_l)\models n$ with
$\lambda(\alpha)=\lambda$.
Consider the left-justified array whose $i$th row from the top consists of $\alpha_i$ boxes.
Fill its boxes with $1,2,\ldots,n$ from left to right starting with the top row.
Let $\ttT_\alpha$ be the tableau obtained from this filling by sliding the boxes in each column upward until no empty positions remain above any box.
For instance, if $\lambda = (4,2,1)$ and $\alpha = (1,4,2)$, $\ttT_\alpha$ is given by
\[
\ttT_\alpha = \begin{array}{l}
\begin{ytableau}
    1 & 3& 4 & 5\\
    2 & 7 \\
    6
\end{ytableau}
\end{array}.
\]
We remark that if $\alpha$ is a partition, then the tableau $\ttT_\alpha$ defined above coincides with the tableau $\ttT_\lambda$ defined in \cref{subsec: comp and diag}.
Using the tableau $\ttT_\alpha$, we define an element $\bfz_\alpha \in \genericH$ by
\[
\bfz_\alpha := \bfy_\alpha T_{\wcol(\ttT_\alpha)}\bfx_{\lambda(\alpha)'}.
\]
This element is a variant of the element $z_\alpha$ introduced in \cite[Section 4]{86DJ}.

\item
For $B \in \sfremov(\lambda)$, let $\ttT_{B} \in \SYT(\lambda)$ be the standard Young tableau obtained from $\ttT_{\lambda\setminus\{B\}}$ by adding the box $B$ and filling it with $n$.
Define the element $\ttg_{B}$ of $X_{B}$ by 
\[
\ttg_{B} := \sum_{B' \leq_{\lambda} B} \ttT_{B'}.
\]
\end{enumerate}

With this preparation, we can state the second main result.

\begin{theorem}\label{thm: generator of S and X as n-1}{\rm (This will be proved in \cref{subsec: proof of main theorem}.)}
Let $\lambda$ be a partition of $n$ and let $B \in \sfremov(\lambda)$. 
\begin{enumerate}[label = {\rm (\arabic*)}]
\item 
Let $\alpha$ be any composition of $n$ such that $\lambda(\alpha) = \lambda$ and $\alpha_{\ell(\lambda)}$ is equal to the column index of $B$. 
As an $H_{\C(v^{1/2});v}(\SG_{n-1})$-module, $\tilS_{B}$ is cyclically generated by $\bfz_\alpha$.
\item
As an $H_{\C;0}(\SG_{n-1})$-module, $X_{B}$  is cyclically generated by $\ttg_{B}$.
\end{enumerate}
\end{theorem}

\begin{remark}\label{rem: removable component is not cyclic}
\cref{thm: generator of S and X as n-1} does not hold for arbitrary skew partitions.
For example, consider $\lambda/\mu=(5,1)/(1)$ and $B = (1,5) \in \sfremov(\lambda/\mu)$.
By \cref{thm: skew Specht poly rest and specialization}(3), we have 
\[
\tilS_B 
\cong
S^{(4)} \oplus
S^{(4,1)/(1)} 
\cong
(S^{(4)})^{\oplus 2}\oplus S^{(3,1)}.
\]
Since 
\[
H_{\C(v^{1/2});v}(\SG_4) \cong \bigoplus_{\lambda \vdash 4} (S^{\lambda})^{\oplus |\SYT(\lambda)|}
\]  
and $|\SYT((4))| = 1$, 
$\tilS_B$ is not a quotient of $H_{\C(v^{1/2});v}(\SG_4)$.
Therefore $\tilS_B$ is not cyclic.
Similarly, 
by \cref{thm: skew Specht poly rest and specialization}(4), we have 
\[
X_B \cong X_{(4)} \oplus X_{(4,1)/(1)}.
\]
By \cite[Remark A.1(2)]{23KLO} and \cite[Theorem 3.3]{16Huang}, $X_{(4)} \cong \bfP_{(4)}$ and $X_{(4,1)/(1)} \cong \bfP_{(4)} \oplus \bfP_{(3,1)}$.
Here, $\bfP_{\alpha}$ is the projective indecomposable module attached to $\alpha \models 4$, defined in \cite{16Huang}.
It follows that
\[
X_B \cong  (\bfP_{(4)})^{\oplus 2} \oplus \bfP_{(3,1)}.
\]
Since
\[
H_{\C;0}(\SG_4) \cong  \bigoplus_{\alpha \models 4} \bfP_{\alpha},
\]
the module $X_B$ is not a quotient  of $H_{\C;0}(\SG_4)$.
Therefore, $X_B$ is not cyclic.
\end{remark}

\subsection{The proof of \cref{thm: generator of S and X as n-1}}
\label{subsec: proof of main theorem}
For the proof, we introduce necessary lemmas. 
The first lemma follows from \cite[Lemma~4.1]{86DJ} by applying the anti-involution of $\genericH$ that fixes $T_i$ for $1 \leq i \leq n-1$.

\begin{lemma}{\rm(cf. \cite[Lemma 4.1]{86DJ})}\label{lem: Dipper--James}
\begin{enumerate}[label = {\rm (\arabic*)}]
\item For $\alpha \models n$,  $\bfz_{\alpha} \neq 0$.
\item Let $\alpha, \beta \models n$ with $\lambda(\alpha) = \lambda(\beta)'$ and $w \in \SG_n$.
If  $\bfy_\alpha T_w \bfx_\beta \neq 0$, then  
$\bfy_\alpha T_w \bfx_\beta = \pm v^i \bfz_\alpha$ for some $i \ge 0$.
\end{enumerate}
\end{lemma}

Let $\alpha,\beta\models n$ with $\lambda(\alpha) = \lambda(\beta)$. 
Consider $n$ consecutive boxes, partitioned from left to right into blocks of sizes $\alpha_1, \alpha_2, \ldots, \alpha_{\ell(\alpha)}$. 
We fill these boxes with the elements of $[n]$ as follows. 
For $1 \leq i \leq \ell(\beta)$, successively fill the leftmost unfilled block of size $\beta_i$ with the $\beta_i$ smallest unused elements of $[n]$, arranged in increasing order from left to right. Since $\lambda(\alpha)=\lambda(\beta)$, this procedure fills all the boxes. 
Define $d_{\alpha,\beta}\in\SG_n$ to be the permutation whose one-line notation is the resulting sequence of entries.
By construction of $d_{\alpha, \beta}$, we have
\begin{align}\label{eq: conjugation of Young subgrp}
\SG_{\beta} = d_{\alpha, \beta} \SG_{\alpha} d_{\alpha, \beta}^{-1}.
\end{align}
For example, if $\alpha = (3,2,2)$ and $\beta = (2,2,3)$, then to obtain $d_{\alpha, \beta}$, we first partition $7$ consecutive boxes according to $\alpha$ as follows:
\[
\boxed{\phantom{5}} \boxed{\phantom{6}} \boxed{\phantom{7}} \mid \boxed{\phantom{1}} \boxed{\phantom{2}} \mid \boxed{\phantom{3}} \boxed{\phantom{4}}.
\]
We then fill these blocks according to $\beta$. First, we fill the leftmost unfilled block of size $2$ with $1,2$. We then fill the next unfilled block of size $2$ with $3,4$, and finally fill the remaining block of size $3$ with $5,6,7$. The resulting filling is
\[
\boxed{5} \boxed{6} \boxed{7} \mid \boxed{1} \boxed{2} \mid \boxed{3} \boxed{4}.
\]
Thus $d_{\alpha,\beta} = 5671234$.%$d_{\alpha,\beta}^{-1} = 4567123$.
%Similarly, we obtain $d_{\beta,\alpha} = 4567123$.

\begin{lemma}\label{lem: T y and y T}
Let $\alpha, \beta \models n$ with $\lambda(\alpha) = \lambda(\beta)$.
Then,
\begin{align*}
T_{d_{\alpha,\beta}}\bfy_\alpha
    =
    \bfy_\beta T_{d_{\alpha,\beta}}.
\end{align*}
\end{lemma}
\begin{proof}
Since 
$d_{\alpha, \beta}(i) < d_{\alpha, \beta}(i+1)$ for all $i \in \set(\alpha)^\rmc$, we have 
\begin{align}\label{eq: d al be property one}
\ell(d_{\alpha, \beta} \sigma) = \ell(d_{\alpha, \beta}) + \ell(\sigma) \quad\text{for all $\sigma \in \SG_{\alpha}$}.
\end{align}
Similarly, since $d_{\alpha, \beta}^{-1}(i) < d_{\alpha, \beta}^{-1}(i+1)$ for all $i \in \set(\beta)^\rmc$, we have
\begin{align}\label{eq: d al be property two}
\ell(\sigma d_{\alpha, \beta}) = \ell(\sigma) + \ell(d_{\alpha,\beta}) \quad\text{for all $\sigma \in \SG_{\beta}$}.
\end{align}
Therefore, 
\begin{align*}
   T_{d_{\alpha, \beta}} \bfy_{\alpha} 
   &= \sum_{\sigma \in \SG_{\alpha}} \sgn(\sigma)v^{\ell(w_0(\alpha)) - \ell(\sigma)}  T_{d_{\alpha, \beta}}T_{\sigma}\\
   &= \sum_{\rho \in d_{\alpha, \beta}\SG_{\alpha}} \sgn(d_{\alpha, \beta})\sgn(\rho)v^{\ell(w_0(\alpha)) - \ell(d_{\alpha,\beta}^{-1} \rho)} T_{\rho} & \text{(by \cref{eq: d al be property one})}\\
   &= \sum_{\gamma \in \SG_{\beta}d_{\alpha, \beta}} \sgn(d_{\alpha, \beta})\sgn(\gamma)v^{\ell(w_0(\alpha)) - \ell(\gamma d_{\alpha, \beta}^{-1})} T_{\gamma} & \text{(by \cref{eq: conjugation of Young subgrp})}\\
   &=  \sum_{\sigma \in \SG_{\beta}} \sgn(\sigma)v^{\ell(w_0(\alpha)) - \ell(\sigma)} T_{\sigma}T_{d_{\alpha, \beta}} & \text{(by \cref{eq: d al be property two})}\\
   &= \bfy_{\beta} T_{d_{\alpha,\beta}}.
\end{align*}
Here, the last equality follows from   $\ell(w_0(\alpha)) = \ell(w_0(\beta))$.
%(see, for example, \cite[p.~33]{99Mat}).
\end{proof}

\begin{lemma}\label{lem: basic step for inclusion}
Let $\alpha, \beta \models n$ with $\lambda(\alpha) = \lambda(\beta)$ and let $l := \ell(\alpha)$.
Then, we have the following.
\begin{enumerate}[label = {\rm (\arabic*)}]
    \item $H_{\C( v^{1/2});v}(\SG_n) \bfv_\alpha = H_{\C( v^{1/2});v}(\SG_n) \bfv_\beta$.
    \item If $\alpha_l = \beta_l$, then $ H_{\C( v^{1/2});v}(\SG_{n-1}) \bfv_\alpha = H_{\C( v^{1/2});v}(\SG_{n-1}) \bfv_\beta$.
    \item If $\alpha_l >\beta_l$, then 
    $\bfz_\beta \in H_{\C(v^{1/2});v}(\SG_{n-1}) \bfz_\alpha$.
\end{enumerate}
\end{lemma}
\begin{proof}
For simplicity, set $\lambda := \lambda(\alpha)$.

(1) To prove the assertion, it suffices to show that $\bfz_{\beta} \in \genericH\bfz_{\alpha}$.
Write 
\[
T_{d_{\alpha, \beta}} T_{\wcol(\ttT_{\alpha})} = \sum_{w \in \SG_n} a_w T_w \quad(a_w \in \Z[v]).
\]
Then, 
we have 
\begin{align*}
T_{d_{\alpha, \beta}} \bfz_\alpha &= 
T_{d_{\alpha, \beta}} \bfy_{\alpha} T_{\wcol(\ttT_\alpha)} \bfx_{\lambda'} \\
%\quad(\text{since $T_{d_{\alpha, \beta}}$ is invertible.})\\
&= \bfy_{\beta} T_{d_{\alpha, \beta}} T_{\wcol(\ttT_\alpha)} \bfx_{\lambda'} \quad(\text{by \cref{lem: T y and y T}})   \\
&= \sum_{w \in \SG_n} a_w \bfy_{\beta} T_{w} \bfx_{\lambda'}.
\end{align*}
Applying \cref{lem: Dipper--James}(2) to each term
$\bfy_\beta T_w\bfx_{\lambda'}$ in the above equation, we derive that 
\[
\sum_{w \in \SG_n} a_w \bfy_{\beta} T_{w} \bfx_{\lambda'} = c \bfz_{\beta}
\]
for some $c \in \Z[v]$.
Since $T_{d_{\alpha, \beta}} \bfz_{\alpha}$ is nonzero,  $c \neq 0$ and thus
\begin{align}\label{eq: bfz beta are in H bfz al}
\bfz_\beta = c^{-1} T_{d_{\alpha,\beta}} \bfz_\alpha \in \genericH \bfz_\alpha.
\end{align}

(2) Suppose that $\alpha_l = \beta_l$.
Then $d_{\alpha, \beta}(n) = n$, and thus $d_{\alpha, \beta}$ can be viewed as an element of the subgroup $\SG_{n-1}$ of $\SG_n$.
By \cref{eq: bfz beta are in H bfz al}, $\bfz_{\beta} \in H_{\C(v^{1/2});v}(\SG_{n-1}) \bfz_\alpha$.
By switching the role of $\alpha$ and $\beta$, we have $\bfz_{\alpha} \in H_{\C(v^{1/2});v}(\SG_{n-1}) \bfz_\beta$.

(3) 
Suppose that $\alpha_l > \beta_l$ and set $x := \alpha_l$ and $y := \beta_l$.
Let $(c_1, c_2, \ldots, c_{l-2})$ be the partition obtained from $\lambda$ by deleting one occurrence of the part $x$ and one occurrence of the part $y$.
%Let $\tilal = (c_1, c_2, \ldots, c_{l-2}, y,x)$ and  $\tilbe = (c_1, c_2, \ldots, c_{l-2}, x,y)$.
By (2), we have 
\begin{align*}
& H_{\C(v^{1/2});v}(\SG_{n-1})\bfz_{\alpha} = H_{\C(v^{1/2});v}(\SG_{n-1})\bfz_{(c_1, c_2, \ldots, c_{l-2}, y,x)}\quad\text{and}\\
& H_{\C(v^{1/2});v}(\SG_{n-1})\bfz_{\beta} = H_{\C(v^{1/2});v}(\SG_{n-1})\bfz_{(c_1, c_2, \ldots, c_{l-2}, x,y)}.
\end{align*}
Thus, we may assume that $\alpha = (c_1, c_2, \ldots, c_{l-2}, y,x)$ and $\beta = (c_1, c_2, \ldots, c_{l-2}, x,y)$.

From now on, we find an element $h \in H_{\C(v^{1/2});v}(\SG_{n-1})$ such that $\bfz_{\beta} = h\bfz_{\alpha}$.
To do this, let $\gamma = (c_1, c_2, \ldots, c_{l-2}, y, x - y, y) \models n$.
Let $\calD_L$ be the set of minimal left coset representatives of $\SG_\beta/\SG_\gamma$ and $\calD_R$ be the set of minimal right coset representatives of $\SG_\gamma\backslash\SG_\alpha$.
Then, we have
\begin{align}\label{eq: alpha beta gamma decomp}
\bfy_{\alpha} = \bfy_\gamma \bfy_{\calD_R} \quad \text{and} \quad \bfy_\beta = \bfy_{\calD_L}  \bfy_\gamma.
\end{align}
Since $\beta_l = y = \gamma_{l+1}$, every permutation in $\calD_L$ fixes $n$, and thus  $\bfy_{\calD_L} \in H_{\C( v^{1/2});v}(\SG_{n-1})$.
We claim that 
\begin{align}\label{eq: y v al equal v be}
\bfv_\beta = \pm \bfy_{\calD_L} \cdot \bfv_\alpha.
\end{align}

{\it Step 1.} We express $\bfy_{\calD_L} \cdot \bfz_\alpha$ as a linear combination of $\bfy_\beta T_{d \wcol(\ttT_\alpha)} \bfx_{\lambda'}$'s, where $d$ ranges over $\calD_R$.
Note that 
\[
\calD_R
=
\SG_{\alpha} \cap
\{\sigma \in \SG_n \mid \emptyset \subseteq \Des_L(\sigma) \subseteq \set(\gamma)\}.
\]
More explicitly, $\calD_R$ is the set of permutations $d \in \SG_n$ satisfying the following:
\begin{enumerate}[label = {\bf D\arabic*.}]
\item $d(i) = i$ for all $1 \leq i \leq n-x$,
\item The one-line notation of $d$ contains 
\[
\sfw_1 := n-x+1 \ \ n-x+2 \ \cdots \ n- y\quad\text{and}\quad \sfw_2 := n-y+1 \ \ n-y+2 \ \cdots\ n
\]
as subwords.
\end{enumerate}
It follows that  
$$\calD_R =   [\id,d_0]_R,$$ 
where $d_0$ is the permutation whose one-line notation is obtained by concatenating the words $1\, 2\, \ldots\, n-x$, $\sfw_2$, and $\sfw_1$ in order.
Also, we note that the one-line notation of $\wcol(\ttT_\alpha)$ contains the concatenation $\sfw_1 \sfw_2$ as a subword.
Combining this with \textbf{D1} yields that 
\begin{align}\label{eq: length cond on branching gen}
\ell(d \wcol(\ttT_\alpha)) = \ell(d) + \ell(\wcol(\ttT_\alpha)) \quad \text{for all $d \in \calD_R$.}
\end{align}
Hence,  
\begin{align*}
\bfy_{\calD_L} \bfz_\alpha 
&= \bfy_{\calD_L} \bfy_\alpha T_{\wcol(\ttT_\alpha)} \bfx_{\lambda'} \\
&= \bfy_{\calD_L} \bfy_\gamma  \bfy_{\calD_R} T_{\wcol(\ttT_\alpha)}  \bfx_{\lambda'} &(\text{by \cref{eq: alpha beta gamma decomp}})\\
&= \bfy_{\beta} \bfy_{\calD_R} T_{\wcol(\ttT_\alpha)}  \bfx_{\lambda'} &(\text{by \cref{eq: alpha beta gamma decomp}})\\
&= 
\bfy_{\beta} 
\left(
\sum_{d \in [\id, d_0]_R} \sgn(d) v^{\ell(d_0) - \ell(d)}T_{d}
\right)
T_{\wcol(\ttT_\alpha)} \bfx_{\lambda'}
\\
&=  \sum_{d \in [\id, d_0]_R} \sgn(d) v^{\ell(d_0) - \ell(d)} \bfy_\beta T_{d \wcol(\ttT_\alpha)} \bfx_{\lambda'} &(\text{by \cref{eq: length cond on branching gen}}).
\end{align*}

{\it Step 2.} Let $d \in [\id, d_0]_R$.
We show that  
\[
\bfy_\beta T_{d \wcol(\ttT_\alpha)} \bfx_{\lambda'} =  
\begin{cases}
   \bfv_\beta  & \text{if $d = d_0$,}\\
   0  & \text{otherwise.}
\end{cases}
\]
For each $1 \leq j \leq x$, we denote $B_j$ the bottommost box in the $j$th column of $\tyd(\lambda)$.
Since $\alpha_i = \beta_i$ for $1 \leq i\leq l-2$, 
$\ttT_{\alpha}$ and $\ttT_{\beta}$ coincide on every box in $\tyd(\lambda)$ except the boxes $B_j$.
In addition, for all $1 \leq j \leq x$, we have
\begin{align}\label{eq: ttT al and ttT be entry compare}
\ttT_{\alpha}(B_j) = n-x+j 
\quad\text{and} \quad
\ttT_{\beta}(B_j) = 
\begin{cases}
    n-y+j & \text{if $1 \leq j \leq y$,}\\
    n-x-y+j & \text{if $y+1 \leq j\leq x$}.
\end{cases}
\end{align}
Therefore, if $d = d_0$, then $d \cdot \ttT_{\alpha} = \ttT_{\beta}$, where $d \cdot \ttT_{\alpha}$ is the bijective tableau of shape $\lambda$ defined by $d \cdot \ttT_{\alpha}(B) = d(\ttT_{\alpha}(B))$ for all $B \in \tyd(\lambda)$.
This implies that
\[
\bfy_{\beta}T_{d \wcol(\ttT_\alpha)} = \bfy_{\beta}T_{ \wcol(d \cdot \ttT_\alpha)}\bfx_{\lambda'} = 
\bfy_{\beta}T_{ \wcol(\ttT_\beta)} \bfx_{\lambda'} = \bfz_{\beta}.
\]

Next, consider the case where $d \neq d_0$.
Set 
\[
j_0 = d^{-1}(n-x+1) - (n - x).
\]
From \textbf{D1} and \textbf{D2}, we derive that $n-x+1 \leq d^{-1}(n-x+1) \leq n-x+y$,
and thus
\begin{align}\label{eq: j0 and y ineq}
1 \le j_0 \le y.
\end{align}

By \cref{eq: ttT al and ttT be entry compare}, $\ttT_{\alpha}(B_{j_0}) = n - x + j_0$.
Moreover, since $d(n-x+j_0) = n-x+1$, we have $d \cdot \ttT_{\alpha}(B_{j_0}) = n - x +1$.
Set 
\[
m = \sum_{1 \leq j\leq j_0}\lambda'_j -1.
\]
Since $j_0 \le \alpha_{l-1} < \alpha_l$, we have $\lambda'_{j_0} \ge 2$.
It follows that $s_m \in \SG_{\lambda'}$, and hence 
\begin{align}\label{eq: x al is decomposed as sm and other sum}
\bfx_{\lambda'} = 
(1+T_m) \left(\sum_{\substack{\sigma \in \SG_{\lambda'} \\ m \notin \Des_L(\sigma)}}T_\sigma\right).
\end{align}
Therefore, to prove $\bfy_{\beta} T_{d \wcol(\ttT_\alpha)} \bfx_{\lambda'} = 0$, it suffices to show that $\bfy_{\beta} T_{d \wcol(\ttT_\alpha)}(1+T_m) = 0$.

Let $B_{j_0}'$ be the box in $\tyd(\lambda)$ immediately above
$B_{j_0}$.
Since $\alpha_{l-1} = y$ and $\ttT_{\alpha}(B_{j_0}) = n - x + j_0$, we have $\ttT_{\alpha}(B_{j_0}') = n - x - y + j_0$. 
Together with \textbf{D1} and \cref{eq: j0 and y ineq}, this gives $d \cdot \ttT_{\alpha}(B'_{j_0}) = n- x- y +j_0$.
Hence,
\begin{align}\label{eq: w col d T m}
\wcol(d \cdot \ttT_{\alpha})(m)
= n-x-y+j_0
\quad \text{and} \quad
\wcol(d \cdot \ttT_{\alpha})(m+1)
= n-x+1.
\end{align}
By \cref{eq: j0 and y ineq} and $y < x$, we have 
\[
n-x-y+j_0,\;  n-x+1 \in [n-x-y+1, \ n-y],
\]
which implies that the transposition $\tau := (n-x-y+j_0 \ \ \  n-x+1)$ is in $\SG_\beta$.
It follows that  $\tau \wcol(d \cdot \ttT_{\alpha})$ and $\wcol(d \cdot \ttT_{\alpha})$ are in the same right coset of $\SG_{\beta}$ in $\SG_n$. 
Then, using \cref{eq: Ti and bfy}, we deduce that
\begin{align}\label{eq: easy fact on bfy}
\bfy_{\beta} T_{\wcol(d \cdot \ttT_{\alpha})} = - \bfy_{\beta} T_{\tau \wcol(d \cdot \ttT_{\alpha})}. 
\end{align}
By \cref{eq: w col d T m}, $\wcol(d \cdot \ttT_{\alpha})s_m = \tau \wcol(d \cdot \ttT_{\alpha})$ and $m \notin \Des_R(\wcol(d \cdot \ttT_{\alpha}))$.
It follows that 
\begin{align}\label{eq: right mul perm left mul perm}
T_{\wcol(d \cdot \ttT_{\alpha})} T_m = 
T_{\tau \wcol(d \cdot \ttT_{\alpha})}.
\end{align}
Putting \cref{eq: easy fact on bfy} and \cref{eq: right mul perm left mul perm} together, we have
\begin{align*}
\bfy_{\beta} T_{ d \wcol(\ttT_\alpha)} (1+T_m) 
= \bfy_{\beta} T_{ \wcol(d \cdot \ttT_\alpha)} (1+T_m) 
= - \bfy_{\beta} T_{\tau \wcol(d \cdot \ttT_\alpha)} + \bfy_{\beta}T_{\tau \wcol(d \cdot \ttT_{\alpha})}
= 0.
\end{align*}

Consequently, the claim \cref{eq: y v al equal v be}
follows by combining 
{\it Step 1} and {\it Step 2}.
\end{proof}

We are now ready to prove \protect\cref{thm: generator of S and X as n-1}.

\begin{proof}[Proof of \protect\cref{thm: generator of S and X as n-1}.]
Set $l := \ell(\lambda)$.
Enumerate the elements of $\sfremov(\lambda)$ as 
$$ 
B_1 <_{\lambda} B_2 <_{\lambda} \cdots <_{\lambda} B_k 
$$
and write $B_i = (r_i, c_i)$ for all $1 \leq i \leq k$.

(1) For each $1 \leq i \leq k$,  let 
\[
\lambda^{\langle i \rangle} := (\lambda_1,  \ldots, {\lambda}_{r_i - 1}, {\lambda}_{r_i + 1}\ldots, \lambda_l, \lambda_{r_i}) \models n.
\]
By \cref{lem: basic step for inclusion}(2),  we have 
\begin{align*}
H_{\C(v^{1/2});v}(\SG_{n-1})\bfz_{\lambda^{\langle i \rangle}} = H_{\C(v^{1/2});v}(\SG_{n-1})\bfz_{\alpha} 
\end{align*}
for all  $\alpha \models n$ with $\lambda(\alpha) = \lambda$ and $\alpha_l = c_i$.
Thus, to prove the assertion, it suffices to show that 
\begin{align}\label{eq: H z lam i SBi}
\tilS_{B_i} = H_{\C(v^{1/2});v}(\SG_{n-1})\bfz_{\lambda^{\langle i \rangle}}.
\end{align}
We prove this by applying mathematical induction on $i$.
Recall that $\sfh_{\ttS} = T_{\wrow(\ttS)} \bfz_{\lambda}$ for  $\ttS \in \SYT(\lambda)$ and $\tilS_{B_i} = \C(v^{1/2})  \{\sfh_{\ttS} \mid \ttS \in \SYT(\lambda)^{B_i}_{\boxn} \}$.

If $i = 1$, then $\lambda^{\langle 1 \rangle} = \lambda$.
Since $\bfz_{\lambda} = \sfh_{\ttT_{\lambda}}  \in \tilS_{B_1}$, we have $\tilS_{B_1} \supseteq H_{\C(v^{1/2});v}(\SG_{n-1})\bfz_{\lambda}$.
The opposite inclusion also holds since $\wrow(\ttS)(n) = n$ for all $\ttS \in \SYT(\lambda)^{B_1}_{\boxn}$.

Now, let $1 < i \leq k$. 
Assume that $\tilS_{B_{i-1}} = H_{\C(v^{1/2});v}(\SG_{n-1}) \bfz_{\lambda^{\langle i-1 \rangle}}$.
We first show that $\tilS_{B_{i-1}}   \subseteq 
H_{\C(v^{1/2});v}(\SG_{n-1}) \bfz_{\lambda^{\langle i \rangle}}   \subseteq \tilS_{B_{i}}$.
The former inclusion  follows from  the assumption together with  \cref{lem: basic step for inclusion}(3).
To prove the latter inclusion, let  $\mathbb{T} \in \RT(\lambda)$ such that $\wrow(\mathbb{T}) = d_{\lambda,  \lambda^{\langle i \rangle}}$.
Since $\mathbb{T}(B_i) = n$, it follows from \cref{lem: expansion of row-standard into standard basis} that  
$\sfh_{\mathbb{T}} \in \tilS_{B_i}$.
By \cref{lem: T y and y T}, we have 
\[
\sfh_{\mathbb{T}} = T_{\wrow(\mathbb{T})}\bfy_{\lambda}T_{\wcol(\ttT_{\lambda})}\bfx_{\lambda'}  = T_{d_{\lambda,\lambda^{\langle i \rangle}}}\bfy_{\lambda}T_{\wcol(\ttT_{\lambda})}\bfx_{\lambda'} = \bfy_{\lambda^{\langle i \rangle}} T_{d_{\lambda,\lambda^{\langle i \rangle}}} T_{\wcol(\ttT_{\lambda})}\bfx_{\lambda'}.
\]
Since $\sfh_{\mathbb{T}} \neq 0$, we deduce from \cref{lem: Dipper--James}(2) that 
\begin{align}\label{eq: h and z scalar diff}
\sfh_{\mathbb{T}} = c \bfz_{\lambda^{\langle i \rangle}}
\quad \text{for some nonzero $c \in \Z[v]$,}
\end{align}
which implies that $\bfz_{\lambda^{\langle i \rangle}} \in \tilS_{B_i}$.

Recall that  $\tilS_{B_i}/\tilS_{B_{i-1}}$ is an irreducible $H_{\C(v^{1/2});v}(\SG_{n-1})$-module (see \cref{thm: skew Specht poly rest and specialization}(3)).
Hence, for the equality \cref{eq: H z lam i SBi}, we have only to show that $\bfz_{\lambda^{\langle i \rangle}} \notin \tilS_{B_{i-1}}$.
Note that 
$w_0(\lambda^{\langle i \rangle}) \cdot \ttT_{\lambda^{\langle i \rangle}}$ is a column-strict tableau of shape $\lambda$ and 
$i \notin \Des_L(\wcol(\ttT_{\lambda^{\langle i \rangle}}))$ for all $i \in \set({\lambda^{\langle i \rangle}})^\rmc$.
Therefore, for all $\sigma \in \SG_{\lambda^{\langle i \rangle}}$ and $\rho \in \SG_{\lambda'}$, we have 
\begin{align}\label{eq: length property of y al T x lam}
\ell(\sigma \wcol(\ttT_{\lambda^{\langle i \rangle}}) \rho) =  \ell(\sigma \wcol(\ttT_{\lambda^{\langle i \rangle}})) + \ell(\rho)=  \ell(\sigma) + \ell(\wcol(\ttT_{\lambda^{\langle i \rangle}})) + \ell(\rho).
\end{align}
Letting 
\begin{align}\label{eq: bfz first one}
\bfz_{\lambda^{\langle i \rangle}} =  \sum_{\sigma \in \SG_n} a_{\sigma} T_{\sigma} \quad(a_{\sigma} \in \C(v^{1/2})),
\end{align}
it follows from 
\cref{eq: length property of y al T x lam} that for $\sigma \in \SG_n$, 
\begin{align}\label{eq: a tau nonzer cond iff}
    a_{\sigma} \neq 0 \quad\text{if and only if}\quad\sigma \in \SG_{\lambda^{\langle i \rangle}}{\wcol(\ttT_{\lambda^{\langle i \rangle}})}\SG_{\lambda'}.
\end{align}
On the other hand, 
by virtue of   \cref{lem: basic step for inclusion}(1) and \cref{lem: skew specht basis},  we can write in the following form:
\begin{align}\label{eq: bfz second one}
\bfz_{\lambda^{\langle i \rangle}} = 
 \sum_{\ttS \in \SYT(\lambda)} b_{\ttS} \sfh_{\ttS} \quad(b_{\ttS} \in \C(v^{1/2})).
\end{align}
By \cref{lem: triangularity of elt in  B}, each $\sfh_{\ttS}$ is written as 
\begin{align}\label{eq: bfz third one}
v^{\ell(\wrow(\ttS))}
T_{\wcol(\ttS)} + \sum_{\substack{\sigma \in \SG_n \\ \ell(\sigma) > \ell(\wcol(\ttS))}} c_{\sigma }T_{\sigma} \quad(c_\sigma \in \C(v^{1/2})).
\end{align}
Let $l_0 := \min\{\ell(\wcol(\ttS)) \mid \ttS \in \SYT(\lambda) \text{ with }b_{\ttS} \neq 0 \}$.
Plugging \cref{eq: bfz third one} into the right hand side of  \cref{eq: bfz second one} and comparing it with \cref{eq: bfz first one}, we observe that 
\begin{enumerate}[label = $\bullet$]
    \item for $\ttS \in \SYT(\lambda)$ with  $\ell(\wcol(\ttS)) = l_0$,  $b_{\ttS} \neq 0$ if and only if
$a_{\wcol(\ttS)} \neq 0$, and 
    \item $\ell(\wcol(\ttT_{\lambda^{\langle i \rangle}})) = l_0$ by \cref{eq: length property of y al T x lam} and  \cref{eq: a tau nonzer cond iff}.
\end{enumerate}
In view of \cref{eq: a tau nonzer cond iff}, $a_{\wcol(\ttT_{\lambda^{\langle i \rangle}})} \neq 0$.
Combining this with the above observation, 
we have $b_{\ttT_{\lambda^{\langle i \rangle}}} \neq 0$.
We also need the fact  that $\sfh_{\ttT_{\lambda^{\langle i \rangle}}}\notin \tilS_{B_{i-1}}$, which follows from $\ttT_{\lambda^{\langle i \rangle}}(B_i) = n$ and \cref{lem: skew specht basis}. 
Consequently, $\bfz_{\lambda^{\langle i \rangle}} \notin \tilS_{B_{i-1}}$.

(2) 
We prove the assertion by mathematical induction on the subscript $i$ of $B_i$.
For $1 \leq i \leq k$, 
let $M_i$ be the cyclic submodule of $X_{B_i}$ generated by $\ttg_{B_i}$.

If $i = 1$, then $\ttg_{B_1} = \ttT_{B_1}$ and the assertion follows  from \cref{eq: left action for X}.

Now, let $1 < i  \leq k$ and assume that  $ X_{B_{i-1}} = M_{i-1}$.
The assertion can be accomplished by showing that $M_i \supseteq \SYT(\lambda)^{B_i}_{\boxn}$.
Let $x = \sum_{1 \leq j \leq r_i - 1} \lambda_j$.
For $1 \leq j \leq i$, we examine the boxes occupied by $x+1, x+2, \ldots, x+{\lambda_{r_i}}$ in $\ttT_{B_j}$.
If $1 \leq j < i$, then these  are precisely the boxes in the $r_i$th row of $\tyd(\lambda)$ as follows:
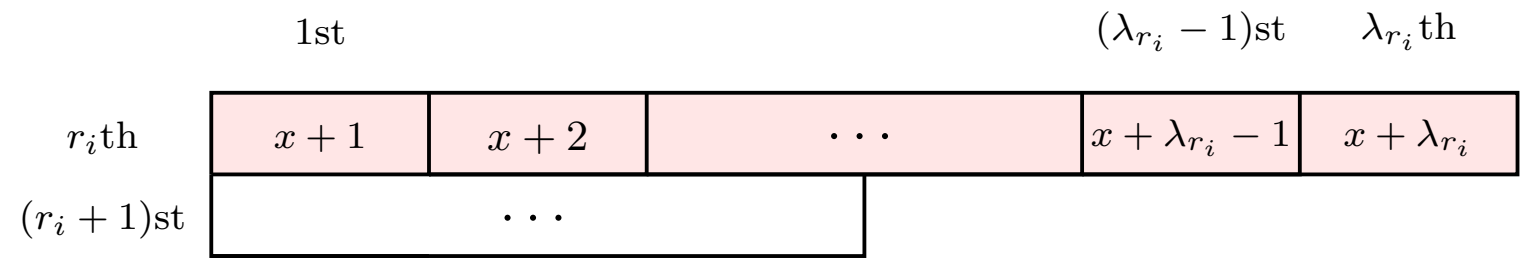
\[
\scalebox{0.9}{$
\def \hd {1.6}
\def \vd {0.6}
\begin{tikzpicture}
\node at (5.5*\hd, 2.75*\vd) {\tiny $1$st};
\node at (9.5*\hd, 2.75*\vd) {\tiny $(\lambda_{r_i} - 1)$st};
\node at (10.5*\hd, 2.75*\vd) {\tiny $\lambda_{r_i}$th};

\draw[thick, fill = red!10] (6*\hd, 1*\vd) -- (11*\hd, 1*\vd) -- (11*\hd, 2*\vd) -- (5*\hd, 2*\vd) -- (5*\hd, 1*\vd);

\draw[thick] (5*\hd, 0*\vd) -- (6*\hd, 0*\vd);
\node at (10.5*\hd, 1.45*\vd) {\tiny $x+\lambda_{r_i}$}; 
\draw[thick] (5*\hd, 1*\vd) -- (5*\hd, 0*\vd) -- (8*\hd, 0*\vd) -- (8*\hd, 1*\vd);
\node at (6.5*\hd, 0.45*\vd) {$\cdots$};

\node at (4.5*\hd, 1.45*\vd) {\tiny $r_i$th};

\node at (4.5*\hd, 0.45*\vd) {\tiny $(r_i + 1)$st};

\draw[thick] (5*\hd, 1*\vd) -- (6*\hd, 1*\vd) -- (6*\hd, 2*\vd);
\node at (5.5*\hd, 1.45*\vd) {\tiny $x + 1$};

\draw[thick] (6*\hd, 1*\vd) -- (7*\hd, 1*\vd) -- (7*\hd, 2*\vd);
\node at (6.5*\hd, 1.45*\vd) {\scriptsize $x +2$};

\draw[thick] (7*\hd, 1*\vd) -- (9*\hd, 1*\vd) -- (9*\hd, 2*\vd);
\node at (8*\hd, 1.45*\vd) {$\cdots$};

\draw[thick] (9*\hd, 1*\vd) -- (10*\hd, 1*\vd) -- (10*\hd, 2*\vd);
\node at (9.5*\hd, 1.45*\vd) 
{\tiny $x + \lambda_{r_i} - 1$};
\end{tikzpicture}
$}
\]
On the other hand, if $j = i$, then 
the  boxes occupied by $x + 1$, $x+2$, $\ldots$, $x+{\lambda_{r_i}}$ in $\ttT_{B_i}$ are as follows:
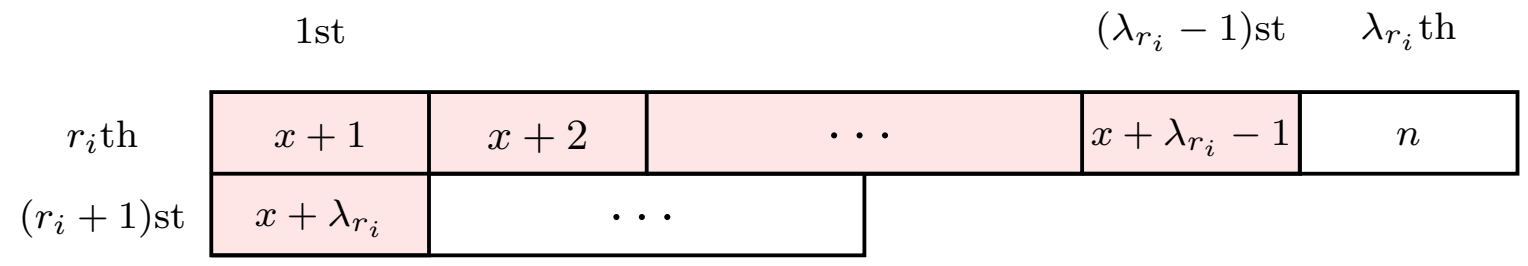
\[
\scalebox{0.9}{$
\def \hd {1.6}
\def \vd {0.6}
\begin{tikzpicture}
\node at (5.5*\hd, 2.75*\vd) {\tiny $1$st};
\node at (9.5*\hd, 2.75*\vd) {\tiny $(\lambda_{r_i} - 1)$st};
\node at (10.5*\hd, 2.75*\vd) {\tiny $\lambda_{r_i}$th};

\draw[thick, fill = red!10] (6*\hd, 1*\vd) -- (10*\hd, 1*\vd) -- (10*\hd, 2*\vd) -- (5*\hd, 2*\vd) -- (5*\hd, 1*\vd) -- (5*\hd, 0*\vd) -- (6*\hd, 0*\vd);
\draw[thick] (10*\hd, 2*\vd) -- (11*\hd, 2*\vd) --(11*\hd, 1*\vd) --  (10*\hd, 1*\vd);

\draw[thick] (5*\hd, 0*\vd) -- (6*\hd, 0*\vd) -- (6*\hd, 1*\vd);
\node at (5.5*\hd, 0.45*\vd) {\tiny $x+\lambda_{r_i}$};
\node at (10.5*\hd, 1.45*\vd) {\tiny $n$}; 
\draw[thick] (6*\hd, 0*\vd) -- (8*\hd, 0*\vd) -- (8*\hd, 1*\vd);
\node at (7*\hd, 0.45*\vd) {$\cdots$};

\node at (4.5*\hd, 1.45*\vd) {\tiny $r_i$th};

\node at (4.5*\hd, 0.45*\vd) {\tiny $(r_i + 1)$st};

\draw[thick] (5*\hd, 1*\vd) -- (6*\hd, 1*\vd) -- (6*\hd, 2*\vd);
\node at (5.5*\hd, 1.45*\vd) {\tiny $x + 1$};

\draw[thick] (6*\hd, 1*\vd) -- (7*\hd, 1*\vd) -- (7*\hd, 2*\vd);
\node at (6.5*\hd, 1.45*\vd) {\scriptsize $x +2$};

\draw[thick] (7*\hd, 1*\vd) -- (9*\hd, 1*\vd) -- (9*\hd, 2*\vd);
\node at (8*\hd, 1.45*\vd) {$\cdots$};

\draw[thick] (9*\hd, 1*\vd) -- (10*\hd, 1*\vd) -- (10*\hd, 2*\vd);
\node at (9.5*\hd, 1.45*\vd) 
{\tiny $x + \lambda_{r_i} - 1$};
\end{tikzpicture}
$}
\]
Since  $x+\lambda_{r_i} < n$, we have   $\opi_{x+1}\opi_{x+2} \cdots 
\opi_{x+\lambda_{r_i} - 1} \in H_{\C;0}(\SG_{n-1})$.
By  \cref{eq: left action for X}, 
\begin{align*}
\opi_{x+1} \opi_{x+2} \cdots \opi_{x + \lambda_{r_i} - 1} \cdot \ttg_{B_i}
& = \sum_{1 \leq j \leq i} \opi_{x+1} \opi_{x+2} \cdots \opi_{x + \lambda_{r_i} - 1} \cdot \ttT_{B_j} \\
&= \sum_{1 \leq j \leq i-1} (-1)^{\lambda_{r_i} - 1} \ttT_{B_j} + \opi_{x+1} \cdot( s_{x+2} \cdots s_{x+\lambda_{r_i} - 1} \cdot \ttT_{B_i})\\
&= (-1)^{\lambda_{r_i} - 1}\ttg_{B_{i-1}}.
\end{align*} 
This implies that  $\ttg_{B_{i-1}}  \in M_i$, and hence  $M_i \supseteq M_{i-1} = X_{B_{i-1}} \supseteq \SYT(\lambda)^{B_{i-1}}_{\boxn}$.
And, since  $\ttT_{B_i} = \ttg_{B_i} - \ttg_{B_{i-1}} \in M_i$, it follows
from 
 \cref{eq: left action for X} that  $M_i \supseteq \SYT(\lambda)^{B_{i}}_{\boxn} \setminus \SYT(\lambda)^{B_{i-1}}_{\boxn}$.
Therefore,  $M_i \supseteq \SYT(\lambda)^{B_i}_{\boxn}$.
\end{proof}

\begin{remark}\label{rem: SB is not known to be cyclic}
Unlike $\tilS_B$ and $X_B$, $S_B$ is not known to be cyclic.
The method used for $\tilS_B$ cannot be applied to $S_B$.
One of the obstacles is that we need to divide by a polynomial
that is not invertible in $\C[v^{1/2}]$.
For example, 
\cref{eq: h and z scalar diff} involves such a polynomial.
Indeed, if $\lambda = (2,1)$ and $i = 2$, then
$\sfh_{\mathbb{T}} = v\bfz_{\lambda^{\langle i \rangle}}$.
Since $v$ is not invertible in $\C[v^{1/2}]$, this equality does not imply that $\bfz_{\lambda^{\langle i \rangle}} \in S_B$.
The same issue arises in the proof of \cref{lem: basic step for inclusion}(1).
\end{remark}

\section{Hom-spaces between dual plactic class modules}\label{subsec: Hom set Gamma Q and Gamma Qp}
As established in \cref{thm: obtain CGammaQ}, for each $Q \in \SYT_n$, we have
\[
\calF_{\iota}(\bfM^{(D)}_Q) \cong \frakM_Q(\C(v^{1/2});v)
\qquad\text{and}\qquad
\calF_{\theta_0}(\bfM^{(D)}_Q) \cong \scrD_Q.
\]
Based on these isomorphisms, we investigate the Hom-spaces between dual plactic class modules and the classification of these modules. 
\medskip

\noindent
{\bf Convention.} 
Let $M$ and $N$ be $H_{R;a}(\SG_n)$-modules, where $R$ is a PID and $a \in R$. When the pair $(R,a)$ is clear from the context, we write $\Hom(M,N)$ for $\Hom_{H_{R;a}(\SG_n)}(M,N)$. We also abbreviate
$\frakM_Q(\C(v^{1/2});v)$
to $\frakM_Q$, $\frakJ_Q(\C(v^{1/2});v)$  to $\frakJ_Q$, and $\hfrakJ_Q(\C(v^{1/2});v)$  to $\hfrakJ_Q$.
\medskip

We begin by recalling the following standard facts on Kazhdan--Lusztig cell modules of $\genericH$ (see, for example, \cite{79KL, 96GL}).
For $Q_1, Q_2 \in \SYT_n$, we have
\begin{align}\label{eq: KL module classification}
\frakM_{Q_1}
\cong
\frakM_{Q_2} \quad\text{if and only if}\quad \sh(Q_1) = \sh(Q_2).
\end{align}
When $\sh(Q_1)=\sh(Q_2)$, an explicit $\genericH$-module isomorphism is given as follows.
Let $\sfu_{Q_1,Q_2}:\Gamma_{Q_1} \ra \Gamma_{Q_2}$ 
be the bijection such that
\[
\Ptab(\sfu_{Q_1, Q_2}(w))=\Ptab(w)
\quad
\text{for all } w\in\Gamma_{Q_1}.
\]
Define a $\C(v^{1/2})$-linear map 
$f_{Q_1,Q_2}:
\frakM_{Q_1}
\ra
\frakM_{Q_2}$  by
\[
f_{Q_1,Q_2}
\bigl(
C_w+\hfrakJ_{Q_1}
\bigr)
=
C_{\sfu_{Q_1,Q_2}(w)}
+\hfrakJ_{Q_2} \quad\text{for all $w \in \Gamma_{Q_1}$}.
\]
It follows from \cite[Theorem 4.10]{96GL} that $f_{Q_1, Q_2}$ is  an $\genericH$-module isomorphism.
Since the Kazhdan--Lusztig cell modules of $\genericH$ are irreducible, combining \cref{eq: KL module classification} with Schur's lemma yields
\begin{align}\label{eq: generic cell module hom sets} 
\Hom\left(
\frakM_{Q_1}, \, \frakM_{Q_2} 
\right) = \begin{cases} \C(v^{1/2})\,\{f_{Q_1,Q_2}\} & \text{if } \sh(Q_1)=\sh(Q_2),\\ 
0 & \text{otherwise.} 
\end{cases} 
\end{align}
If $\sh(Q_1) = \sh(Q_2)$, then for any $w \in \Gamma_{Q_1}$,
\begin{align*}
f_{Q_1,Q_2}(D_w+\hfrakJ_{Q_1})
&=
f_{Q_1,Q_2}(v^{\ell(w)/2}C_w+\hfrakJ_{Q_1})\\
&=
v^{\ell(w)/2}C_{\sfu_{Q_1,Q_2}(w)}+\hfrakJ_{Q_2}\\
&=
v^{(\ell(w)-\ell(\sfu_{Q_1,Q_2}(w)))/2}
D_{\sfu_{Q_1,Q_2}(w)}
+\hfrakJ_{Q_2}.
\end{align*}

We define a  $\C[v^{1/2}]$-linear map $f^{\rm poly}_{Q_1,Q_2}:\bfM^{(D)}_{Q_1} \ra \bfM^{(D)}_{Q_2}$ by
\[
f^{\rm poly}_{Q_1,Q_2}
\bigl(
D_w+\hbfJ^{(D)}_{Q_1}
\bigr)
=
v^{(\ell(w)-\ell(\sfu_{Q_1,Q_2}(w))-m_{Q_1,Q_2})/2}
D_{\sfu_{Q_1,Q_2}(w)}
+\hbfJ^{(D)}_{Q_2} \quad\text{for $w \in \Gamma_{Q_1}$},
\]
where
\[
m_{Q_1,Q_2}
:=
\min\{
\ell(w)-\ell(\sfu_{Q_1, Q_2}(w))
\mid
w\in\Gamma_{Q_1}
\}.
\]
As a $\C(v^{1/2})$-linear map, $\id \otimes f^{\rm poly}_{Q_1,Q_2}$ is a scalar multiple of $f_{Q_1, Q_2}$ if we regard $\C(v^{1/2}) \otimes_{\C[v^{1/2}]} \bfM^{(D)}_{Q_1}$ and $\C(v^{1/2}) \otimes_{\C[v^{1/2}]} \bfM^{(D)}_{Q_2}$ as
$\frakM_{Q_1}$ and $\frakM_{Q_2}$, respectively.
Therefore,  $f^{\rm poly}_{Q_1, Q_2}$ is an  $\polyH$-module homomorphism.

\begin{lemma}\label{thm: Hom set of MdQ and MdQp}
Let $Q_1, Q_2 \in \SYT_n$.
Then,
\[
\Hom(\bfM^{(D)}_{Q_1}, \bfM^{(D)}_{Q_2}) = \begin{cases}
\C[v^{1/2}] \{f^{\rm poly}_{Q_1, Q_2} \}    & \text{if $\sh(Q_1) = \sh(Q_2)$,}\\
0    & \text{otherwise.}
\end{cases}
\]
\end{lemma}
\begin{proof}
Let $g: \bfM^{(D)}_{Q_1} \ra \bfM^{(D)}_{Q_2}$ be an $\polyH$-module homomorphism.
For $w_1 \in \Gamma_{Q_1}$, write
\[
g(D_{w_1} + \hbfJ^{(D)}_{Q_1}) = \sum_{w_2 \in \Gamma_{Q_2}} a_{w_1, w_2} D_{w_2} + \hbfJ^{(D)}_{Q_2} \quad(a_{w_1, w_2} \in \C[v^{1/2}]).
\]
Then $\calF_\iota(g): \calF_\iota(\bfM^{(D)}_{Q_1})  \ra \calF_\iota(\bfM^{(D)}_{Q_2})$ is an $\genericH$-module homomorphism such that
\[
\calF_\iota(g) ( 1 \otimes C_{w_1}) = \sum_{w_2 \in \Gamma_{Q_2}} v^{(\ell(w_2) - \ell(w_1))/2} a_{w_1, w_2} \otimes C_{w_2} \quad \text{for all $w_1 \in \Gamma_{Q_1}$.}
\]

{\bf Case 1:} $\sh(Q_1) = \sh(Q_2)$.
For simplicity,
set $\sfu := \sfu_{Q_1, Q_2}$ and $m := m_{Q_1, Q_2}$. 
By \cref{thm: obtain CGammaQ}(1) and \cref{eq: generic cell module hom sets}, 
we have the following.
\begin{enumerate}[label = $\bullet$]
    \item $a_{w_1, w_2} = 0$ for all $w_1 \in \Gamma_{Q_1}$ and $w_2 \in \Gamma_{Q_2}$ with $w_2 \neq \sfu(w_1)$.
    \item There exists  $c \in \C(v^{1/2})$ such that
    \[
    c = v^{(\ell(\sfu(w_1)) - \ell(w_1))/2}a_{w_1, \sfu(w_1)} \quad\text{for all $w_1 \in \Gamma_{Q_1}$.}
    \]
\end{enumerate}
It follows that for all $w_1 \in \Gamma_{Q_1}$,
\begin{align*}
g(D_{w_1} + \hbfJ^{(D)}_{Q_1}) 
&= a_{w_1, \sfu(w_1)} D_{\sfu(w_1)} + \hbfJ^{(D)}_{Q_2} =  cv^{(\ell(w_1) - \ell(\sfu(w_1)))/2}  D_{\sfu(w_1)} + \hbfJ^{(D)}_{Q_2}.
\end{align*}
Since $g$ is an $\polyH$-module homomorphism, each coefficient $c v^{(\ell(w_1) - \ell(\sfu(w_1)))/2}$ lies in $\C[v^{1/2}]$. 
By the definition of $m$, so does $c v^{m/2}$.
Therefore, for all $w_1 \in \Gamma_{Q_1}$, 
\begin{align*}
g(D_{w_1} + \hbfJ^{(D)}_{Q_1}) 
&=  c v^{(\ell(w_1) - \ell(\sfu(w_1)))/2}  D_{\sfu(w_1)} + \hbfJ^{(D)}_{Q_2} \\
&= cv^{m/2} v^{(\ell(w_1) - \ell(\sfu(w)) -m)/2}  D_{\sfu(w_1)} + \hbfJ^{(D)}_{Q_2}\\
&= cv^{m/2} f^{\rm poly}_{Q_1, Q_2}(D_w + \hbfJ^{(D)}_{Q_1}).
\end{align*}
This shows that $g = cv^{m/2} f^{\rm poly}_{Q_1, Q_2}$. 

{\bf Case 2: }$\sh(Q_1) \neq \sh(Q_2)$.
By \cref{thm: obtain CGammaQ}(1) and \cref{eq: generic cell module hom sets},
$a_{w_1, w_2} = 0$ for all $w_1 \in \Gamma_{Q_1}$ and $w_2 \in \Gamma_{Q_2}$.
This shows that $g = 0$.
\end{proof}

We now investigate $\Hom(\scrD_{Q_1}, \scrD_{Q_2})$ when $\sh(Q_1) = \sh(Q_2)$.
Let $g_{Q_1, Q_2}: \scrD_{Q_1} \ra \scrD_{Q_2}$ be the $\C$-linear map defined by
\[
g_{Q_1, Q_2}(w)
=
\begin{cases}
\sfu_{Q_1,Q_2}(w)
& \text{if }
\ell(w)-\ell(\sfu_{Q_1,Q_2}(w)) = m_{Q_1,Q_2},\\
0
& \text{otherwise}
\end{cases}
\]
for $w \in \Gamma_{Q_1}$.
By definition, $g_{Q_1, Q_2}$ is nonzero.
We will show that $g_{Q_1,Q_2}$ is an $\compH{0}$-module homomorphism and characterize precisely when it is an isomorphism.
To this end, we associate with each subset $S\subseteq \SG_n$ a colored digraph as follows.
Its vertices are the permutations in $S$ and its $[n-1]$-colored arrows are given as follows. 
For $w,w'\in S$ and $i\in[n-1]$, we put $w\overset{i}{\longrightarrow}w'$ if and only if
\begin{align}\label{eq: colored digraph str def}
w=w' \text{ and } i\in\Des_L(w),
\quad\text{or}\quad  s_iw=w' \text{ and } i \not\in \Des_L(w).
\end{align}
In particular, for $Q\in\SYT_n$, the colored digraph associated with
$\Gamma_Q$ encodes the $\compH{0}$-action on $\scrD_Q$.
More precisely, for $w\in\Gamma_Q$ and
$i\in[n-1]$, the following hold:
\begin{enumerate}[label={\bf E\arabic*.}]
\item $w\overset{i}{\longrightarrow}w$ if and only if
$\opi_i\cdot w=-w$.
\item $w\overset{i}{\longrightarrow}s_i w$ if and only if
$\opi_i\cdot w=s_i w$.
\item No $i$-colored arrow starts at $w$ if and only if
$\opi_i\cdot w = 0$.
\end{enumerate}

\begin{lemma}\label{lem: new lem for g Q Q}
Let $Q_1, Q_2 \in \SYT_n$ with $\sh(Q_1) = \sh(Q_2)$.
\begin{enumerate}[label = {\rm (\arabic*)}]
\item $g_{Q_1, Q_2}$ is an $\compH{0}$-module homomorphism. 
\item $g_{Q_1, Q_2}$ is an $\compH{0}$-module isomorphism if and only if  $\sfu_{Q_1,Q_2}$ is an isomorphism of colored digraphs.
\end{enumerate}
\end{lemma}
\begin{proof}
Set
$\sfu:=\sfu_{Q_1,Q_2}$, $m:=m_{Q_1,Q_2}$, and $g := g_{Q_1, Q_2}$.

(1) 
By \cref{cor: algebra homomorphism vs functor sub of regular}(2) and
\cref{thm: obtain CGammaQ}(3), for each
$Q \in \SYT_n$, there is an $\compH{0}$-module isomorphism
\[
\uppsi_{Q}: \calF_{\theta_0}(\bfM^{(D)}_Q) \ra \scrD_{Q}, \quad 1 \otimes(D_w+\hbfJ^{(D)}_Q) \mapsto w\quad (w \in \Gamma_Q).
\]
For any $w \in \Gamma_{Q_1}$, we have
\begin{align*}
\uppsi_{Q_2}\circ \calF_{\theta_0}(f^{\rm poly}_{Q_1, Q_2}) \circ \uppsi_{Q_1}^{-1}(w) 
& = \uppsi_{Q_2}\circ \calF_{\theta_0}(f^{\rm poly}_{Q_1, Q_2})
\left(
1 \otimes \bigl(D_w + \hbfJ^{(D)}_{Q_1} \bigr) \right) \\
& = 
\begin{cases}
\uppsi_{Q_2}\left(
1 \otimes \bigl(D_{\sfu(w)} + \hbfJ^{(D)}_{Q_2}\bigr)
\right)
& \text{if $\ell(w) - \ell(\sfu(w)) = m$} ,
\\
0 & \text{otherwise}
\end{cases} \\
& =
\begin{cases}
\sfu(w)
& \text{if $\ell(w) - \ell(\sfu(w)) = m$} ,
\\
0 & \text{otherwise}
\end{cases} \\
& = g(w).
\end{align*}
Therefore, $g = \uppsi_{Q_2}
\circ \calF_{\theta_0}(f^{\rm poly}_{Q_1,Q_2})
\circ \uppsi_{Q_1}^{-1}$ is an $\compH{0}$-module homomorphism.

(2)
We first prove the ``only if'' part.
Suppose that $g$ is an $\compH{0}$-module isomorphism.
Then $g(w)=\sfu(w)$ for all $w\in\Gamma_{Q_1}$.
We show that  for all $w,w' \in \Gamma_{Q_1}$ and $i \in [n-1]$,
\begin{align}\label{eq: graph implication}
w \overset{i}{\longrightarrow} w' 
\quad \text{implies} \quad
\sfu(w) \overset{i}{\longrightarrow} \sfu(w').
\end{align}
Fix $w,w' \in \Gamma_{Q_1}$ and $i \in [n-1]$ such that $w \overset{i}{\longrightarrow} w'$.
By the definition of $[n-1]$-colored arrows in \cref{eq: colored digraph str def}, $w'$ is either $w$ or $s_i w$.
If $w' = w$, then $\opi_i \cdot w = -w$ by \textbf{E1}, and thus
\[
\opi_i \cdot \sfu(w) = \opi_i \cdot g(w) = g(\opi_i \cdot w) = g(-w) = -\sfu(w).
\]
Again by \textbf{E1}, we have $\sfu(w) \overset{i}{\longrightarrow} \sfu(w)$.
We next assume that $w' = s_i w$.
By \textbf{E2}, we have $\opi_i \cdot w = s_i w$, and thus 
\[
\opi_i \cdot \sfu(w) = \opi_i \cdot g(w) =  g(\opi_i \cdot w) = g(s_i w) = \sfu(s_i w).
\]
Since $\sfu(s_i w)$ is neither $-\sfu(w)$ nor $0$, it follows from \cref{eq: action on DQ} that $\sfu(s_i w) = \opi_i \cdot \sfu(w) = s_i \sfu(w)$.
Combining this with  \textbf{E2}, we have $\sfu(w) \overset{i}{\longrightarrow} \sfu(s_i w)$.
For the converse of \cref{eq: graph implication}, note that
$\sfu^{-1}=\sfu_{Q_2,Q_1}$ and $g^{-1}=g_{Q_2,Q_1}$.
Applying the preceding argument with $Q_1$ and $Q_2$ interchanged proves the converse. 
Therefore, $\sfu$ is an isomorphism of colored digraphs.

For the ``if''  part, assume that $\sfu$ is an isomorphism of colored digraphs.
By (1), it remains to show that $g$ is bijective.
By definition of $g$, this can be shown by proving that for any distinct $w, w' \in \Gamma_{Q_1}$, 
\begin{align*}
\ell(w) - \ell(\sfu(w)) = \ell(w') - \ell(\sfu(w')).
\end{align*}
Let $w$ and $w'$ be any distinct elements in $\Gamma_{Q_1}$.
Note that the dual Knuth equivalence can also be defined as the transitive closure of the dual Knuth relations of the first and second kinds (see \cite[Definition 3.6.8]{91Sag}).
It follows that there exist $\gamma_0, \gamma_1, \ldots, \gamma_k \in \Gamma_{Q_1}$ 
and $j_0, j_1, \ldots, j_{k-1} \in [n-1]$ such that
\[
\gamma_0=w, \quad
\gamma_k=w', \quad \text{and} \quad
\gamma_{i+1}=s_{j_i}\gamma_i
\quad \text{for all $0\leq i\leq k-1$.}
\]
Since $\sfu$ is an isomorphism of colored digraphs, for $0 \le i \leq k-1$,
\[
\gamma_i \overset{j_i}{\longrightarrow}\gamma_{i+1}
\quad \text{if and only if} \quad
\sfu(\gamma_i) \overset{j_i}{\longrightarrow} \sfu(\gamma_{i+1})
\]
and
\[
\gamma_{i+1} \overset{j_i}{\longrightarrow}\gamma_{i}\quad \text{if and only if} \quad
\sfu(\gamma_{i+1}) \overset{j_i}{\longrightarrow} \sfu(\gamma_{i}).
\]
It follows that 
\[
\ell(\gamma_i)-\ell(\sfu(\gamma_i)) = \ell(\gamma_{i+1})-\ell(\sfu(\gamma_{i+1})) 
\quad \text{for all $0\leq i\leq k-1$,}
\]
thus
$\ell(w) - \ell(\sfu(w)) = \ell(w') - \ell(\sfu(w'))$. 
\end{proof}

Before stating the next lemma,
we introduce two functors that will be used in its proof.
By \cite[Proposition~3.2]{05Fayers}, $\compH{0}$ has an automorphism
$\autophi$ and an anti-automorphism $\hautotheta$ defined by
\[
\autophi(\opi_i)=\opi_{n-i}
\quad\text{and}\quad
\hautotheta(\opi_i)=-1-\opi_i
\quad\text{for all $1\leq i\leq n-1$}.
\]
These maps induce the equivalence and duality 
\begin{align*}
\bfT_{\autophi}^{+}
& \colon
\compH{0}\text{-}\module\longrightarrow\compH{0}\text{-}\module\\
\bfT_{\hautotheta}^{-}
&\colon
\compH{0}\text{-}\module\longrightarrow\compH{0}\text{-}\module,
\end{align*}
called the $\autophi$-twist and the $\hautotheta$-twist, respectively.
For the precise definitions of them, see \cite[Subsection 3.4]{22JKLO}.
By \cite[Proposition~3.3]{05Fayers}, for an $\compH{0}$-module $M$, we have
\begin{align}\label{eq: same characteristic and twist}
\ch([\bfT^{+}_{\autophi}(M)]) = \uprho \circ \ch([M]) \quad\text{and}\quad \ch([\bfT^{-}_{\hautotheta}(M)]) = \uppsi \circ \ch([M]).
\end{align}
Here, $\uprho$ and $\uppsi$ are automorphisms of $\Qsym$ defined by
\[
\uprho(F_{\alpha}) = F_{\alpha^\rmr}
\quad\text{and}\quad
\uppsi(F_{\alpha}) = F_{\alpha^\rmc}
\]
for every composition $\alpha$.

\begin{lemma}\label{lem: basic lemma for Hom}
Let $\lambda \vdash n$, and let $M$ and $N$ be $\compH{0}$-modules such that $\ch([M]) = \ch([N])=s_\lambda$.
If at least one of $M$ and $N$ is isomorphic to either $X_\lambda$ or $X_{\lambda^\circ}$, then
\[
\dim \Hom(M,N)\leq 1.
\]
\end{lemma}
\begin{proof}
We first consider the case where $M$ is isomorphic to $X_{\lambda}$.
Let
\[
0 = N_0 \subsetneq N_1 \subsetneq \cdots \subsetneq N_l = N
\]
be a composition series of $N$.
Choose a basis $\{v_1,v_2,\ldots,v_l\}$ of $N$  so that $\{v_1,v_2,\ldots,v_j\}$ is a basis of $N_j$
for each $1 \leq j \leq l$.
For all $1 \leq i \leq n-1$ and $1 \leq j \leq l$, write 
\begin{align}\label{eq: action on the basis}
\opi_i \cdot v_j = \sum_{1 \leq k \leq j} c^{i}_{j, k} v_{k} \quad(c^{i}_{j, k} \in \C).
\end{align}
Note that $\bfDes(\ttT_{\lambda}) = \set(\lambda)^\rmc$ and $\bfDes(\ttT) \not\supseteq \set(\lambda)^\rmc$ for all $\ttT \in \SYT(\lambda) \setminus \{\ttT_{\lambda}\}$. Thus, we have \begin{align}\label{eq: descent property of normal SYT} \bfDes(\ttT) \not\supseteq \bfDes(\ttT_{\lambda}) \quad\text{for all $\ttT \in \SYT(\lambda) \setminus \{\ttT_{\lambda}\}$.} 
\end{align}
Since every $N_i/N_{i-1}$ is irreducible and 
\[
\sum_{1 \leq i \leq l} \ch([N_i/N_{i-1}]) = 
\ch([N])
=
s_{\lambda}
= 
\sum_{\ttT\in\SYT(\lambda)}
F_{\comp(\bfDes(\ttT)^\rmc)},
\]
we may enumerate the elements of $\SYT(\lambda)$ as
$\ttT_1,\ttT_2,\ldots,\ttT_l$ so that
\[
N_j/N_{j-1}
\cong
\bfF_{\comp(\bfDes(\ttT_j)^{\rmc})}
\qquad
\text{for each }1\leq j\leq l.
\]
It follows that for $1\leq i\leq n-1$ and  $1\leq j\leq l$,
\begin{align}\label{eq: nonzero leading condition}
c^i_{j,j}
=
\begin{cases}
-1 & \text{if } i\in\bfDes(\ttT_j),\\
0  & \text{if } i\notin\bfDes(\ttT_j).
\end{cases}
\end{align}
Denote by $r$ the index $1 \leq j \leq l$ such that
$\ttT_{j}=\ttT_\lambda$.

Let $f \in\Hom(X_\lambda,N)$ be nonzero.
By \cite[Lemma 3.11]{19Searles}, $\ttT_\lambda$ is a cyclic generator
of $X_\lambda$, and hence $f(\ttT_\lambda)\neq 0$.
Moreover,  it follows from
\cref{eq: left action for X} that
\begin{align}\label{eq: opi action on g(T lambda)}
\opi_i\cdot f(\ttT_\lambda)
=
f(\opi_i\cdot\ttT_\lambda)
=
-f(\ttT_\lambda)
\quad
\text{for all }i\in\bfDes(\ttT_{\lambda}).
\end{align}
Write
\[
f(\ttT_\lambda)=\sum_{j=1}^l a_jv_j
\quad (a_j\in\C),
\]
and let $m$ be the largest index such that $a_{m} \neq 0$.
We claim that $m=r$.
Suppose, to the contrary, that $m\neq r$.
Then $\ttT_{m}\neq\ttT_\lambda$, so
\cref{eq: descent property of normal SYT} yields an element $1 \leq i \leq n-1$ such that $i\in\bfDes(\ttT_\lambda)\setminus\bfDes(\ttT_{m})$.
Hence, by \cref{eq: nonzero leading condition}, 
 $c^i_{m,m}  = 0$.
By \cref{eq: action on the basis},
\[
\opi_i\cdot f(\ttT_\lambda)
=
\opi_i\cdot\left(\sum_{j=1}^{m}a_jv_j\right)
=
\sum_{j=1}^{m}\sum_{k=1}^j
a_j c^i_{j,k} v_k.
\]
The coefficient of $v_{m}$ in this expression is
$a_{m} \, c^i_{m,m} = 0$.
On the other hand, since
$i\in\bfDes(\ttT_\lambda)$,
\cref{eq: opi action on g(T lambda)} shows that this coefficient is
$-a_{m} \neq 0$, a contradiction.
Therefore, $m = r$ and thus
$f(\ttT_\lambda)
=
\sum_{j=1}^{r}a_jv_j$
with $a_r \neq 0$.

Choose any nonzero $g \in \Hom(X_{\lambda}, M)$. 
By the preceding argument, we may write
\[
g(\ttT_\lambda)
=
\sum_{j=1}^{r}b_jv_j
\quad  (b_j\in\C),
\]
where $b_{r}\neq0$.
Set $h := f- (a_{r} / b_{r})g$.
Then
\[
h(\ttT_\lambda)
=
\sum_{j=1}^{r-1}
\left(
a_j-\frac{a_{r}}{b_{r}}b_j
\right)v_j.
\]
If $h$ is nonzero, by the preceding argument, the
coefficient of $v_{r}$ in $h(\ttT_\lambda)$ is nonzero.
But this is a contradiction.
It follows that $f$ is a scalar multiple of $g$.
Thus,
\begin{align}\label{eq: dim X lam N leq 1}
\dim\Hom(X_\lambda,N)\leq1.
\end{align}

For the remaining cases, we will use  the $\autophi$-twist and $\hautotheta$-twist.
By
\cite[Proposition A.2(3)]{23KLO},
\begin{align*}
\bfT_{\autophi}^+(X_{\lambda})\cong X_{\lambda^\circ}
\quad\text{and}\quad
\bfT_{\hautotheta}^-(X_{\lambda})\cong X_{\lambda'}.
\end{align*}
Using these isomorphisms, we have
\begin{equation}\label{eq: three Hom space iso}
\begin{aligned}
& \Hom(X_{\lambda^\circ}, N)
\cong \Hom\Bigl(\bfT^{+}_{\autophi}( X_{\lambda^\circ}), \bfT^{+}_{\autophi}(N)\Bigr)
\cong \Hom(X_{\lambda}, \bfT^{+}_{\autophi}(N)), \\
& \Hom(M,  X_{\lambda})  
\cong \Hom\Bigl(\bfT^{-}_{\hautotheta}(X_{\lambda}), \bfT^{-}_{\hautotheta}(M)\Bigr) 
\cong \Hom(X_{\lambda'}, \bfT^{-}_{\hautotheta}(M)), \quad \text{and} \\
& \Hom(M,  X_{\lambda^\circ}) 
\cong
\Hom\Bigl(
\bfT_{\hautotheta}^-\circ \bfT^{+}_{\autophi}(X_{\lambda^\circ}), \bfT_{\hautotheta}^- \circ \bfT^{+}_{\autophi}(M) \Bigr)
\cong \Hom(X_{\lambda'}, \bfT_{\hautotheta}^- \circ \bfT^{+}_{\autophi}(M)).
\end{aligned}
\end{equation}
By \cref{eq: same characteristic and twist}, we have 
\begin{equation*}
\begin{array}{c}
\ch([\bfT^{+}_{\autophi}(N)]) = \ch([X_{\lambda}]) = s_{\lambda}, \quad
\ch([\bfT^{-}_{\hautotheta}(M)]) = \ch([X_{\lambda'}]) = s_{\lambda'}, \\[1ex]\quad\text{and}\quad\ch([\bfT^{-}_{\hautotheta} \circ \bfT^{+}_{\autophi}(M)]) = \ch([X_{\lambda'}]) = s_{\lambda'}.
\end{array}
\end{equation*}
Therefore, the desired results for the remaining cases can be obtained by combining \cref{eq: dim X lam N leq 1} and \cref{eq: three Hom space iso}.
\end{proof}

Given two subsets $S_1$ and $ S_2$ of  $\SG_n$, we write $S_1 \Deq S_2$ if they are isomorphic as colored digraphs.
We are now ready to state the main result of this section.
\begin{theorem}\label{thm: new thm for D Q iso to X}
Let $\lambda\vdash n$ and $Q, Q_1, Q_2 \in \SYT(\lambda')$. 
\begin{enumerate}[label={\rm (\arabic*)}]
\item
$\Gamma_Q$ is a left weak Bruhat interval if and only if $\scrD_Q$ is isomorphic to one of $X_\lambda$ and $X_{\lambda^\circ}$.

\item If at least one of $\Gamma_{Q_1}$ and $\Gamma_{Q_2}$ is a left weak Bruhat interval, then
\[
\Hom(\scrD_{Q_1},\scrD_{Q_2}) = \C\{g_{Q_1,Q_2}\}.
\]

\item Suppose that $ \Hom(\scrD_{Q_1}, \scrD_{Q_2}) = \C\{g_{Q_1, Q_2}\}$. 
Then, 
\begin{align*}
\scrD_{Q_1} \cong \scrD_{Q_2} \quad\text{if and only if}\quad \Gamma_{Q_1} \Deq \Gamma_{Q_2}.
\end{align*}
\end{enumerate}
\end{theorem}
\begin{proof}
(1) We first prove the ``if'' part.
We begin the case where $\scrD_{Q} \cong X_{\lambda}$.
Choose any $P \in \sfRSP_n$ such that $M_{P} \cong X_{\lambda}$.
It follows from \cite[Theorem 6.4]{23KLO} that $\Sigma_L(P)$ is a dual plactic class.  
Fix any $\sigma \in \Sigma_L(P)$.
Then $\Gamma_{\Qtab(\sigma)} = \Sigma_L(P)$, and thus
\[
\scrD_{\Qtab(\sigma)} = M_{P}  \cong X_{\lambda} \cong \scrD_Q.
\]
Since $\ch([\scrD_{\Qtab(\sigma)}]) = s_{\lambda}$, it follows from \cref{eq: ch of D Q} that  $\sh(\Qtab(\sigma)) = \lambda'$.
Combining \cref{lem: new lem for g Q Q}(1) with  \cref{lem: basic lemma for Hom}, we have that  $g_{Q, \Qtab(\sigma)}$ is an $\compH{0}$-module isomorphism.
By \cref{lem: new lem for g Q Q}(2), $\sfu_{Q, \Qtab(\sigma)}$ is an isomorphism of colored digraphs.
And,  $\Gamma_{\Qtab(\sigma)} = \Sigma_L(P)$ is a left weak Bruhat interval since $P \in \sfRSP_n$ (see \cite[Theorem 6.8]{91BW}).
Therefore, $\Gamma_{Q}$ is a left weak Bruhat interval.
The case where  $\scrD_{Q} \cong X_{\lambda^\circ}$  can be proved in the same manner.

Next, we prove the ``only if'' part.
By \cite[Theorem 6.8]{91BW}, there exists a poset $\mathbb{P} \in \sfP_n$ such that  $\Gamma_Q = \Sigma_L(\mathbb{P})$.
Thus, $\scrD_{Q}$ and $M_{\mathbb{P}}$ are identical as  $\compH{0}$-modules (see  the definitions of $\scrD_Q$ and $M_{\mathbb{P}}$ in \cref{subsec: two kinds of notable 0 Hecke modules}).
Since $\Gamma_Q$ is  closed under $\dKeq$, it follows from \cite[Theorem 6.4]{23KLO} that $\mathbb{P} \in \sfRSP_n$.
Recall that $\{X_{\nu/\kappa} \mid \nu/\kappa \text{ is a skew partition of $n$}\}$ forms a complete set of representatives for the isomorphism classes of $M_{P}$ ($P \in \sfRSP_n$).
Letting $X_{\nu/\kappa} \cong M_{\mathbb{P}}$, we have
\[
s_{\nu/\kappa} 
= \ch([X_{\nu/\kappa}])
= \ch([M_{\mathbb{P}}])
= s_{\lambda}. 
\]
By \cite[Theorem 2.2]{05W}, $\nu/\kappa$ is $\lambda$ or $\lambda^\circ$.
Therefore, $M_{\mathbb{P}}$ is isomorphic to $X_{\lambda}$ or $X_{\lambda^\circ}$.

(2) The assertion can be derived by putting   \cref{lem: new lem for g Q Q}(1),  \cref{lem: basic lemma for Hom}, and (1) together.

(3) We prove the ``if'' part. 
Suppose that $\Gamma_{Q_1} \Deq \Gamma_{Q_2}$.
Let $f: \Gamma_{Q_1} \ra \Gamma_{Q_2}$ be an isomorphism of colored digraphs, and let $\tf$ be the $\C$-linear extension of $f$.
Since $f$ is an isomorphism of colored digraphs, it follows from \textbf{E1}, \textbf{E2}, and \textbf{E3} that $\tf$ is an $\compH{0}$-module isomorphism. 

We next prove the ``only if'' part.
Suppose that $\scrD_{Q_1} \cong \scrD_{Q_2}$.
Since $\Hom(\scrD_{Q_1}, \scrD_{Q_2}) = \C\{g_{Q_1, Q_2}\}$, $g_{Q_1, Q_2}$ is an $\compH{0}$-module isomorphism.
Therefore, by \cref{lem: new lem for g Q Q}(2), $\sfu_{Q_1, Q_2}$ is an isomorphism of colored digraphs.
\end{proof}

\begin{figure}[h]
\centering
\begin{minipage}{0.4\textwidth}
\centering
\def \vp {1.3}
\def \hp {1.5}
\begin{tikzpicture}
%%%node
\node[] at (-1*\hp, 0*\vp) {$41325$};
\node[] at (1*\hp, 0*\vp) {$32415$};

\node[] at (-1*\hp, -1*\vp) {$51324$};
\node[] at (1*\hp, -1*\vp) {$42315$};

\node[] at (-1*\hp, -2*\vp) {$51423$};
\node[] at (1*\hp, -2*\vp) {$52314$};

\node[] at (0.4*\hp, -2.8*\vp){$\Gamma_{Q_1}$};

%%%%%%%%%%%%%%%%%%loop%%%%%%%%%%%%%%%%%%%%%%%%%%%%
\node at (1*\hp + 0.25*\hp, 0*\vp) {} edge [out=40,in=320, loop] ();
\node[right] at (1*\hp + 0.6*\hp, 0*\vp) {\footnotesize $1, 2$};

\node at (-1*\hp + 0.25*\hp, 0*\vp) {} edge [out=40,in=320, loop] ();
\node[right] at (-1*\hp + 0.6*\hp, 0*\vp) {\footnotesize $2, 3$};

\node at (-1*\hp + 0.25*\hp, -1*\vp) {} edge [out=40,in=320, loop] ();
\node[right] at (-1*\hp + 0.6*\hp, -1*\vp) {\footnotesize $2, 4$};
\node at (1*\hp + 0.25*\hp, -1*\vp) {} edge [out=40,in=320, loop] ();
\node[right] at (1*\hp + 0.6*\hp, -1*\vp) {\footnotesize $1, 3$};

\node at (-1*\hp + 0.25*\hp, -2*\vp) {} edge [out=40,in=320, loop] ();
\node[right] at (-1*\hp + 0.6*\hp, -2*\vp) {\footnotesize $3, 4$};
\node at (1*\hp + 0.25*\hp, -2*\vp) {} edge [out=40,in=320, loop] ();
\node[right] at (1*\hp + 0.6*\hp, -2*\vp) {\footnotesize $1, 4$};

%%%%%%%%%%%%%%%%%%%%%%EDGE%%%%%%%%%%%%%%%%%%%%%%%%
\draw[->] (-1*\hp, -0.25*\vp) -- (-1*\hp, -0.75*\vp);
\node[left] at (-1*\hp, -0.5*\vp) {\footnotesize $4$};
\draw[->] (-0.85*\hp, -0.25*\vp) -- (0.85*\hp, -0.75*\vp);
\node[above] at (0*\hp, -0.55*\vp) {\footnotesize $1$};

\draw[->] (-1*\hp, -1.25*\vp) -- (-1*\hp, -1.75*\vp);
\node[left] at (-1*\hp, -1.5*\vp) {\footnotesize $3$};
\draw[->] (-0.85*\hp, -1.25*\vp) -- (0.85*\hp, -1.75*\vp);
\node[above] at (0*\hp, -1.55*\vp) {\footnotesize $1$};

%%%%
\draw[->] (1*\hp, -0.25*\vp) -- (1*\hp, -0.75*\vp);
\node[right] at (1*\hp, -0.5*\vp) {\footnotesize $3$};
\draw[->] (1*\hp, -1.25*\vp) -- (1*\hp, -1.75*\vp);
\node[right] at (1*\hp, -1.5*\vp) {\footnotesize $4$};
\end{tikzpicture}
\end{minipage}
\hfill
\begin{minipage}{0.56\textwidth}
\centering
\def \vp {1.3}
\def \hp {1.5}
\begin{tikzpicture}
%%%node
\node[] at (-2*\hp, 0*\vp) {$51243$};
\node[] at (0*\hp, 0*\vp) {$41352$};
\node[] at (2*\hp, 0*\vp) {$32451$};

\node[] at (-1*\hp, -1*\vp) {$51342$};
\node[] at (1*\hp, -1*\vp) {$42351$};

\node[] at (0*\hp, -2*\vp) {$52341$};

\node[] at (0.3*\hp, -2.8*\vp){$\Gamma_{Q_2}$};

%%%%%%%%%%%%%%%%%%loop%%%%%%%%%%%%%%%%%%%%%%%%%%%%
\node at (-2*\hp + 0.25*\hp, 0*\vp) {}
    edge [out=40,in=320, loop] ();
\node[right] at (-2*\hp + 0.6*\hp, 0*\vp)
    {\footnotesize $3, 4$};

\node at (0*\hp + 0.25*\hp, 0*\vp) {}
    edge [out=40,in=320, loop] ();
\node[right] at (0*\hp + 0.6*\hp, 0*\vp)
    {\footnotesize $2, 3$};

\node at (2*\hp + 0.25*\hp, 0*\vp) {}
    edge [out=40,in=320, loop] ();
\node[right] at (2*\hp + 0.6*\hp, 0*\vp)
    {\footnotesize $1, 2$};

\node at (-1*\hp + 0.25*\hp, -1*\vp) {}
    edge [out=40,in=320, loop] ();
\node[right] at (-1*\hp + 0.6*\hp, -1*\vp)
    {\footnotesize $2, 4$};

\node at (1*\hp + 0.25*\hp, -1*\vp) {}
    edge [out=40,in=320, loop] ();
\node[right] at (1*\hp + 0.6*\hp, -1*\vp)
    {\footnotesize $1, 3$};

\node at (0*\hp + 0.25*\hp, -2*\vp) {}
    edge [out=40,in=320, loop] ();
\node[right] at (0*\hp + 0.6*\hp, -2*\vp)
    {\footnotesize $1, 4$};

%%%%%%%%%%%%%%%%%%%%%%EDGE%%%%%%%%%%%%%%%%%%%%%%%%
\draw[->]
    (-2*\hp, -0.25*\vp) -- (-1.15*\hp, -0.75*\vp);
\node[right] at (-1.55*\hp, -0.45*\vp)
    {\footnotesize $2$};

\draw[->]
    (-0.15*\hp, -0.25*\vp) -- (-0.85*\hp, -0.75*\vp);
\node[left] at (-0.5*\hp, -0.45*\vp)
    {\footnotesize $4$};

\draw[->]
    (0.15*\hp, -0.25*\vp) -- (0.85*\hp, -0.75*\vp);
\node[right] at (0.5*\hp, -0.45*\vp)
    {\footnotesize $1$};

\draw[->]
    (2*\hp, -0.25*\vp) -- (1.15*\hp, -0.75*\vp);
\node[left] at (1.55*\hp, -0.45*\vp)
    {\footnotesize $3$};

\draw[->]
    (-1*\hp, -1.25*\vp) -- (-0.15*\hp, -1.75*\vp);
\node[right] at (-0.55*\hp, -1.45*\vp)
    {\footnotesize $1$};

\draw[->]
    (1*\hp, -1.25*\vp) -- (0.15*\hp, -1.75*\vp);
\node[left] at (0.55*\hp, -1.45*\vp)
    {\footnotesize $4$};
\end{tikzpicture}
\end{minipage}
\hfill
\begin{minipage}{0.98\textwidth}
\centering
\def \vp {1.3}
\def \hp {1.5}
\begin{tikzpicture}
%%%node
\node[] at (0*\hp, 0*\vp) {$14325$};

\node[] at (-1*\hp, -1*\vp) {$15324$};
\node[] at (1*\hp, -1*\vp) {$24315$};

\node[] at (-2*\hp, -2*\vp) {$15423$};
\node[] at (0*\hp, -2*\vp) {$25314$};
\node[] at (2*\hp, -2*\vp) {$34215$};

\node[] at (0.3*\hp, -2.8*\vp){$\Gamma_{Q_3}$};
%%%%%%%%%%%%%%%%%%loop%%%%%%%%%%%%%%%%%%%%%%%%%%%%
\node at (0*\hp + 0.25*\hp, 0*\vp) {}
    edge [out=40,in=320, loop] ();
\node[right] at (0*\hp + 0.6*\hp, 0*\vp)
    {\footnotesize $2, 3$};

\node at (-1*\hp + 0.25*\hp, -1*\vp) {}
    edge [out=40,in=320, loop] ();
\node[right] at (-1*\hp + 0.6*\hp, -1*\vp)
    {\footnotesize $2, 4$};

\node at (1*\hp + 0.25*\hp, -1*\vp) {}
    edge [out=40,in=320, loop] ();
\node[right] at (1*\hp + 0.6*\hp, -1*\vp)
    {\footnotesize $1, 3$};

\node at (-2*\hp + 0.25*\hp, -2*\vp) {}
    edge [out=40,in=320, loop] ();
\node[right] at (-2*\hp + 0.6*\hp, -2*\vp)
    {\footnotesize $3, 4$};

\node at (0*\hp + 0.25*\hp, -2*\vp) {}
    edge [out=40,in=320, loop] ();
\node[right] at (0*\hp + 0.6*\hp, -2*\vp)
    {\footnotesize $1, 4$};

\node at (2*\hp + 0.25*\hp, -2*\vp) {}
    edge [out=40,in=320, loop] ();
\node[right] at (2*\hp + 0.6*\hp, -2*\vp)
    {\footnotesize $1, 2$};

%%%%%%%%%%%%%%%%%%%%%%EDGE%%%%%%%%%%%%%%%%%%%%%%%%
\draw[->]
    (-0.15*\hp, -0.25*\vp) -- (-0.85*\hp, -0.75*\vp);
\node[left] at (-0.5*\hp, -0.45*\vp)
    {\footnotesize $4$};

\draw[->]
    (0.15*\hp, -0.25*\vp) -- (0.85*\hp, -0.75*\vp);
\node[right] at (0.5*\hp, -0.45*\vp)
    {\footnotesize $1$};

\draw[->]
    (-1.15*\hp, -1.25*\vp) -- (-1.85*\hp, -1.75*\vp);
\node[left] at (-1.5*\hp, -1.45*\vp)
    {\footnotesize $3$};

\draw[->]
    (-0.85*\hp, -1.25*\vp) -- (-0.15*\hp, -1.75*\vp);
\node[right] at (-0.5*\hp, -1.45*\vp)
    {\footnotesize $1$};

\draw[->]
    (0.85*\hp, -1.25*\vp) -- (0.15*\hp, -1.75*\vp);
\node[left] at (0.5*\hp, -1.45*\vp)
    {\footnotesize $4$};

\draw[->]
    (1.15*\hp, -1.25*\vp) -- (1.85*\hp, -1.75*\vp);
\node[right] at (1.5*\hp, -1.45*\vp)
    {\footnotesize $2$};
\end{tikzpicture}
\end{minipage}
\caption{The colored digraphs associated with
$\Gamma_{Q_i}$ ($i = 1,2,3$)}
\label{fig: Hom dimension two example}
\end{figure}
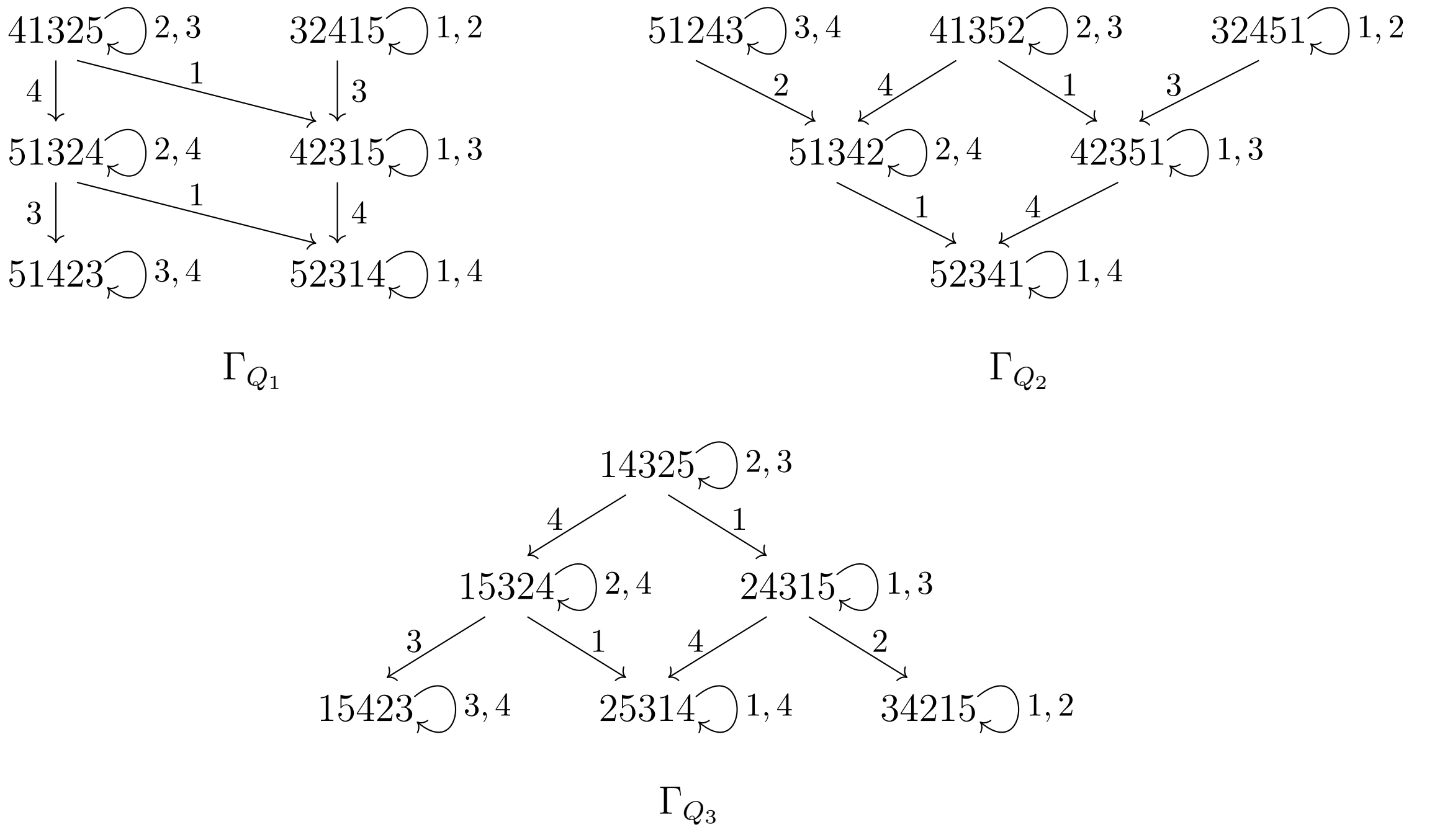

\begin{remark}\label{rem: new remark for theorem on Hom}
Let $Q_1, Q_2 \in \SYT_n$ with $\sh(Q_1) = \sh(Q_2)$.

\noindent (1) 
Consider the case where exactly one of $\Gamma_{Q_1}$ and $\Gamma_{Q_2}$ is a left weak Bruhat interval.
Since $\Gamma_{Q_1} \not\Deq\Gamma_{Q_2}$,  it follows from \cref{thm: new thm for D Q iso to X}(2) and (3) that  $\scrD_{Q_1} \not\cong \scrD_{Q_2}$.

\noindent (2) Consider the case where  neither $\Gamma_{Q_1}$ nor $\Gamma_{Q_2}$ is a left weak Bruhat interval.
In this case, the dimension of $\Hom(\scrD_{Q_1},  \scrD_{Q_2})$ can be strictly greater than $1$.
For example, let 
\[
Q_1 = 
\begin{array}{l}
\begin{ytableau}
    1 & 3 & 5\\
    2\\
    4
\end{ytableau}
\end{array},\quad Q_2 = \begin{array}{l}\begin{ytableau}
    1 & 3 & 4\\
    2\\
    5
\end{ytableau}
\end{array}, \quad\text{and}\quad
Q_3 = 
\begin{array}{l}
\begin{ytableau}
    1 & 2 & 5\\
    3\\
    4
\end{ytableau}
\end{array}.
\]
As seen in \cref{fig: Hom dimension two example}, none of $\Gamma_{Q_i}$'s ($i = 1,2,3$) are  left weak Bruhat intervals.
After simple computation, we see that
\[ 
\dim \Hom(\scrD_{Q_1}, \scrD_{Q_2}) = 1 \quad\text{and}\quad 
\dim \Hom(\scrD_{Q_2}, \scrD_{Q_3}) = 2.
\]
\end{remark}

\appendix

\section{Lemmas for Section 4}\label[appendix]{sec: Appendix 1}

In this appendix, we introduce \cref{lem: triangularity of elt in B} and \cref{lem: Garnir rel}. 
These results are deferred to the appendix to avoid interrupting the flow of the main argument in \cref{sec: 0-Hecke modules from specializations}.

We begin with a lemma needed for the proof of \cref{lem: triangularity of elt in B}.

\begin{lemma}\label{lem: mult of T_sigma and T_rho}
For $\sigma, \rho \in \SG_n$, let
\[
T_\sigma T_\rho = \sum_{\tau \in \SG_n} a^{\tau}_{\sigma, \rho} T_{\tau} \quad(a^{\tau}_{\sigma, \rho} \in \Z[v^{1/2}]).
\]
Then we have the following.
\begin{enumerate}[label = {\rm (\arabic*)}]
\item 
If $a^{\tau}_{\sigma, \rho} \neq 0$, then $\ell(\tau) \geq \ell(\rho) - \ell(\sigma)$.

\item If $\sigma^{-1} \preceq_R \rho$ and $\ell(\tau) = \ell(\rho) - \ell(\sigma)$, then 
\[
a^\tau_{\sigma,\rho} \neq 0
\quad \text{if and only if} \quad
\tau = \sigma \rho.
\]
Moreover, $a^{\sigma \rho}_{\sigma,\rho} = v^{\ell(\sigma)}$.
\end{enumerate}
\end{lemma}

\begin{proof}
(1) Let $\sigma, \rho \in \SG_n$.
We prove the assertion by mathematical induction on $\ell(\sigma)$.
If $\ell(\sigma) = 0$, i.e., $\sigma = \id$, then $T_{\sigma}T_{\rho} = T_{\rho}$ and thus the assertion holds.

Now, suppose that $\ell(\sigma) > 0$ and that the assertion holds for $\sigma' \in \SG_n$ with $\ell(\sigma') < \ell(\sigma)$.
Let $\tau \in \SG_n$ with $a^{\tau}_{\sigma, \rho} \neq 0$.
We have two cases.
\smallskip

{\bf Case 1:} $\Des_R(\sigma) \cap \Des_L(\rho) = \emptyset$.
Let $i \in \Des_R(\sigma)$.
Then, we have $T_\sigma T_\rho = T_{\sigma s_i} T_{s_i \rho}$, and thus 
$a^\tau_{\sigma, \rho} = a^{\tau}_{\sigma s_i, s_i \rho}$.
Since  $a^{\tau}_{\sigma, \rho} \neq 0$,
by the induction hypothesis,
\begin{equation}\label{eq: case one tau has strict longer legnth cond}
\ell(\tau) \geq \ell(s_i \rho) - \ell(\sigma s_i) = \ell(\rho) - \ell(\sigma) + 2 > \ell(\rho) - \ell(\sigma).
\end{equation}  
Thus the assertion holds.
\smallskip

{\bf Case 2:} $\Des_R(\sigma) \cap \Des_L(\rho) \neq \emptyset$.
Let $i \in \Des_R(\sigma) \cap \Des_L(\rho)$.
Then, we have 
\[
T_\sigma T_\rho = 
T_{\sigma s_i} T_i T_i T_{s_i \rho} = T_{\sigma s_i} ((v-1) T_i + v) T_{s_i \rho} = (v-1)T_{\sigma s_i} T_\rho + v T_{\sigma s_i} T_{s_i \rho}.
\]
It follows that
\begin{equation}\label{eq: a sig tau sum of two}
a^{\tau}_{\sigma, \rho} = (v-1) a^{\tau}_{\sigma s_i, \rho} + v a^{\tau}_{\sigma s_i, s_i \rho}.
\end{equation}
By \cref{eq: a sig tau sum of two}, at least one of $a^{\tau}_{\sigma s_i , \rho}$ and $a^{\tau}_{\sigma s_i , s_i \rho}$ is nonzero.
By the induction hypothesis, if $a^{\tau}_{\sigma s_i, \rho} \neq 0$, then
\begin{align}\label{eq: first coeff is zero}
\ell(\tau) \geq  \ell(\rho) - \ell(\sigma s_i) = \ell(\rho) - \ell(\sigma) + 1,
\end{align}
and if $a^{\tau}_{\sigma s_i, s_i \rho} \neq 0$, then
\[
\ell(\tau) \geq \ell(s_i\rho) - \ell(\sigma s_i) = \ell(\rho)-\ell(\sigma).
\]
Thus the assertion holds.

(2) Let $\sigma, \rho \in \SG_n$ with $\sigma^{-1} \preceq_R \rho$.
We prove the assertion by mathematical induction on $\ell(\sigma)$.
If $\ell(\sigma) = 0$,  then $T_{\sigma}T_{\rho} = T_{\rho}$ and thus the assertion holds.

Now, suppose that $\ell(\sigma) > 0$ and that the assertion holds for $\sigma' \in \SG_n$ with $\ell(\sigma') < \ell(\sigma)$.
Let $\tau \in \SG_n$ with $\ell(\tau) = \ell(\rho) - \ell(\sigma)$.

Let $i \in \Des_R(\sigma)$.
Since $\sigma^{-1} \preceq_R \rho$, $\Des_R(\sigma) \subseteq \Des_L(\rho)$.
Following the proof of {\bf Case 2} of (1), we see that \cref{eq: a sig tau sum of two} and \cref{eq: first coeff is zero} hold.
Since $\ell(\tau) = \ell(\rho) - \ell(\sigma)$, \cref{eq: first coeff is zero} implies that $a^{\tau}_{\sigma s_i, \rho} = 0$.
Combining this with \cref{eq: a sig tau sum of two} yields 
\begin{equation}\label{eq: a is v a in the case where minimal length}
a^{\tau}
_{\sigma, \rho} = v a^{\tau}_{\sigma s_i, s_i \rho}.
\end{equation}
Since $\sigma^{-1} \preceq_R \rho$ and $i \in \Des_L(\sigma^{-1}) \subseteq \Des_L(\rho)$, we have  $s_i \sigma^{-1} \preceq_R s_i \rho$.
By the induction hypothesis,
\[
a^{\tau}_{\sigma s_i, s_i \rho} \neq 0 \quad\text{if and only if}\quad \tau = \sigma \rho,
\]
and $a^{\sigma\rho}_{\sigma s_i, s_i \rho} = v^{\ell(\sigma s_i)}$.
Now, the assertion  follows from \cref{eq: a is v a in the case where minimal length}.
\end{proof}

\begin{lemma}\label{lem: triangularity of elt in  B}
Let $\lambda/\mu$ be a skew partition of $n$.
For $\ttS \in \SYT(\lambda/\mu)$, 
let 
\[
\sfh_{\ttS} = \sum_{\tau \in \SG_n} c_\tau T_\tau \quad\text{$(c_\tau \in \C[v^{1/2}])$}.
\]
Then, we have the following.
\begin{enumerate}[label = {\rm (\arabic*)}]
    \item If $c_\tau \neq 0$, then $\ell(\tau) \geq \ell(\wcol(\ttS))$.
    \item If $\ell(\tau) = \ell(\wcol(\ttS))$, then 
    \[
c_\tau \neq 0    \quad\text{if and only if}\quad \tau = \wcol(\ttS).
    \]
Moreover,  $c_{\wcol(\ttS)} = v^{\ell(w_0(\alpha(\lambda/\mu))) + \ell(\wrow(\ttS))}$.
\end{enumerate}
\end{lemma}
\begin{proof}
Set $\alpha := \alpha(\lambda/\mu)$ and $\beta := \alpha(\lambda'/\mu')$.
We first show that for $\tau \in \SG_n$,
\begin{align}\label{eq: c tau coeff is sum of a tau sig rho}
c_\tau = \sum_{\substack{
x \in \SG_\alpha\\ y \in \SG_\beta}} \sgn(x) v^{\ell(w_0(\alpha)) - \ell(x)} a^\tau_{\wrow(\ttS), x \wcol(\ttT_{\lambda/\mu}) y}.
\end{align}
For the definition of
$a^\tau_{\sigma,\rho}$, see
\cref{lem: mult of T_sigma and T_rho}.
This can be accomplished by proving that
\begin{align}\label{eq: sfh ttS in another form}
\sfh_{\ttS}  =  \sum_{\substack{
x \in \SG_\alpha\\ y \in \SG_\beta}
}\sgn(x) v^{\ell(w_0(\alpha)) - \ell(x)} T_{\wrow(\ttS)} T_{x \wcol(\ttT_{\lambda/\mu}) y}.
\end{align}
Since
\[
\sfh_{\ttS} = T_{\wrow(\ttS)} \bfy_{\alpha}T_{\wcol(\ttT_{\lambda/\mu})}\bfx_{\beta} =  \sum_{\substack{
x \in \SG_\alpha\\ y \in \SG_\beta}
}\sgn(x) v^{\ell(w_0(\alpha)) - \ell(x)} T_{\wrow(\ttS)} T_x  T_{\wcol(\ttT_{\lambda/\mu})} T_y,
\]
it suffices to show that for all $x \in \SG_\alpha$ and $y \in \SG_\beta$, 
\begin{align}\label{eq: appendix length preserving cond}
\ell(x \wcol(\ttT_{\lambda/\mu})y) = \ell(x) + \ell(\wcol(\ttT_{\lambda/\mu})) + \ell(y).    
\end{align}
Let $x \in \SG_{\alpha}$ and $y \in \SG_{\beta}$.
Since every column of $x\cdot\ttT_{\lambda/\mu}$ is strictly increasing from top to bottom, we have
\begin{align}\label{eq: ell x dot tt y}
\ell(\wcol(x \cdot \ttT_{\lambda/\mu}) y) = 
\ell(\wcol(x \cdot \ttT_{\lambda/\mu})) + \ell(y).
\end{align}
Moreover, since $x$ permutes the entries within each row of $\ttT_{\lambda/\mu}$, we have
\begin{align}\label{eq: ell x dot tt}
\ell(x \wcol(\ttT_{\lambda/\mu})) = \ell(x) + \ell(\wcol(\ttT_{\lambda/\mu})).
\end{align}
Combining \cref{eq: ell x dot tt y} and \cref{eq: ell x dot tt} yields \cref{eq: appendix length preserving cond}.

Now, we prove the assertion (1). 
 Suppose that $c_\tau\neq0$. Then there exist $x\in\SG_\alpha$ and $y\in\SG_\beta$ such that 
$a^\tau_{\wrow(\ttS),\,x\wcol(\ttT_{\lambda/\mu})y} \neq 0$.
By 
\cref{lem: mult of T_sigma and T_rho}(1) and \cref{eq: appendix length preserving cond}, 
\begin{align}\label{eq: length inequality condition}
\ell(\tau) \geq \ell(x \wcol(\ttT_{\lambda/\mu}) y) - \ell(\wrow(\ttS)) = \ell(\wcol(\ttT_{\lambda/\mu})) - \ell(\wrow(\ttS)) + \ell(x) + \ell(y).
\end{align}
We claim that 
\begin{align}\label{eq: wcol ttS is equal to wrow ttS and wcol ttT}
\ell(\wcol(\ttT_{\lambda/\mu})) - \ell(\wrow(\ttS)) = 
\ell(\wcol(\ttS)).
\end{align}
We first observe that
\[
\wrow(\ttS) \preceq_L \wrow(\ttT'_{\lambda/\mu})\quad\text{and}\quad  \wrow(\ttT'_{\lambda/\mu}) = \wcol(\ttT_{\lambda/\mu})^{-1},
\]
which implies that
\begin{align}\label{eq: right weak Bruhat order rel}
\wrow(\ttS)^{-1} \preceq_R \wrow(\ttT'_{\lambda/\mu})^{-1} = \wcol(\ttT_{\lambda/\mu}).
\end{align}
We also observe that 
$\wcol(\ttS)
=
\wrow(\ttS)\wcol(\ttT_{\lambda/\mu})$.
Combining these two observations yields \cref{eq: wcol ttS is equal to wrow ttS and wcol ttT}.
Finally,  plugging \cref{eq: wcol ttS is equal to wrow ttS and wcol ttT} to \cref{eq: length inequality condition} 
completes the proof.

Next, we prove the assertion (2).
Suppose that $\ell(\tau) = \ell(\wcol(\ttS))$.
In the case where  $x\in\SG_\alpha$ and $y\in\SG_\beta$ satisfy $a^{\tau}_{\wrow(\ttS), x\wcol(\ttT_{\lambda/\mu})y} \neq 0$, 
\cref{eq: length inequality condition}  remains still valid.
Plugging 
\cref{eq: wcol ttS is equal to wrow ttS and wcol ttT} to \cref{eq: length inequality condition} gives $x = y = \id$.
Thus, by \cref{eq: c tau coeff is sum of a tau sig rho}, we have
\begin{align}\label{eq: c tau is a signle sum of a tau sig rho}
c_\tau = v^{\ell(w_0(\alpha))} a^{\tau}_{\wrow(\ttS), \wcol(\ttT_{\lambda/\mu})}.
\end{align}
Due to \cref{eq: right weak Bruhat order rel},
we can use \cref{lem: mult of T_sigma and T_rho}(2) to obtain
\[
a^{\tau}_{\wrow(\ttS), \wcol(\ttT_{\lambda/\mu})} \neq 0\quad\text{if and only if}\quad \tau = \wrow(\ttS)\wcol(\ttT_{\lambda/\mu}) = \wcol(\ttS)
\]
and
\[ 
a^{\wcol(\ttS)}_{\wrow(\ttS), \wcol(\ttT_{\lambda/\mu})} = v^{\ell(\wrow(\ttS))}.
\]
Combining these with \cref{eq: c tau is a signle sum of a tau sig rho} completes the proof.
\end{proof}

\begin{lemma}\label{lem: Garnir rel}
Let $\lambda/\mu$ be a skew partition of  $n$, and let $(i, j) \in \tyd(\lambda/\mu)$ with $(i, j+1) \in \tyd(\lambda/\mu)$.
Then,
\begin{equation*}
G^{\lambda/\mu}_{i, j} \bfz_{\lambda/\mu}(\C[v^{1/2}];v) = 0.
\end{equation*}
\end{lemma}
\begin{proof}
Set $\alpha := \alpha(\lambda/\mu)$ and $\beta := \alpha(\lambda'/\mu')$.
Let $\nu$ be the composition of $n$ with length $\ell(\lambda) + 1$ defined by 
\[
\nu_t = \begin{cases}
  \alpha_t  & \text{if $1 \leq t < i$,}\\
  j - \mu_i - 1 & \text{if  $t = i$,}\\
  \lambda_i - \mu_{i+1} + 1 & \text{if  $t = i+1$,}\\
  \lambda_{i+1} - j & \text{if  $t = i+2$,}\\
  \alpha_{t - 1} & \text{if $i+2 < t \leq \ell(\lambda) + 1$.}
\end{cases}
\]
That is, $\nu$ is obtained from $\alpha$ by replacing $\alpha_i$ and $\alpha_{i+1}$ with three parts.
The part $\nu_{i+1}$ counts the boxes of $\tyd(\lambda/\mu)$ in the $i$th row weakly to the right of the $j$th column together with those in the $(i+1)$st row weakly to the left of the $j$th column.
The parts $\nu_i$ and $\nu_{i+2}$ count the remaining boxes in the $i$th and $(i+1)$st rows, respectively.
Let $\calD_L^\alpha$ be the set of minimal coset representatives of $\SG_n/\SG_\alpha$ and $\calD_R^\nu$ the set of minimal coset representatives of $\SG_\nu\backslash\SG_n$. 
One can easily see that
\begin{align}\label{eq: disjoint union of left and right cosets}
\left( \calD_L^\alpha \cap \SG_\nu \right) \SG_\alpha
= \{\sigma \rho \mid \sigma \in \SG_\nu,\, \rho \in \SG_\alpha \}
= \SG_\nu \left( \calD_R^\nu \cap \SG_\alpha \right)
\end{align}
Considering the expression for $\SG^{\lambda/\mu}_{i,j}$ in \cref{eq: Garnir permutations alternating def}, one can also see that $\SG^{\lambda/\mu}_{i,j} = \calD_L^\alpha \cap \SG_\nu$. 
Combining this with \cref{eq: disjoint union of left and right cosets} yields
\[
G^{\lambda/\mu}_{i, j} \bfy_\alpha 
= \bfy_{\SG^{\lambda/\mu}_{i,j}} \, \bfy_\alpha 
= \bfy_\nu \, \bfy_{\calD_R^\nu \cap \SG_\alpha}.
\]
Thus, for the assertion, it is enough to prove that  for every $\sigma \in \calD_R^\nu \cap \SG_\alpha$,
\begin{align}\label{eq: want to show for A3}
\bfy_\nu T_\sigma T_{\wcol(\ttT_{\lambda/\mu})} \bfx_\beta = 0.
\end{align}

In what follows, we will verify \cref{eq: want to show for A3}.
Let $\sigma \in \calD_R^\nu \cap \SG_\alpha$ and
set $\ttU_{\sigma} := \sigma \cdot \ttT_{\lambda/\mu}$.
Since $\wcol(\ttT_{\lambda/\mu})$ is the minimal  coset representative of the right coset $\SG_{\alpha} \wcol(\ttT_{\lambda/\mu})$, we have
\[
\ell(\sigma \wcol(\ttT_{\lambda/\mu})) =  \ell(\sigma) + \ell(\wcol(\ttT_{\lambda/\mu})).
\]
Combining this identity with  $\wcol(\ttU_\sigma) = \sigma \wcol(\ttT_{\lambda/\mu})$, we obtain
\begin{align}\label{eq: first step for A3}
T_{\sigma} T_{\wcol(\ttT_{\lambda/\mu})} = T_{\wcol(\ttU_\sigma)}.
\end{align}
Set $x := \sum_{1 \leq t \leq i-1} \alpha_t$ and $I := [x+\nu_i +1,x+\nu_i+\nu_{i+1}]$.
Since $\sigma \in \SG_{\alpha}$ and the elements of $I$ lie in the $i$th or the $(i+1)$st row of $\sfT_{\lambda/\mu}$, these elements lie in the same two rows of $\ttU_{\sigma}$.
Moreover, the entries in the $i$th and $(i+1)$st rows of $\ttU_\sigma$ occupy only $\lambda_i-\mu_{i+1}$ distinct columns and $|I|=\nu_{i+1}=\lambda_i-\mu_{i+1}+1$.
By the pigeonhole principle, two distinct elements in $I$ occur in the same column of $\ttU_\sigma$.
Hence, there exists $k \in [\mu_i+1, \lambda_{i+1}]$ such that
\[
\ttU_\sigma(i,k),\ \ttU_\sigma(i+1,k)\in I.
\]
Let $m \in [n]$ such that $\wcol(\ttU_\sigma)(m) = \ttU_\sigma(i, k)$.
Then $s_m \in \SG_\beta$, and hence 
\[
\bfx_{\beta} = (1+T_m) \sum_{\substack{\rho \in \SG_{\beta}\\m \notin \Des_L(\rho)}} T_{\rho}.
\] 
By virtue of this identity together with  \cref{eq: first step for A3}, \cref{eq: want to show for A3} can be derived if we show
\begin{align}\label{eq: bfy T (1+T) zero}
\bfy_\nu T_{\wcol(\ttU_\sigma)}(1+T_m)=0.
\end{align}
Set $a := \ttU_{\sigma}(i, k)$ and  $b := \ttU_{\sigma}(i+1, k)$.
Then $(a \  b) \wcol(\ttU_{\sigma}) = \wcol(\ttU_\sigma) s_m$.
Moreover, since $\sigma \in \SG_{\alpha}$ and each element in the $i$th row of $T_{\lambda/\mu}$ is less than every element in the $(i+1)$st row, we have
\begin{align*}
a = \ttU_{\sigma}(i, k) < \ttU_{\sigma}(i+1, k) = b.   
\end{align*}
Thus,
$\bfy_\nu T_{\wcol(\ttU_\sigma)} T_m = \bfy_{\nu} T_{\wcol(\ttU_{\sigma}) s_m} = \bfy_{\nu} T_{(a \; b)\wcol(\ttU_{\sigma})}$.
Since $(a \  b) \in \SG_{\nu}$, it follows from \cref{eq: Ti and bfy} that $\bfy_{\nu} T_{(a \;  b)\wcol(\ttU_{\sigma})} = -\bfy_{\nu} T_{\wcol(\ttU_{\sigma})}$.
This proves \cref{eq: bfy T (1+T) zero} and thus completes the proof.
\end{proof}

\paragraph{{\bf Funding statement.}}
The first author was supported by 
Basic Science Research Program through the National Research Foundation of Korea(NRF) funded by the Ministry of Education (No. RS-2023-00240377) and
NRF grant funded by the Korean Government(MSIT) (No. RS-2026-25595299).
The second author was supported by NRF grant funded by the Korean Government (No. NRF-2020R1F1A1A01071055), Basic Science Research Program through NRF funded by the Ministry of Education (No. RS-2023-00271282), and NRF grant funded by the Korea government(MSIT) (No. RS-2024-00342349).
The third author was supported by NRF grant funded by the Korean Government (No. NRF-2020R1F1A1A01071055) and by NRF grant funded by the Korea government(MSIT) (No. RS-2024-00342349).

\bibliographystyle{abbrv}
\bibliography{references}

@article {23KLO,
    AUTHOR = {Kim, Young-Hun and Lee, So-Yeon and Oh, Young-Tak},
     TITLE = {Regular {S}chur labeled skew shape posets and their 0-{H}ecke
              modules},
   JOURNAL = {Forum Math. Sigma},
  FJOURNAL = {Forum of Mathematics. Sigma},
    VOLUME = {12},
      YEAR = {2024},
     PAGES = {Paper No. e110, 40}
}

@article {05W,
    AUTHOR = {van Willigenburg, Stephanie},
     TITLE = {Equality of {S}chur and skew {S}chur functions},
   JOURNAL = {Ann. Comb.},
  FJOURNAL = {Annals of Combinatorics},
    VOLUME = {9},
      YEAR = {2005},
    NUMBER = {3},
     PAGES = {355--362},
      ISSN = {0218-0006,0219-3094},
   MRCLASS = {05E05 (05E10)},
  MRNUMBER = {2176598},
MRREVIEWER = {Grant\ Walker},
       DOI = {10.1007/s00026-005-0263-5},
       URL = {https://doi.org/10.1007/s00026-005-0263-5},
}

@article {86DJ,
    AUTHOR = {Dipper, Richard and James, Gordon},
     TITLE = {Representations of {H}ecke algebras of general linear groups},
   JOURNAL = {Proc. London Math. Soc. (3)},
  FJOURNAL = {Proceedings of the London Mathematical Society. Third Series},
    VOLUME = {52},
      YEAR = {1986},
    NUMBER = {1},
     PAGES = {20--52},
      ISSN = {0024-6115,1460-244X},
   MRCLASS = {20G05},
  MRNUMBER = {812444},
MRREVIEWER = {G.\ I.\ Lehrer},
       DOI = {10.1112/plms/s3-52.1.20},
       URL = {https://doi.org/10.1112/plms/s3-52.1.20},
}

@book {99Mat,
    AUTHOR = {Mathas, Andrew},
     TITLE = {Iwahori-{H}ecke algebras and {S}chur algebras of the symmetric
              group},
    SERIES = {University Lecture Series},
    VOLUME = {15},
 PUBLISHER = {American Mathematical Society, Providence, RI},
      YEAR = {1999},
     PAGES = {xiv+188},
      ISBN = {0-8218-1926-7},
   MRCLASS = {20C08 (20C20 20G99)},
  MRNUMBER = {1711316},
MRREVIEWER = {Meinolf\ Geck},
       DOI = {10.1090/ulect/015},
       URL = {https://doi.org/10.1090/ulect/015},
}

@article {79Norton,
    AUTHOR = {P. Norton},
     TITLE = {{$0$}-{H}ecke algebras},
   JOURNAL = {J. Austral. Math. Soc. Ser. A},
    VOLUME = {27},
      YEAR = {1979},
    NUMBER = {3},
     PAGES = {337--357}
}

@article {96DKLT,
    AUTHOR = {G. Duchamp and D. Krob and B. Leclerc and
              J.-Y. Thibon},
     TITLE = {Fonctions quasi-sym\'{e}triques, fonctions sym\'{e}triques non
              commutatives et alg\`ebres de {H}ecke \`a {$q=0$}},
   JOURNAL = {C. R. Acad. Sci. Paris S\'{e}r. I Math.},
    VOLUME = {322},
      YEAR = {1996},
    NUMBER = {2},
     PAGES = {107--112}
}

@article{19Searles,
    AUTHOR = {Searles, Dominic},
     TITLE = {Indecomposable {$0$}-{H}ecke modules for extended {S}chur
              functions},
   JOURNAL = {Proc. Amer. Math. Soc.},
  FJOURNAL = {Proceedings of the American Mathematical Society},
    VOLUME = {148},
      YEAR = {2020},
    NUMBER = {5},
     PAGES = {1933--1943}
}

@article {05Fayers,
    AUTHOR = {M. Fayers},
     TITLE = {{$0$}-{H}ecke algebras of finite {C}oxeter groups},
   JOURNAL = {J. Pure Appl. Algebra},
    VOLUME = {199},
      YEAR = {2005},
    NUMBER = {1-3},
     PAGES = {27--41}
}

@article {16Huang,
    AUTHOR = {J. Huang},
     TITLE = {A tableau approach to the representation theory of {$0$}-{H}ecke
              algebras},
   JOURNAL = {Ann. Comb.},
    VOLUME = {20},
      YEAR = {2016},
    NUMBER = {4},
     PAGES = {831--868}
}

@article {02DHT,
    AUTHOR = {G. Duchamp and F. Hivert and J.-Y. Thibon},
     TITLE = {Noncommutative symmetric functions. {VI}. {F}ree
              quasi-symmetric functions and related algebras},
   JOURNAL = {Internat. J. Algebra Comput.},
    VOLUME = {12},
      YEAR = {2002},
    NUMBER = {5},
     PAGES = {671--717},
}

@book {97Ful,
    AUTHOR = {Fulton, William},
     TITLE = {Young tableaux},
    SERIES = {London Mathematical Society Student Texts},
    VOLUME = {35},
      NOTE = {With applications to representation theory and geometry},
 PUBLISHER = {Cambridge University Press, Cambridge},
      YEAR = {1997},
     PAGES = {x+260},
}

@book{91Sag,
  AUTHOR = {Sagan, Bruce E.},
  TITLE = {The symmetric group - representations, combinatorial algorithms, and symmetric functions.},
  PAGES = {I-XV, 1-197},
  PUBLISHER = {Wadsworth},
  SERIES = {Wadsworth \& Brooks / Cole mathematics series},
  YEAR = 1991
}

@article {91BW,
    AUTHOR = {Bj\"{o}rner, Anders and Wachs, Michelle L.},
     TITLE = {Permutation statistics and linear extensions of posets},
   JOURNAL = {J. Combin. Theory Ser. A},
    VOLUME = {58},
      YEAR = {1991},
    NUMBER = {1},
     PAGES = {85--114}
}

@book {05BB,
    AUTHOR = {A. Bj\"{o}rner and F. Brenti},
     TITLE = {Combinatorics of {C}oxeter groups},
    SERIES = {Graduate Texts in Mathematics},
    VOLUME = {231},
 PUBLISHER = {Springer, New York},
      YEAR = {2005}
}

@book {72Stanley,
    AUTHOR = {Stanley, Richard},
     TITLE = {Ordered structures and partitions},
      NOTE = {Memoirs of the American Mathematical Society, No. 119},
 PUBLISHER = {American Mathematical Society, Providence, R.I.},
      YEAR = {1972},
     PAGES = {iii+104}
}

@book {99Stanley,
    AUTHOR = {Stanley, Richard},
     TITLE = {Enumerative combinatorics. {V}ol. 2},
    SERIES = {Cambridge Studies in Advanced Mathematics},
    VOLUME = {62},
 PUBLISHER = {Cambridge University Press, Cambridge},
      YEAR = {1999},
     PAGES = {xii+581}
}

@article {09BL,
    AUTHOR = {Bergeron, Nantel and Li, Huilan},
     TITLE = {Algebraic structures on {G}rothendieck groups of a tower of
              algebras},
   JOURNAL = {J. Algebra},
  FJOURNAL = {Journal of Algebra},
    VOLUME = {321},
      YEAR = {2009},
    NUMBER = {8},
     PAGES = {2068--2084},
      ISSN = {0021-8693,1090-266X},
   MRCLASS = {16T30 (05E05 05E15)},
  MRNUMBER = {2501510},
MRREVIEWER = {Mark\ J.\ Wildon},
       DOI = {10.1016/j.jalgebra.2008.12.005},
       URL = {https://doi.org/10.1016/j.jalgebra.2008.12.005},
}

@article {96GL,
    AUTHOR = {Graham, J. J. and Lehrer, G. I.},
     TITLE = {Cellular algebras},
   JOURNAL = {Invent. Math.},
  FJOURNAL = {Inventiones Mathematicae},
    VOLUME = {123},
      YEAR = {1996},
    NUMBER = {1},
     PAGES = {1--34},
      ISSN = {0020-9910,1432-1297},
   MRCLASS = {20C99},
  MRNUMBER = {1376244},
       DOI = {10.1007/BF01232365},
       URL = {https://doi.org/10.1007/BF01232365},
}

@article{03Geordie,
  title={Mind your {$P$} and {$Q$}-symbols: {W}hy the {K}azhdan-{L}usztig basis of the {H}ecke algebra of {T}ype {$A$} is cellular},
  author={Williamson, Geordie},
  year={2003}
}

@article {89GW,
    AUTHOR = {Garsia, A. M. and Wachs, M. L.},
     TITLE = {Combinatorial aspects of skew representations of the symmetric
              group},
   JOURNAL = {J. Combin. Theory Ser. A},
  FJOURNAL = {Journal of Combinatorial Theory. Series A},
    VOLUME = {50},
      YEAR = {1989},
    NUMBER = {1},
     PAGES = {47--81},
      ISSN = {0097-3165,1096-0899},
   MRCLASS = {20C30 (05A17)},
  MRNUMBER = {978066},
MRREVIEWER = {Michael\ H.\ Peel},
       DOI = {10.1016/0097-3165(89)90005-8},
       URL = {https://doi.org/10.1016/0097-3165(89)90005-8},
}

@article {81Lusztig,
    AUTHOR = {Lusztig, George},
     TITLE = {On a theorem of {B}enson and {C}urtis},
   JOURNAL = {J. Algebra},
  FJOURNAL = {Journal of Algebra},
    VOLUME = {71},
      YEAR = {1981},
    NUMBER = {2},
     PAGES = {490--498},
      ISSN = {0021-8693},
   MRCLASS = {20G05 (14L30)},
  MRNUMBER = {630610},
MRREVIEWER = {S.\ I.\ Gel\cprime fand},
       DOI = {10.1016/0021-8693(81)90188-5},
       URL = {https://doi.org/10.1016/0021-8693(81)90188-5},
}

@article {86Gyoja,
    AUTHOR = {Gyoja, Akihiko},
     TITLE = {A {$q$}-analogue of {Y}oung symmetrizer},
   JOURNAL = {Osaka J. Math.},
  FJOURNAL = {Osaka Journal of Mathematics},
    VOLUME = {23},
      YEAR = {1986},
    NUMBER = {4},
     PAGES = {841--852},
      ISSN = {0030-6126},
   MRCLASS = {20C15 (20C30)},
  MRNUMBER = {873212},
MRREVIEWER = {David\ B.\ Surowski},
       URL = {http://projecteuclid.org/euclid.ojm/1200779724},
}

@article {79PJ,
    AUTHOR = {Peel, M. H. and James, G. D.},
     TITLE = {Specht series for skew representations of symmetric groups},
   JOURNAL = {J. Algebra},
  FJOURNAL = {Journal of Algebra},
    VOLUME = {56},
      YEAR = {1979},
    NUMBER = {2},
     PAGES = {343--364},
      ISSN = {0021-8693},
   MRCLASS = {20C30},
  MRNUMBER = {528580},
MRREVIEWER = {J\"{u}rgen Tappe},
       DOI = {10.1016/0021-8693(79)90342-9},
       URL = {https://doi-org.libproxy.sogang.ac.kr/10.1016/0021-8693(79)90342-9},
}

@article {11DY,
    AUTHOR = {Deng, BangMing and Yang, GuiYu},
     TITLE = {Representation type of 0-{H}ecke algebras},
   JOURNAL = {Sci. China Math.},
    VOLUME = {54},
      YEAR = {2011},
    NUMBER = {3},
     PAGES = {411--420}
}

@article {22JKLO,
    AUTHOR = {Jung, Woo-Seok and Kim, Young-Hun and Lee, So-Yeon and Oh,
              Young-Tak},
     TITLE = {Weak {B}ruhat interval modules of the 0-{H}ecke algebra},
   JOURNAL = {Math. Z.},
    VOLUME = {301},
      YEAR = {2022},
    NUMBER = {4},
     PAGES = {3755--3786}
}

@article {23CKO,
    AUTHOR = {Choi, Seung-Il and Kim, Young-Hun and Oh, Young-Tak},
     TITLE = {{P}oset modules of the 0-{H}ecke algebras and related quasisymmetric power sum expansions},
   JOURNAL = {European J. Combin.},
    VOLUME = {120},
      YEAR = {2024},
     PAGES = {Paper No. 103965, 34}
}

@article {25Searles,
    AUTHOR = {Searles, Dominic},
     TITLE = {Diagram supermodules for 0-{H}ecke-{C}lifford algebras},
   JOURNAL = {Math. Z.},
  FJOURNAL = {Mathematische Zeitschrift},
    VOLUME = {310},
      YEAR = {2025},
    NUMBER = {3},
     PAGES = {Paper No. 43, 36},
      ISSN = {0025-5874,1432-1823},
   MRCLASS = {05E10 (20C08)},
  MRNUMBER = {4900538},
       DOI = {10.1007/s00209-025-03750-6},
       URL = {https://doi.org/10.1007/s00209-025-03750-6},
}

@article {91Ram,
    AUTHOR = {Ram, Arun},
     TITLE = {A {F}robenius formula for the characters of the {H}ecke
              algebras},
   JOURNAL = {Invent. Math.},
  FJOURNAL = {Inventiones Mathematicae},
    VOLUME = {106},
      YEAR = {1991},
    NUMBER = {3},
     PAGES = {461--488},
      ISSN = {0020-9910,1432-1297},
   MRCLASS = {20C30 (05E05 16W30 17B37)},
  MRNUMBER = {1134480},
MRREVIEWER = {Jie\ Du},
       DOI = {10.1007/BF01243921},
       URL = {https://doi.org/10.1007/BF01243921},
}

@incollection {80KL,
    AUTHOR = {Kazhdan, David and Lusztig, George},
     TITLE = {Schubert varieties and {P}oincar\'e{} duality},
 BOOKTITLE = {Geometry of the {L}aplace operator ({P}roc. {S}ympos. {P}ure
              {M}ath., {U}niv. {H}awaii, {H}onolulu, {H}awaii, 1979)},
    SERIES = {Proc. Sympos. Pure Math.},
    VOLUME = {XXXVI},
     PAGES = {185--203},
 PUBLISHER = {Amer. Math. Soc., Providence, RI},
      YEAR = {1980},
      ISBN = {0-8218-1439-7},
   MRCLASS = {14M15 (55N35 57R19)},
  MRNUMBER = {573434},
}

@incollection {00Ariki,
    AUTHOR = {Ariki, Susumu},
     TITLE = {Robinson-{S}chensted correspondence and left cells},
 BOOKTITLE = {Combinatorial methods in representation theory ({K}yoto,
              1998)},
    SERIES = {Adv. Stud. Pure Math.},
    VOLUME = {28},
     PAGES = {1--20},
 PUBLISHER = {Kinokuniya, Tokyo},
      YEAR = {2000},
      ISBN = {4-314-10141-5},
   MRCLASS = {20C08 (05E15 20C30)},
  MRNUMBER = {1855588},
MRREVIEWER = {Cheng\ Dong\ Chen},
       DOI = {10.2969/aspm/02810001},
       URL = {https://doi.org/10.2969/aspm/02810001},
}

@article {79KL,
    AUTHOR = {Kazhdan, David and Lusztig, George},
     TITLE = {Representations of {C}oxeter groups and {H}ecke algebras},
   JOURNAL = {Invent. Math.},
  FJOURNAL = {Inventiones Mathematicae},
    VOLUME = {53},
      YEAR = {1979},
    NUMBER = {2},
     PAGES = {165--184},
      ISSN = {0020-9910,1432-1297},
   MRCLASS = {20H15 (17B35 20G05 22E47)},
  MRNUMBER = {560412},
MRREVIEWER = {Vinay\ V.\ Deodhar},
       DOI = {10.1007/BF01390031},
       URL = {https://doi-org.vnlib.kias.re.kr/10.1007/BF01390031},
}

@book {00GP,
    AUTHOR = {Geck, Meinolf and Pfeiffer, G\"otz},
     TITLE = {Characters of finite {C}oxeter groups and {I}wahori-{H}ecke
              algebras},
    SERIES = {London Mathematical Society Monographs. New Series},
    VOLUME = {21},
 PUBLISHER = {The Clarendon Press, Oxford University Press, New York},
      YEAR = {2000},
     PAGES = {xvi+446},
      ISBN = {0-19-850250-8},
   MRCLASS = {20C15 (20C08 20F55)},
  MRNUMBER = {1778802},
MRREVIEWER = {Jian-yi\ Shi},
}

@book {89NT,
    AUTHOR = {Nagao, Hirosi and Tsushima, Yukio},
     TITLE = {Representations of finite groups},
      NOTE = {Translated from the Japanese},
 PUBLISHER = {Academic Press, Inc., Boston, MA},
      YEAR = {1989},
     PAGES = {xviii+424},
      ISBN = {0-12-513660-9},
   MRCLASS = {20Cxx (20-01 20C20)},
  MRNUMBER = {998775},
MRREVIEWER = {Roderick\ Gow},
}

@article {92KW,
    AUTHOR = {King, R. C. and Wybourne, B. G.},
     TITLE = {Representations and traces of the {H}ecke algebras {$H_n(q)$} of type {$A_{n-1}$}},
   JOURNAL = {J. Math. Phys.},
    VOLUME = {33},
      YEAR = {1992},
    NUMBER = {1},
     PAGES = {4--14}
}

@article {88GM,
    AUTHOR = {Garsia, A. M. and McLarnan, T. J.},
     TITLE = {Relations between {Y}oung's natural and the
              {K}azhdan-{L}usztig representations of {$S_n$}},
   JOURNAL = {Adv. in Math.},
  FJOURNAL = {Advances in Mathematics},
    VOLUME = {69},
      YEAR = {1988},
    NUMBER = {1},
     PAGES = {32--92},
      ISSN = {0001-8708},
   MRCLASS = {20C30},
  MRNUMBER = {937317},
MRREVIEWER = {Michael\ H.\ Peel},
       DOI = {10.1016/0001-8708(88)90060-6},
       URL = {https://doi-org.libproxy.sogang.ac.kr/10.1016/0001-8708(88)90060-6},
}

@article {81BK,
    AUTHOR = {Brylinski, J.-L. and Kashiwara, M.},
     TITLE = {Kazhdan-{L}usztig conjecture and holonomic systems},
   JOURNAL = {Invent. Math.},
  FJOURNAL = {Inventiones Mathematicae},
    VOLUME = {64},
      YEAR = {1981},
    NUMBER = {3},
     PAGES = {387--410},
      ISSN = {0020-9910,1432-1297},
   MRCLASS = {22E47 (17B10 20G05)},
  MRNUMBER = {632980},
MRREVIEWER = {H.\ H.\ Andersen},
       DOI = {10.1007/BF01389272},
       URL = {https://doi.org/10.1007/BF01389272},
}

@article {06Geck,
    AUTHOR = {Geck, Meinolf},
     TITLE = {Kazhdan-{L}usztig cells and the {M}urphy basis},
   JOURNAL = {Proc. London Math. Soc. (3)},
  FJOURNAL = {Proceedings of the London Mathematical Society. Third Series},
    VOLUME = {93},
      YEAR = {2006},
    NUMBER = {3},
     PAGES = {635--665},
      ISSN = {0024-6115,1460-244X},
   MRCLASS = {20C08 (20G05)},
  MRNUMBER = {2266962},
MRREVIEWER = {Andrew\ Mathas},
       DOI = {10.1017/S0024611506015930},
       URL = {https://doi.org/10.1017/S0024611506015930},
}

@article {14EBW,
    AUTHOR = {Elias, Ben and Williamson, Geordie},
     TITLE = {The {H}odge theory of {S}oergel bimodules},
   JOURNAL = {Ann. of Math. (2)},
  FJOURNAL = {Annals of Mathematics. Second Series},
    VOLUME = {180},
      YEAR = {2014},
    NUMBER = {3},
     PAGES = {1089--1136},
      ISSN = {0003-486X,1939-8980},
   MRCLASS = {20C08 (20F55)},
  MRNUMBER = {3245013},
MRREVIEWER = {Chi\ Kin\ Mak},
       DOI = {10.4007/annals.2014.180.3.6},
       URL = {https://doi.org/10.4007/annals.2014.180.3.6},
}

@article {96Ariki,
    AUTHOR = {Ariki, Susumu},
     TITLE = {On the decomposition numbers of the {H}ecke algebra of
              {$G(m,1,n)$}},
   JOURNAL = {J. Math. Kyoto Univ.},
  FJOURNAL = {Journal of Mathematics of Kyoto University},
    VOLUME = {36},
      YEAR = {1996},
    NUMBER = {4},
     PAGES = {789--808},
      ISSN = {0023-608X},
   MRCLASS = {20C20 (20G99)},
  MRNUMBER = {1443748},
MRREVIEWER = {Meinolf\ Geck},
       DOI = {10.1215/kjm/1250518452},
       URL = {https://doi-org.libproxy.swu.ac.kr/10.1215/kjm/1250518452},
}

@article {14WW,
    AUTHOR = {Wan, Jinkui and Wang, Weiqiang},
     TITLE = {Frobenius map for the centers of {H}ecke algebras},
   JOURNAL = {Trans. Amer. Math. Soc.},
  FJOURNAL = {Transactions of the American Mathematical Society},
    VOLUME = {367},
      YEAR = {2015},
    NUMBER = {8},
     PAGES = {5507--5520},
      ISSN = {0002-9947,1088-6850},
   MRCLASS = {20C08 (05E05)},
  MRNUMBER = {3347181},
MRREVIEWER = {Selvaraj\ Chelliah},
       DOI = {10.1090/S0002-9947-2014-06211-9},
       URL = {https://doi.org/10.1090/S0002-9947-2014-06211-9},
}

@article {21Korff,
    AUTHOR = {Korff, Christian},
     TITLE = {Cylindric {H}ecke characters and {G}romov-{W}itten invariants
              via the asymmetric six-vertex model},
   JOURNAL = {Comm. Math. Phys.},
  FJOURNAL = {Communications in Mathematical Physics},
    VOLUME = {381},
      YEAR = {2021},
    NUMBER = {2},
     PAGES = {591--640},
      ISSN = {0010-3616,1432-0916},
   MRCLASS = {14N35 (05E14 20C08)},
  MRNUMBER = {4207452},
       DOI = {10.1007/s00220-020-03906-x},
       URL = {https://doi.org/10.1007/s00220-020-03906-x},
}

@article {15SY,
    AUTHOR = {Savage, Alistair and Yacobi, Oded},
     TITLE = {Categorification and {H}eisenberg doubles arising from towers
              of algebras},
   JOURNAL = {J. Combin. Theory Ser. A},
  FJOURNAL = {Journal of Combinatorial Theory. Series A},
    VOLUME = {129},
      YEAR = {2015},
     PAGES = {19--56},
      ISSN = {0097-3165,1096-0899},
   MRCLASS = {16T05 (05E05 16W50 20C08)},
  MRNUMBER = {3275113},
MRREVIEWER = {Volodymyr\ Mazorchuk},
       DOI = {10.1016/j.jcta.2014.09.002},
       URL = {https://doi-org.vnlib.kias.re.kr/10.1016/j.jcta.2014.09.002},
}

@book {02Bourbaki,
    AUTHOR = {Bourbaki, Nicolas},
     TITLE = {Lie groups and {L}ie algebras. {C}hapters 4--6},
    SERIES = {Elements of Mathematics (Berlin)},
 PUBLISHER = {Springer-Verlag, Berlin},
      YEAR = {2002},
     PAGES = {xii+300},
      ISBN = {3-540-42650-7},
   MRCLASS = {17-01 (00A05 20E42 20F55 22-01)},
  MRNUMBER = {1890629},
       DOI = {10.1007/978-3-540-89394-3},
       URL = {https://doi.org/10.1007/978-3-540-89394-3},
}

@article {09BK,
    AUTHOR = {Brundan, Jonathan and Kleshchev, Alexander},
     TITLE = {Blocks of cyclotomic {H}ecke algebras and {K}hovanov-{L}auda
              algebras},
   JOURNAL = {Invent. Math.},
  FJOURNAL = {Inventiones Mathematicae},
    VOLUME = {178},
      YEAR = {2009},
    NUMBER = {3},
     PAGES = {451--484},
      ISSN = {0020-9910,1432-1297},
   MRCLASS = {20C08 (16G99 20F55)},
  MRNUMBER = {2551762},
MRREVIEWER = {Andrew\ Mathas},
       DOI = {10.1007/s00222-009-0204-8},
       URL = {https://doi.org/10.1007/s00222-009-0204-8},
}

@article{26KS,
      title={Poset modules of the $0$-{H}ecke algebras of type $B$}, 
      author={Young-Hun Kim and Dominic Searles},   
      JOURNAL = {},
      YEAR = {arXiv:2601.22926, 2026}
}

@article {03McNamara,
    AUTHOR = {McNamara, Peter},
     TITLE = {E{L}-labelings, supersolvability and 0-{H}ecke algebra actions
              on posets},
   JOURNAL = {J. Combin. Theory Ser. A},
  FJOURNAL = {Journal of Combinatorial Theory. Series A},
    VOLUME = {101},
      YEAR = {2003},
    NUMBER = {1},
     PAGES = {69--89},
      ISSN = {0097-3165,1096-0899},
   MRCLASS = {05E25 (20C08)},
  MRNUMBER = {1953281},
MRREVIEWER = {Jian-yi\ Shi},
       DOI = {10.1016/S0097-3165(02)00019-5},
       URL = {https://doi-org.libproxy.swu.ac.kr/10.1016/S0097-3165(02)00019-5},
}

\end{document}